\documentclass[%
 aip,
 amsmath,amssymb,
 reprint,%
]{revtex4-2}

\usepackage{graphicx}
\usepackage{dcolumn}
\usepackage{bm}

\usepackage[colorlinks,allcolors=blue]{hyperref} 

\usepackage[utf8]{inputenc}
\usepackage[T1]{fontenc}
\usepackage{mathptmx}

\usepackage{subcaption}
\usepackage{caption}

\usepackage{calrsfs}

\DeclareMathAlphabet{\pazocal}{OMS}{zplm}{m}{n}
\DeclareMathAlphabet\mathbfcal{OMS}{cmsy}{b}{n}
\SetMathAlphabet\pazocal{bold}{OMS}{zplm}{bx}{n}

\newcommand{\pza}{\pazocal{A}}

\newcommand{\pzc}{\pazocal{C}}

\newcommand{\pzs}{\pazocal{S}}

\newcommand{\MATHUCC}{School of Mathematical Sciences, University College Cork, Cork T12 XF62, Ireland}
\newcommand{\CART}{Collins Aerospace – Applied Research \& Technology, Cork T23 XN53, Ireland}

\begin{document}

\preprint{AIP/123-QED}

\title[Seeing double with a multifunctional reservoir computer]{Seeing double with a multifunctional reservoir computer}

\author{Andrew Flynn}
\email{andrew\_flynn@umail.ucc.ie}
\affiliation{ 
\MATHUCC
}%

\author{Vassilios A. Tsachouridis}
\affiliation{%
\CART
}%

\author{Andreas Amann}
\affiliation{ 
\MATHUCC
}%

\date{\today}

\begin{abstract}

Multifunctional biological neural networks exploit multistability in order to perform multiple tasks without changing any network properties. Enabling artificial neural networks (ANNs) to obtain certain multistabilities in order to perform several tasks, where each task is related to a particular attractor in the network's state space, naturally has many benefits from a machine learning perspective. Given the association to multistability, in this paper we explore how the relationship between different attractors influences the ability of a reservoir computer (RC), which is a dynamical system in the form of an ANN, to achieve multifunctionality. We construct the `seeing double' problem in order to systematically study how a RC reconstructs a coexistence of attractors when there is an overlap between them. As the amount of overlap increases, we discover that for multifunctionality to occur, there is a critical dependence on a suitable choice of the spectral radius for the RC's internal network connections. A bifurcation analysis reveals how multifunctionality emerges and is destroyed as the RC enters a chaotic regime that can lead to chaotic itinerancy. 

\end{abstract}

\maketitle

\begin{quotation}
Multifunctionality is an emerging area of interest in reservoir computing that opens the door to many new machine learning application areas. Multifunctionality describes the ability of a neural network to perform multiple tasks without changing any network properties. However, there is still relatively little known about how multifunctional reservoir computers (RCs) perform certain tasks, such as reconstructing coexistences of attractors. Questions regarding the relationship between the attractors and how many attractors can coexist become important. In this paper we shed some light on these issues by exploring the dynamics of a RC when trained to solve the `seeing double' problem. To be more specific, this problem involves training a RC to reconstruct a coexistence of a pair of two-dimensional circular orbits that rotate in opposite directions and can be moved closer together and overlap in state space. We find that the range of RC training parameters where multifunctionality is achieved greatly diminishes as the amount of overlap between the orbits increases. In the extreme case when the orbits are completely overlapping, we examine the nature of the basins of attraction and explore the bifurcation structure of the trained RC. By virtue of its simplicity, the seeing double problem enables us to study the fundamentals and improve the formalism behind how a RC achieves multifunctionality.

\end{quotation}

\section{\label{sec:Intro}Introduction}

A `reservoir computer' (RC) is a dynamical system that can be realised as an artificial neural network (ANN) and trained to solve certain machine learning (ML) problems. 
As outlined in \cite{nakajima2021RCbook}, \textcite{verstraeten2005firstRCmention} coined the term RC to unify the following shared philosophies of \textcite{Jaeger01ESN} and \textcite{Maass02_LSM}: instead of training all the weights in a network it is sufficient to optimise only the weights of a readout layer to solve a given task. This ideological shift in training ANNs stems from the design of a RC's internal layer, known as the `reservoir', whose role is to enable the RC to respond to a sequence of driving inputs that are related to a particular problem. Depending on the nature of the problem, before the RC is trained it is described as an open-loop (driven) system and its response to a sequence of driving inputs is used to train the RC by producing a readout layer that can, for instance, transform the open-loop RC to a closed-loop (autonomous) system or enable the open-loop RC to generate a desired response to a different input.
Using this approach, RCs have been trained to perform tasks such as, decoding neuronal information in brain-machine interfaces \cite{sussillo12_BCI}, model-free control of nonlinear systems \cite{waegeman2012feedbackcontrol,CanadayGauthier20modelfreecontrol}, predicting critical transitions \cite{lim20pred_crit_trans,kong21PredCritTrans}, and attractor reconstruction \cite{JaegerHaas04ESN,LuHuntOtt18RC}.

Within the context of attractor reconstruction, which involves training an open-loop RC to become a closed-loop RC that reconstructs the long-term dynamics of a given attractor, the following question naturally emerges, what is the `reconstruction capacity' of a given RC design? In other words, can a RC be trained to have more than one function and become multifunctional? If so, how many attractors can a given RC reconstruct simultaneously and what are the limiting factors? It is questions like these that motivated the translation of `multifunctionality' from biological neural networks (BNNs) to ANNs using a RC \cite{flynn2021multifunctionality}. Here it was shown that,
in a similar fashion to the operational principles of multifunctional neural networks in the brains of humans and other animals, 
a RC can also be trained to harness its inherent ability to exhibit multistable dynamics and achieve multifunctionality by performing multiple tasks without needing to change any network properties to do so. It is important to make clear that, in the present paper, multifunctionality is regarded as a feature of only the trained closed-loop multistable RC.

While it has been demonstrated that multifunctional RCs can possess exotic multistabilities like, for example, reconstructing coexistences of chaotic attractors from different dynamical systems \cite{flynn2021multifunctionality,flynn2021symmetry,Flynn22_LimitsMF,herteux2021MScthesis}, several issues arise that do not ordinarily occur in the case of reconstructing only a single attractor.
For instance, as shown in \textcite{Flynn22_LimitsMF}, when training different types of RCs to reconstruct a coexistence of the Lorenz and Halvorsen chaotic attractors, the RCs fail to achieve multifunctionality when the training data describing these attractors are too close together and begin to overlap, i.e., share common regions of state space. 
However, there is a limit to the level of insight about the relationship between multifunctionality and overlapping training data that can be gained from the above results given this arbitrary choice of chaotic attractors. This subsequently motivates the need for a reductionist approach to study this issue of overlap in a systematic manner. 

In this paper we introduce the `seeing double' problem as a means to explore the dynamics of a RC when trained to achieve multifunctionality in a paradigmatic case of overlapping training data. The seeing double problem involves training a RC to reconstruct a coexistence of two circular orbits that rotate in opposite directions. While these are relatively simple dynamical objects for a RC to reconstruct, we find that the closer together these orbits are the more problematic it becomes for the closed-loop RC to achieve multifunctionality.

In Sec.\,\ref{sec:MFRC} we further discuss the basics of multifunctionality and outline the procedure we use to train a RC to become multifunctional. In Sec.\,\ref{sec:SeeingDouble} we describe how the training data is generated. The results of our numerical experiments are detailed in Secs.\,\ref{sec:SD_results} and \ref{sec:SDdynamicsanal}. Here we show that as the orbits are moved closer together, the range of spectral radii of the RC's internal layer, $\rho$, where multifunctionality is achieved greatly diminishes. Despite these difficulties, we find that for a small range of $\rho$ values the closed-loop RC achieves multifunctionality when the orbits are completely overlapping. We examine the nature of the trained closed-loop RC's basins of attraction in this extreme case and unearth the existence of several `untrained attractors', attractors which exist in the trained closed-loop RC's state space but were not present during the training. We explore the interplay between these untrained attractors and the reconstructed orbits and identify the particular bifurcations responsible for the rise and fall of multifunctionality. The relationship between symmetry and the seeing double problem is discussed in Appendix\,\ref{sssec:SymmUAs}. Several other dynamical phenomena are also encountered during our analysis, such as routes to chaos and evidence of chaotic itinerancy. We provide some concluding remarks in Sec.\,\ref{sec:DisConc}.

\section{\label{sec:MFRC}Multifunctional Reservoir Computers}

\subsection{\label{sec:Multifunctionality}Introduction to Multifunctionality}

Multifunctionality is the term used to describe neural networks that are capable of changing their dynamics on
demand of a given duty without needing to alter their synaptic properties. This is a fundamental feature of many neural architectures and is considered to be key to the survival of certain species over time. Multifunctionality is typically found in BNNs with a small number of neurons which are used to perform a set of mutually exclusive tasks. For example, in switching between swimming and crawling motions \cite{briggman2006imaging} or regular breathing, sighing, and gasping \cite{lieske00ReconfigurationBreathing, Popescu02MFNN}. Multifunctionality has been an active area of research in neuroscience since the mid-1980s with seminal work published by \citet{Mpitsos86MFNN} and \citet{getting89Principles}, followed by review papers by \citet{Dickinson95MF} and \citet{Marder96MFprinciples} and more recently reviewed by \citet{briggman08multifunctional}.


At the core of the above examples is a BNN that is capable of changing its activity patterns based on a given input without altering any synapses. From the perspective of dynamical systems, we can say that a multifunctional neural network possesses a multistability. Each task that the network performs there is, in this sense, an attractor associated with it. This attractor is in coexistence with several other attractors in the network's state space, each of which are distinctly related to the tasks that the network performs.

This multistability interpretation of events dates back to the seminal work of \citet{Mpitsos86MFNN}, where the researchers use the term multifunctionality to describe networks in which neurons have multiple dynamical regimes but without altering the strength of the synapses. Through Definition\,1 in Sec.\,\ref{sec:trainingMFRC} we provide a more precise definition of multifunctionality in the context of reservoir computing.

It is natural to consider that many more network systems may have the ability to achieve multifunctionality. Where this becomes immediately relevant is in the domain of ML. By harnessing multifunctionality from a ML perspective, it can be used to unlock additional computational capabilities of certain ANNs which would otherwise have remained dormant. 

Moreover, training ANNs to achieve multifunctionality is highly advantageous from a practicality point of view as it has the potential to broaden the networks functional capacity to solve multiple tasks using the same set of trained weights. As multifunctionality brings multistability into the world of ML, this in turn opens the door for many new applications areas. 

The development of multifunctional RCs came as a result of recognising that since a multifunctional neural network in principle resembles a system with a coexistence of attractors, and that a RC can be trained to reconstruct the dynamics of an attractor, then by extending this formalism a RC was trained to exhibit multifunctionality by reconstructing a coexistence of attractors \cite{flynn2021multifunctionality}.

\subsection{\label{sec:RCintro}Reservoir computing}

The RC that is studied throughout this paper was presented by \textcite{LuHuntOtt18RC}. To outline the steps in how this RC is trained to reconstruct a coexistence of attractors and achieve multifunctionality, let's first consider the case of training this RC to reconstruct the dynamics of a single attractor, $\pza \subset \mathbb{R}^{D}$, given access to a trajectory on $\pza$ described by a vector $\boldsymbol{u}(t) \in \pza$.

There are two stages involved in training this RC, a listening stage and a training stage and we outline the specifics of these stages after introducing the RC. In both stages, the open-loop RC, defined as the nonautonomous dynamical system in Eq.\,\eqref{eq:ListenRes} (referred to as the open-loop system in Sec.\,\ref{sec:Intro}), is driven by $\boldsymbol{u}(t)$ for $0 \leq t \leq t_{\text{train}}$ and its response to this driving input is found by generating solutions of the following,
\begin{align}
    \dot{\boldsymbol{r}}(t) &= \gamma \left[ - \boldsymbol{r}(t) + \tanh{\left( \textbf{M} \, \boldsymbol{r}(t) + \sigma \textbf{W}_{in} \, \boldsymbol{u}(t) \right)} \right], \label{eq:ListenRes}\\
    \boldsymbol{r}\left( 0 \right) &= \boldsymbol{0}^{T}.\label{eq:ListenResIC}
\end{align}
In Eq.\,\eqref{eq:ListenRes}, $\boldsymbol{r}(t) \in \mathbb{R}^{N}$ describes the state of the open-loop RC at a given time $t$ and $N$ is the number of artificial neurons in the network. Solutions of Eq.\,\eqref{eq:ListenRes} are computed using the 4$^{th}$ order Runge-Kutta method with time step $\tau$. $\gamma$ is a decay-rate parameter arising from the derivation of this RC from the initial discrete-time design proposed by \textcite{Jaeger01ESN}. The $\tanh$ `activation function' is a pointwise operation and is defined as $\tanh\left( \cdot \right) : \mathbb{R}^{N} \to \mathbb{R}^{N}$. The topology of the network is described by the adjacency matrix, $\textbf{M} \in \mathbb{R}^{N \times N}$. The input strength parameter, $\sigma$, and the input matrix, $\textbf{W}_{in} \in \mathbb{R}^{N \times D}$, when multiplied together represent the weight given to $\boldsymbol{u}(t)$, as it is projected into the open-loop RC. The steps involved in constructing $\textbf{M}$ and $\textbf{W}_{in}$ are outlined in Appendix\,\ref{app:RCdesign}.

In our results, the spectral radius of $\textbf{M}$, which is denoted by $\rho$, is shown to play a crucial role in training this RC design to solve the seeing double problem and achieve multifunctionality. $\rho$ is tuned by rescaling the elements of $\textbf{M}$ such that the maximum of the absolute values of its eigenvalues is $\rho$. $\rho$ has also been a key parameter in our previous results on training a RC to achieve multifunctionality \cite{flynn2021multifunctionality,flynn2021symmetry,Flynn22_LimitsMF}.

We now describe the specifics of the listening and training stages we mentioned earlier. The listening stage is described as the duration of time where Eq.\,\eqref{eq:ListenRes} is driven by $\boldsymbol{u}(t)$ for $0 \leq t \leq t_{\text{listen}}$ where $t_{\text{listen}}<t_{\text{train}}$ is chosen such that $\boldsymbol{r}(t)$ is determined by a history of driving inputs and is no longer dependent on its initial condition.

The training stage is the time where Eq.\,\eqref{eq:ListenRes} is driven by $\boldsymbol{u}(t)$ for $t_{\text{listen}} \leq t \leq t_{\text{train}}$. In this paper, training consists of finding a `readout function/layer' defined as, $\hat{\boldsymbol{\psi}}\left( \cdot \right): \mathbb{R}^{2 N} \to \mathbb{R}^{D}$, that replaces $\boldsymbol{u}(t)$ in Eq.\,\eqref{eq:ListenRes} and is written as, 
\begin{align}
    \hat{\boldsymbol{\psi}}\left(\boldsymbol{r}(t)\right) = \textbf{W}_{out} \boldsymbol{q}( \boldsymbol{r}(t) ),\label{eq:ReadoutFunction}
\end{align}
where $\textbf{W}_{out} \in \mathbb{R}^{D \times 2N}$ is the `readout matrix' and $\boldsymbol{q}( \boldsymbol{r}(t) ) \in \mathbb{R}^{2 N}$ is given by,
\begin{align}
    \boldsymbol{q}(\boldsymbol{r}(t))=\left(\begin{array}{c}
\boldsymbol{r}(t)\\
\boldsymbol{r}^{2}(t)
\end{array}\right),\label{eq:q_square}
\end{align}
where $\boldsymbol{r}^{2}(t) = \left( r_{1}^{2}(t), r_{2}^{2}(t), \ldots, r_{N}^{2}(t) \right)^{T}$. From Eq.\,\eqref{eq:q_square} it then follows that Eq.\,\eqref{eq:ReadoutFunction} can be rewritten as,
\begin{align}
    \hat{\boldsymbol{\psi}}\left(\boldsymbol{r}(t)\right) = \textbf{W}_{out}^{(1)} \, \boldsymbol{r}(t) + \textbf{W}_{out}^{(2)} \, \boldsymbol{r}^{2}(t), \label{eq:BreakSym_r2}
\end{align}
where $\textbf{W}_{out}^{(1)}$ is the `linear readout matrix', and $\textbf{W}_{out}^{(2)}$ is the `square readout matrix' which breaks the symmetry in Eq.\,\eqref{eq:ListenRes} when replacing $\boldsymbol{u}(t)$ with $\hat{\boldsymbol{\psi}}\left(\boldsymbol{r}(t)\right)$ after the training and prevents the occurrence of `mirror-attractors' which can impede the ability of the RC to reconstruct attractors\cite{herteux2020Symm,flynn2021symmetry}.

$\textbf{W}_{out}$ is determined by the ridge regression technique which consists of minimising the following expression,
\begin{align}
    \frac{1}{t^{*}-l^{*}}\left[\sum_{i=l^{*}}^{t^{*}} || \textbf{W}_{out}\,\boldsymbol{q}(\boldsymbol{r}[i]) - \boldsymbol{u}[i]||_{2}^{2} + \beta \, || \textbf{W}_{out} ||_{2}^{2}\right],\label{eq:minimiseWout}
\end{align}
with respect to $\textbf{W}_{out}$, where $l^{*} = t_{\text{listen}}/\tau$ and $t^{*} = t_{\text{train}}/\tau$, $\boldsymbol{r}[i]$ and $\boldsymbol{u}[i]$ are discrete-time samples of the continuous-time variables $\boldsymbol{r}(t)$ and $\boldsymbol{u}(t)$ at discrete time $i = t/\tau$ where, in this paper, $\tau$ is chosen as the time step of the integration. The corresponding time series, $\left[ \boldsymbol{r}[i]\right]_{i=l^{*}}^{t^{*}}$ and $\left[ \boldsymbol{u}[i]\right]_{i=l^{*}}^{t^{*}}$, are constructed by sampling $\boldsymbol{r}(t)$ and $\boldsymbol{u}(t)$ at intervals of length $\tau$ (discrete-time intervals of length $1$). $\beta$ is the regularisation parameter and the purpose of the $\beta \, || \textbf{W}_{out} ||_{2}$ term in Eq.\,\eqref{eq:minimiseWout} is to modify the linear least-squares regression to reduce the magnitudes of elements in $\textbf{W}_{out}$ in order to discourage overfitting.

Minimising Eq.\,\eqref{eq:minimiseWout} involves computing its partial derivative with respect to $\textbf{W}_{out}$ and setting the resulting expression equal to $0$. The $\textbf{W}_{out}$ that minimises Eq.\,\eqref{eq:minimiseWout} is given by,
\begin{align}
    \textbf{W}_{out} = \textbf{Y} \textbf{X}^{T} \left( \textbf{X} \textbf{X}^{T} + \beta \, \textbf{I} \right)^{-1},\label{eq:WoutRegression}
\end{align}
where,
\begin{align}
    \textbf{X} = \left[ \begin{array}{cccc}
    \boldsymbol{q}(\boldsymbol{r}[l^{*}]) & \boldsymbol{q}(\boldsymbol{r}[l^{*}+1]) 
    &
    \cdots 
    &
    \boldsymbol{q}(\boldsymbol{r}[t^{*}])
    \end{array} \right],\label{eq:Xmat}
\end{align}
is the `response data matrix' and,
\begin{align}
    \textbf{Y} = \left[ \begin{array}{cccc}
        \boldsymbol{u}[l^{*}] & \boldsymbol{u}[l^{*}+1] & \cdots & \boldsymbol{u}[t^{*}]
    \end{array} \right],\label{eq:Ymat}
\end{align}
is the `input data matrix' and $\textbf{I}$ is the identity matrix.

After the training, $\boldsymbol{u}(t)$ in Eq.\,\eqref{eq:ListenRes} is replaced by $\hat{\boldsymbol{\psi}} \left( \boldsymbol{r}(t) \right)$ and in Eq.\,\eqref{eq:PredRes} we now define the closed-loop RC (referred to as the trained closed-loop system in Sec.\,\ref{sec:Intro}) as the following autonomous dynamical system,
\begin{align}
    \hspace{-0.1cm}\dot{\hat{\boldsymbol{r}}}(t) &= \gamma \left[ - \hat{\boldsymbol{r}}(t) + \tanh{\left( \textbf{M} \, \hat{\boldsymbol{r}}(t) + \sigma \textbf{W}_{in}  \textbf{W}_{out} \, \boldsymbol{
    q}(\hat{\boldsymbol{r}}(t)) \right)} \right], \label{eq:PredRes}\\
    \hspace{-0.1cm}\hat{\boldsymbol{r}}\left( 0 \right) &= \boldsymbol{r}\left( t_{train} \right),\label{eq:PredResIC}
\end{align}
where $\hat{\boldsymbol{r}}(t)$ denotes the state of the closed-loop RC at a given time $t$. While $\hat{\boldsymbol{r}}(t)$ and $\boldsymbol{r}(t)$ are both $N$-dimensional vectors, to distinguish the dynamics of $\hat{\boldsymbol{r}}(t)$ from $\boldsymbol{r}(t)$, $\hat{\boldsymbol{r}}(t)$ is defined as $\hat{\boldsymbol{r}}(t) \in \mathbb{S}$ where $\mathbb{S}$ is referred to as the `RC's state space' and is used henceforth when discussing the dynamics of the closed-loop RC. By computing solutions of Eq.\,\eqref{eq:PredRes}, predictions of $\boldsymbol{u}(t)$ for $t>t_{\text{train}}$, denoted as $\hat{\boldsymbol{u}}(t)$, are given by,
\begin{align}
    \hat{\boldsymbol{u}}(t) = \hat{\boldsymbol{\psi}}\left( \hat{\boldsymbol{r}}(t)\right).\label{eq:uhat}
\end{align}
Again, while both $\boldsymbol{u}(t)$ and $\hat{\boldsymbol{u}}(t)$ are $D$-dimensional vectors, to distinguish the dynamics of $\hat{\boldsymbol{u}}(t)$ from $\boldsymbol{u}(t)$, $\hat{\boldsymbol{u}}(t)$ is defined as $\hat{\boldsymbol{u}}(t) \in \mathbb{P}$ where $\mathbb{P}$ is referred to as the `projected state space' and is used henceforth when discussing the dynamics of the closed-loop RC when projected via $\hat{\boldsymbol{\psi}}$.

We say the closed-loop RC has achieved attractor reconstruction when the long-term dynamical characteristics of $\hat{\boldsymbol{u}}(t)$ are indistinguishable from $\boldsymbol{u}(t)$. In this scenario there exists an attractor $\pzs \subset \mathbb{S}$ such that when the state of the closed-loop RC approaches $\pzs$ and is projected from $\mathbb{S}$ to $\mathbb{P}$ via $\hat{\boldsymbol{\psi}}\left( \cdot \right)$, the dynamics of the `reconstructed attractor', $\hat{\pza} \subset \mathbb{P}$, will resemble the dynamics of $\pza$ in the long-term. By resembling the long-term dynamics it is meant that, for instance, $\pza$ and $\hat{\pza}$ will have nearly identical Poincar{\'e} sections when computed for the same region of $\mathbb{R}^{D}$ and $\mathbb{P}$ as $t \to \infty$ for $t > t_{\text{train}}$.

\subsection{\label{sec:trainingMFRC}Training a RC to become multifunctional}

Using the terminology we have established, we now define multifunctionality in the context of having successfully trained the closed-loop RC in Eq.\,\eqref{eq:PredRes} to reconstruct a coexistence of $n$ attractors, $\pza_{1}\subset\mathbb{R}^{D}$, $\pza_{2}\subset\mathbb{R}^{D}$, $\ldots$, $\pza_{n}\subset\mathbb{R}^{D}$, given input training data from each of, $\boldsymbol{u}_{\left(\pza_{1}\right)}(t) \in \pza_{1}$, $\boldsymbol{u}_{\left(\pza_{2}\right)}(t) \in \pza_{2}$, $\ldots$, $\boldsymbol{u}_{\left(\pza_{n}\right)}(t) \in \pza_{n}$, for $0 < t \leq t_{\text{train}}$.

\textbf{Definition\,1:}\,\textit{ The closed-loop RC, defined in Eq.\,\eqref{eq:PredRes}, achieves multifunctionality if for every $\pza_{j}$, for $j=1,\, \ldots,\, n$, there exists a corresponding attractor $\pzs_{j} \subset \mathbb{S}$ such that the projection of each $\pzs_{j}$ from $\mathbb{S}$ to $\mathbb{P}$ via $\hat{\boldsymbol{\psi}}$ resembles the long-term dynamics of the respective $\pza_{j}$.}

The projection of a given $\pzs_{j}$ from $\mathbb{S}$ to $\mathbb{P}$ is also referred to as the reconstructed attractor $\hat{\pza}_{j} \subset \mathbb{P}$.
What Definition\,1 says it that in order for a RC to achieve multifunctionality it is necessary for it to inherit the ability to \textbf{perform multiple tasks using the same readout matrix}, $\textbf{W}_{out}$, \textbf{and without changing any other structural properties of the RC}.


The steps involved in computing $\textbf{W}_{out}$ from Sec.\,\ref{sec:RCintro} are adapted for the case of multifunctionality. For simplicity, let's first consider the case of reconstructing a coexistence of two attractors, $\pza_{1}\subset\mathbb{R}^{D}$ and $\pza_{2}\subset\mathbb{R}^{D}$. 

We first drive the open-loop RC in Eq.\,\eqref{eq:ListenRes} with input $\boldsymbol{u}_{\left(\pza_{1}\right)}(t) \in \pza_{1}$ for $0 < t \leq t_{\text{train}}$ and then repeat for $\boldsymbol{u}_{\left(\pza_{2}\right)}(t) \in \pza_{2}$. The responses of the open-loop RC to these driving inputs are denoted by $\boldsymbol{r}_{\left(\pza_{1}\right)}(t)$ and $\boldsymbol{r}_{\left(\pza_{2}\right)}(t)$.
It is important to highlight that $\textbf{M}$, $\textbf{W}_{in}$, and all training parameters remain identical when generating $\boldsymbol{r}_{\left(\pza_{1}\right)}(t)$ and $\boldsymbol{r}_{\left(\pza_{2}\right)}(t)$ in order to be consistent with the descriptions of multifunctionality from the biological perspective and the reservoir computing analogue in Sec.\,\ref{sec:Multifunctionality}.
The same readout function design as in Eq.\,\eqref{eq:ReadoutFunction} is used and the ridge regression approach for computing $\textbf{W}_{out} \in \mathbb{R}^{D \times 2N}$ in Eq.\,\eqref{eq:minimiseWout} is modified to account for additional attractors. This modified ridge regression approach consists of minimising the following expression with respect to $\textbf{W}_{out}$,
\begin{align}
    &\frac{1}{t^{*}-l^{*}}\Bigg[\sum_{i=l^{*}}^{t^{*}} || \textbf{W}_{out}\,\boldsymbol{q}(\boldsymbol{r}_{\left(\pza_{1}\right)}[i]) - \boldsymbol{u}_{\left(\pza_{1}\right)}[i]||_{2}^{2} \nonumber \\
    & + \sum_{i=l^{*}}^{t^{*}} || \textbf{W}_{out}\,\boldsymbol{q}(\boldsymbol{r}_{\left(\pza_{2}\right)}[i]) - \boldsymbol{u}_{\left(\pza_{2}\right)}[i]||_{2}^{2} + \beta \, || \textbf{W}_{out} ||_{2}^{2}\Bigg],\label{eq:minimiseWoutMF}
\end{align}
where $\boldsymbol{r}_{\left(\pza_{1}\right)}[i]$, $\boldsymbol{u}_{\left(\pza_{1}\right)}[i]$, $\boldsymbol{r}_{\left(\pza_{2}\right)}[i]$, and $\boldsymbol{u}_{\left(\pza_{2}\right)}[i]$ are discrete-time samples of $\boldsymbol{r}_{\left(\pza_{1}\right)}(t)$, $\boldsymbol{u}_{\left(\pza_{1}\right)}(t)$, $\boldsymbol{r}_{\left(\pza_{2}\right)}(t)$, and $\boldsymbol{u}_{\left(\pza_{2}\right)}(t)$ and are constructed using the same convention outlined in Sec.\,\ref{sec:RCintro}. The regularisation parameter $\beta$ plays the same role as discussed in Sec.\,\ref{sec:RCintro}. The $\textbf{W}_{out}$ that minimises Eq.\,\eqref{eq:minimiseWoutMF} is,
\begin{align}
    \textbf{W}_{out} = \textbf{Y}_{C} \textbf{X}_{C}^{T} \left( \textbf{X}_{C} \textbf{X}_{C}^{T} + \beta \, \textbf{I} \right)^{-1},\label{eq:WoutRegressionMF}
\end{align}
where $\textbf{X}_{C} = \left[ \textbf{X}_{\left(\pza_{1}\right)}, \, \textbf{X}_{\left(\pza_{2}\right)} \right]$, $\textbf{Y}_{C} = \left[  \textbf{Y}_{\left(\pza_{1}\right)}, \, \textbf{Y}_{\left(\pza_{2}\right)} \right]$ with $\textbf{X}_{\left(\pza_{1}\right)}$ and $\textbf{X}_{\left(\pza_{2}\right)}$ constructed as in Eq.\,\eqref{eq:Xmat} for the corresponding $\boldsymbol{r}_{\left(\pza_{1}\right)}[i]$ and $\boldsymbol{r}_{\left(\pza_{2}\right)}[i]$ and similarly for $\textbf{Y}_{\left(\pza_{1}\right)}$ as in Eq.\,\eqref{eq:Ymat} and $\textbf{Y}_{\left(\pza_{2}\right)}$ for the associated $\boldsymbol{u}_{\left(\pza_{1}\right)}[i]$ and $\boldsymbol{u}_{\left(\pza_{2}\right)}[i]$. In Eq.\,\eqref{eq:WoutRegressionMF}, $\textbf{I}$ is the identity matrix of the appropriate size.

The closed-loop RC is then described by the same autonomous continuous-time dynamical system as in Eq.\,\eqref{eq:PredRes} using the $\textbf{W}_{out}$ that is obtained from Eq.\,\eqref{eq:WoutRegressionMF}. If multifunctionality is achieved then we refer to the resulting multistable closed-loop RC as the `multifunctional RC'. 

As per Definition 1, for each $\pza_{1}$ and $\pza_{2}$ there exists two coexisting attractors $\pzs_{1}, \pzs_{2} \subset \mathbb{S}$ such that when the state of the multifunctional RC approaches either $\pzs_{1}$ or $\pzs_{2}$ and is projected from $\mathbb{S}$ to $\mathbb{P}$ using $\hat{\boldsymbol{\psi}}\left(\cdot\right)$ as in Eq.\,\eqref{eq:ReadoutFunction}, then the long-term dynamics of the reconstructed attractors, $\hat{\pza}_{1}$ and $\hat{\pza}_{2}$, will resemble the long-term dynamics of $\pza_{1}$ and $\pza_{2}$. 

Therefore, to reconstruct the dynamics of either $\pza_{1}$ or $\pza_{2}$, the multifunctional RC needs to be initialised with the corresponding $\hat{\boldsymbol{r}}\left( 0 \right) = \boldsymbol{r}_{\left(\pza_{1}\right)}\left( t_{train} \right)$ or $\boldsymbol{r}_{\left(\pza_{2}\right)}\left( t_{train} \right)$, or some other initial condition that is in the basin of attraction of the corresponding $\pzs_{1}$ or $\pzs_{2}$.


The same steps as outlined in the present Sec.\,\ref{sec:trainingMFRC} can be extended in order to produce a multifunctional RC that reconstructs a coexistence of $n$ attractors, $\pza_{1}, \, \pza_{2}, \, \ldots, \, \pza_{n}$. All that is required is to produce the corresponding $\textbf{X}_{C} = \left[ \textbf{X}_{\left(\pza_{1}\right)}, \, \textbf{X}_{\left(\pza_{2}\right)}, \, \ldots, \, \textbf{X}_{\left(\pza_{n}\right)} \right]$ and $\textbf{Y}_{C} = \left[  \textbf{Y}_{\left(\pza_{1}\right)}, \, \textbf{Y}_{\left(\pza_{2}\right)}, \, \ldots, \, \textbf{Y}_{\left(\pza_{n}\right)} \right]$ and solve for $\textbf{W}_{out}$ as in Eq.\,\eqref{eq:WoutRegressionMF}.

\subsection{Summary of previous work}\label{ssec:PreviousWork}

As mentioned previously, multifunctionality was first translated from BNNs to ANNs in \textcite{flynn2021multifunctionality}. Here it was shown that the same RC design as in Eq.\,\eqref{eq:ListenRes} can be trained to reconstruct a coexistence of chaotic attractors from a multistable system, different parameter settings of the same system, and two systems with different nonlinearities. The issues of choosing appropriate training parameters, differences in the timescales of the chaotic attractor's dynamics, and the existence of untrained attractors are also investigated. Independent of these results, \textcite{herteux2020Symm} also demonstrated how to train a RC to reconstruct a coexistence of chaotic attractors using the same training technique described in Sec.\,\ref{sec:trainingMFRC}. The phenomena which arise when training Eq.\,\eqref{eq:ListenRes} to reconstruct a coexistence of attractors that are related by certain symmetries are explored in \textcite{flynn2021symmetry}, the analytical and numerical results presented here show how a RC can adapt to and take on the symmetry of the trained system spontaneously. In \textcite{Flynn22_LimitsMF} the ability of different types of RCs to achieve multifunctionality in the presence of overlapping training data are assessed. The present paper extends on the above results in many ways. Moreover, it is thanks to the seeing double problem's simplicity that we are able to more rigorously examine the dynamics and improve the formalism behind how a RC achieves multifunctionality.

Since the first results regarding multifunctionality in a RC\cite{flynn2021multifunctionality}, it has been discovered that a given $\textbf{W}_{out}$ can enable a RC to perform more than one task using different principles to those involving multifunctionality. In this context, the so-called `parameter-aware' RCs have been trained to interpolate and extrapolate the dynamics of a system such as those studied in \textcite{kim2021GlobalLocal}, \textcite{kong21FWRC,kong2023digitaltwins}, and \textcite{fan2021anticipatingsynch}. To do this, an additional parameter input channel is incorporated into the RC's architecture during training so that, in this case, the closed-loop RC performs a task associated with a given parameter value and $\textbf{W}_{out}$ remains fixed after the training. In these examples of parameter-aware RCs, it is shown that the RC can go through similar bifurcations to the attractor from the dynamical system used to generate the training data without necessarily needing to be trained on data after the bifurcation takes place. In future work, it would be interesting to investigate the dynamical principles behind how this is achieved and how this relates to the phenomenon of untrained attractors as presented in \textcite{flynn2021multifunctionality}. These attractors form a key aspect of our analysis in Sec.\,\ref{sec:SDdynamicsanal}. Furthermore, while multifunctionality has so far not been studied using these parameter-aware RCs, we comment that a `multifunctional parameter-aware' RC may further enhance the reconstruction capacity of a given RC. 

\section{\label{sec:SeeingDouble}Seeing Double}

The `seeing double' problem is introduced in this section. This numerical experiment is designed as a means to systematically induce the issues related to multifunctionality and overlapping training data.

This task involves training a RC to reconstruct a coexistence of attractors that follow trajectories on two circular orbits of equal radius, rotate in opposite directions, and the centres of these orbits can be moved either closer together or further apart. When these orbits are overlapping, the RC is therefore required to distinguish between points in $\mathbb{R}^{D}$ that are common to both cycles in order to exhibit multifunctionality. 

By virtue of the seeing double problem's simplicity, a greater emphasis can be placed on examining how multifunctionality arises in a RC. At the same time, despite this simplicity, an enormity of elaborate dynamics are encountered and studied in Secs.\,\ref{sec:SD_results} and \ref{sec:SDdynamicsanal}.

\subsection[Numerical experiment setup]{Numerical experiment setup}\label{sec:SD_Experiment}

The training data is computed via,
\begin{equation}
\boldsymbol{u}(t)=
    \left( \begin{array}{c}
        x(t) \\
        y(t)
    \end{array} \right)
    = \left( \begin{array}{c}
        b_{x} \cos{\left( t \right)} + x_{cen}\\
        b_{y} \sin{\left( t \right)} + y_{cen}
    \end{array} \right),
    \label{eq:InputSys}
\end{equation}
for $t=0, \tau, 2 \tau, \ldots,$ using the time-step $\tau = 0.01$. The resultant $\boldsymbol{u}(t)$ generated from Eq.\,\eqref{eq:InputSys} resembles a trajectory around a circle of radius $b=|b_{x}|=|b_{y}|$ and centered at $\left( x_{cen}, \, y_{cen} \right)$. The sets, $\pzc_{A}$ and $\pzc_{B}$, are produced using Eq.\,\eqref{eq:InputSys} for specific values of $b_{x}, b_{y}, x_{cen},$ and $y_{cen}$. The corresponding input signals that are used to drive the open-loop RC in Eq.\,\eqref{eq:ListenRes} are denoted by $\boldsymbol{u}_{\left(\pzc_{A}\right)}(t)$ and $\boldsymbol{u}_{\left(\pzc_{B}\right)}(t)$. The open-loop RC's response to $\boldsymbol{u}_{\left(\pzc_{A}\right)}(t)$ and $\boldsymbol{u}_{\left(\pzc_{B}\right)}(t)$ are denoted as $\boldsymbol{r}_{\left(\pzc_{A}\right)}(t)$ and $\boldsymbol{r}_{\left(\pzc_{B}\right)}(t)$ for $0 \leq t \leq t_{train}$. The training procedure outlined in Sec.\,\ref{sec:trainingMFRC} is used to produce the corresponding $\textbf{X}_{\left(\pzc_{A}\right)}, \, \textbf{X}_{\left(\pzc_{B}\right)}, \, \textbf{Y}_{\left(\pzc_{A}\right)},$ and $\textbf{Y}_{\left(\pzc_{B}\right)}$ that are used to compute $\textbf{W}_{out}$ in Eq.\,\eqref{eq:WoutRegressionMF}. This $\textbf{W}_{out}$ is then used in the closed-loop RC setup given by Eq.\,\eqref{eq:PredRes}.

We say that this closed-loop RC achieve multifunctionality and solves the seeing double problem once it reconstructs a coexistence of two attractors, $\pzs_{A}$ and $\pzs_{B}$, that exist in $\mathbb{S}$ and resemble $\pzc_{A}$ and $\pzc_{B}$ when projected to $\mathbb{P}$ using $\hat{\boldsymbol{\psi}}\left(\cdot\right)$ as in Eq.\,\eqref{eq:ReadoutFunction} with the $\textbf{W}_{out}$ as mentioned above. As per the same convention used earlier, these projected dynamics of $\pzs_{A}$ and $\pzs_{B}$ are referred to as the reconstructed attractors and denoted by $\hat{\pzc}_{A}$ and $\hat{\pzc}_{B}$. To reconstruct the dynamics of $\pzc_{A}$ or $\pzc_{B}$ using this multifunctional RC we initialise Eq.\,\eqref{eq:PredRes} with $\hat{\boldsymbol{r}}(0)=\hat{\boldsymbol{r}}_{\left(\pzc_{A}\right)}(0)=\boldsymbol{r}_{\left(\pzc_{A}\right)}(t_{\text{listen}})$ or $\hat{\boldsymbol{r}}(0)=\hat{\boldsymbol{r}}_{\left(\pzc_{B}\right)}(0)=\boldsymbol{r}_{\left(\pzc_{B}\right)}(t_{\text{listen}})$ or some known point in the basin of attraction of $\pzs_{A}$ or $\pzs_{B}$. The subsequent states of Eq.\,\eqref{eq:PredRes} when approaching $\pzs_{A}, \pzs_{B} \subset \mathbb{S}$ ($\hat{\pzc}_{A}, \hat{\pzc}_{B} \subset \mathbb{P}$) are written as $\hat{\boldsymbol{r}}_{\left(\pzc_{A}\right)}(t)$ and $\hat{\boldsymbol{r}}_{\left(\pzc_{B}\right)}(t)$.

In this paper, $b_{x}$ and $b_{y}$ are set to $b_{x}=b_{y}=5$ to create $\pzc_{A}$ and for $\pzc_{B}$, $b_{x}=-5$ and $b_{y}=5$. The centres of $\pzc_{A}$ and $\pzc_{B}$ are brought closer together by changing the values of $\left( x_{cen}, \, y_{cen} \right)$. In order to simplify the experiment further, $\pzc_{A}$ and $\pzc_{B}$ are designed to be centered at $\left( x_{cen}, \, 0 \right)$ and $\left( -x_{cen}, \, 0 \right)$ respectively. Therefore, by changing $x_{cen}$ the centres of these cycles are moved equidistantly along the line $y=0$. In this case, an overlapping region between $\pzc_{A}$ and $\pzc_{B}$ exists whenever $|x_{cen}| \leq b = 5$, i.e., $\pzc_{A} \cap \pzc_{B} \neq \emptyset$ $\forall \, |x_{cen}| \leq 5$. Furthermore, $\pzc_{A}$ and $\pzc_{B}$ are said to be `entirely/completely overlapping' when $x_{cen}=0$. In this extreme case,
the only difference between $\pzc_{A}$ and $\pzc_{B}$ is the direction of rotation on both cycles.

\subsection{Illustrating issues of overlap}\label{sssec:SD_illus_overlap}

\begin{figure*}
    \centering
    \begin{subfigure}{0.49\textwidth}
        \centering
        \includegraphics[width=0.98\textwidth]{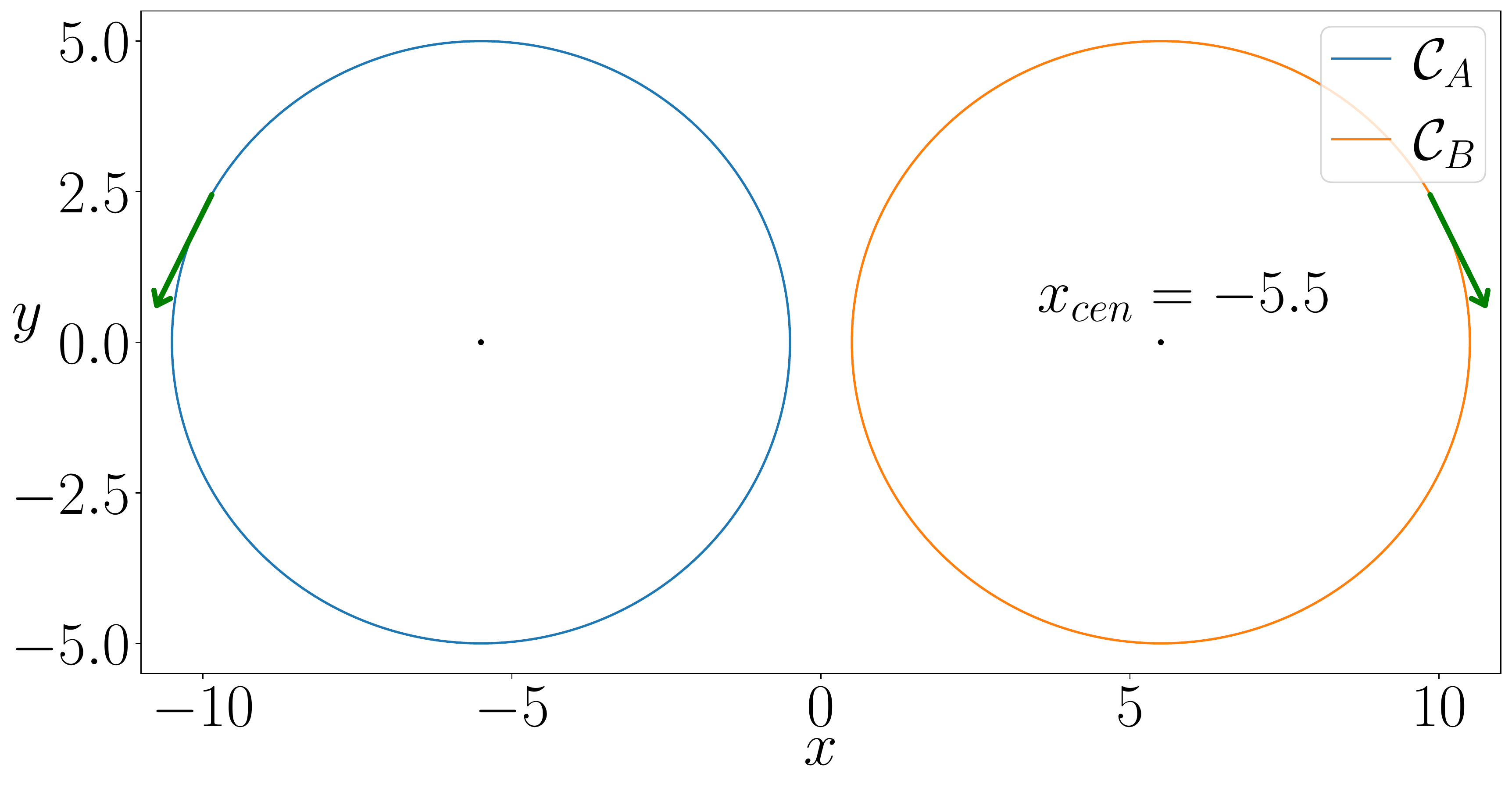}
        \caption{$\pzc_{A} \cap \pzc_{B} = \emptyset$: no overlapping}
        \label{fig:Case2_Disjoint_cycles}
    \end{subfigure}
    \hfill
    \begin{subfigure}{0.49\textwidth}
        \centering
        \includegraphics[width=0.98\textwidth]{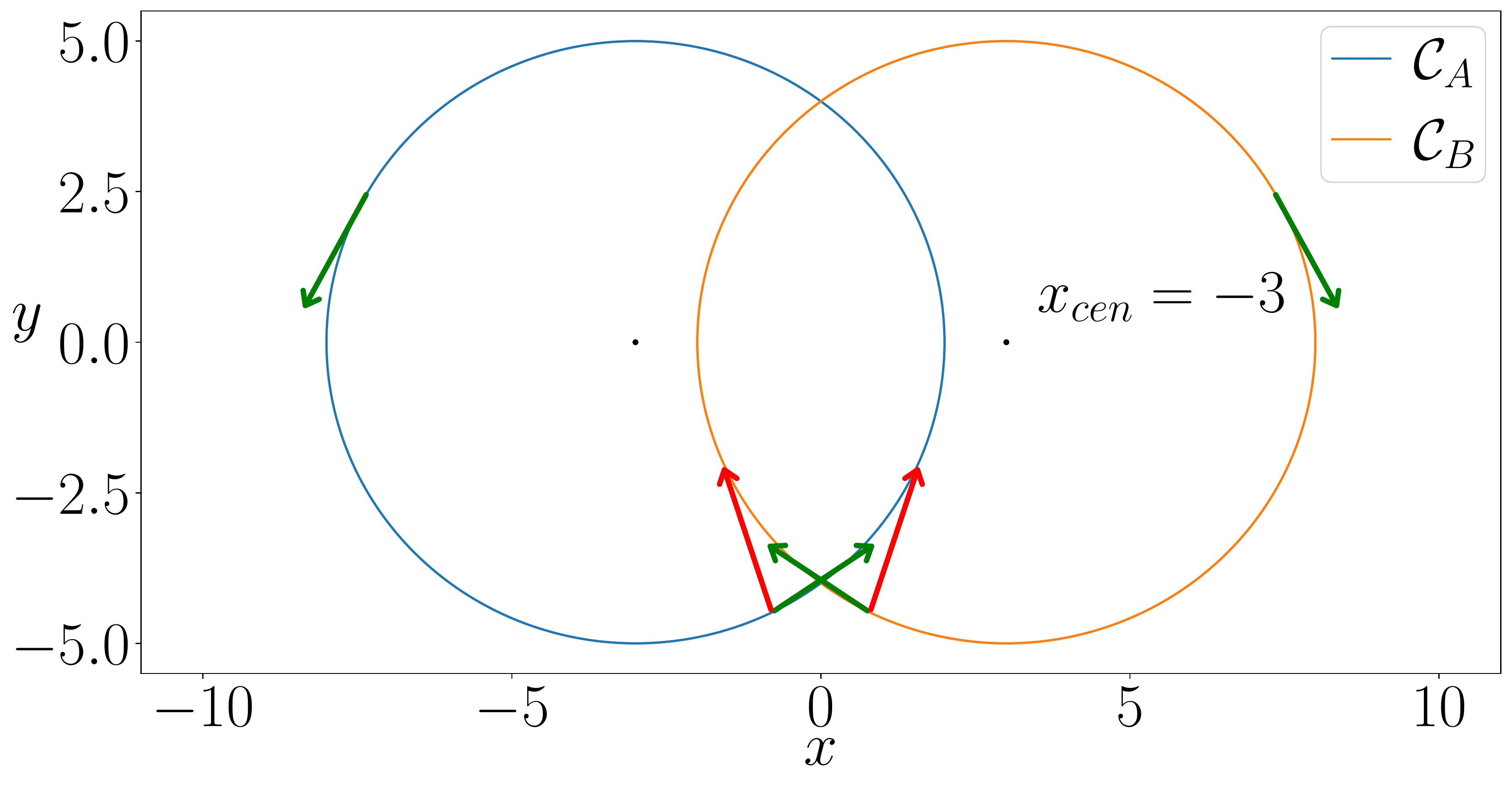}
        \caption{$\pzc_{A} \cap \pzc_{B} \neq \emptyset$: overlapping}
        \label{fig:Case2_Intersect_cycles}
    \end{subfigure}
    \caption{Illustrating how the centres of $\pzc_{A}$ and $\pzc_{B}$ move in $\mathbb{R}^{D}$ for changes in $x_{cen}$. In (a) $\pzc_{A} \cap \pzc_{B} = \emptyset$ for $x_{cen} = -5.5$, and in (b) $\pzc_{A} \cap \pzc_{B} \neq \emptyset$ when $x_{cen} = -3$. Green arrows indicate a trajectory on these cycles. Red arrows indicate a potential path to failure.}
    \label{fig:Case2_MovingCycles}
\end{figure*}

The orbits $\pzc_{A}$ and $\pzc_{B}$ produced from Eq.\,\eqref{eq:InputSys} are depicted in Fig.\,\ref{fig:Case2_MovingCycles} with $\boldsymbol{u}_{\left(\pzc_{A}\right)}$ and $\boldsymbol{u}_{\left(\pzc_{B}\right)}$ plotted in the $(x,y)$-plane for a given $x_{cen}$. The green arrows indicate the direction of rotation in which an orbit along these circles evolves upon in time.

Fig.\,\ref{fig:Case2_Disjoint_cycles} shows that the cycles do not overlap when $x_{cen}=-5.5$. However, when $|x_{cen}| \leq 5$, like in Fig.\,\ref{fig:Case2_Intersect_cycles} for $x_{cen}=-3$, the orbits share similar regions of the $\left( x, y \right)$-space. This presents a challenge for the resulting closed-loop RC to remain on the correct path when its state approaches the point where these cycles intersect as indicated by the green arrows at the crossing point in Fig.\,\ref{fig:Case2_Intersect_cycles}. If the RC fails to distinguish between these overlapping cycles then, for example, the state of the RC may move onto a different path (as indicated by the red arrow in Fig.\,\ref{fig:Case2_Intersect_cycles}) or potentially result in an attractor merging crisis scenario. 

Therefore, when there is an overlap between $\pzc_{A}$ and $\pzc_{B}$, the state of the open-loop RC must contain sufficient knowledge of its previous states during the training so that after the training the closed-loop RC remains on the correct path when approaching a crossing between $\pzc_{A}$ and $\pzc_{B}$. In Sec.\,\ref{sec:SD_results} it is found that by changing the spectral radius, $\rho$, of the RC's internal connections, which as a result tunes the weight given to the RC's previous states, it is possible to overcome the issue of overlap. 


\section{Seeing Double with a Multifunctional RC}\label{sec:SD_results}

The focus of this section is to explore the dynamics of the closed-loop RC in Eq.\,\eqref{eq:PredRes} after it is has trained to solve the seeing double problem. Here we show that as the difficulty of the seeing double problem is varied according to changes in $x_{cen}$ (the amount of overlap between $\pzc_{A}$ and $\pzc_{B}$ and the nearer they are to one another), the more `memory' is required of the RC and consequently the spectral radius, $\rho$, plays a more significant role. In this instance we do not use the word memory in a quantitative sense but in how it relates to a property of this RC design that, during the training, by increasing $\rho$ up to a certain point the state of the RC depends not only on the current input but also previous inputs. In other words, the RC has a greater ability to recall information from the past.

The results in this section are generated by training the RC using the method outlined in Sec.\,\ref{sec:trainingMFRC} using the input data, $\boldsymbol{u}_{\left(\pzc_{A}\right)}(t)$ and $\boldsymbol{u}_{\left(\pzc_{B}\right)}(t)$ constructed as described in Sec.\,\ref{sec:SD_Experiment}. The RC training and design parameters are specified in Table\,\ref{tab:RC_SDparams} and the same $\textbf{M}$ and $\textbf{W}_{in}$ matrices are used to produce the results presented throughout this paper.

While the results that follow are computed using just one random realisation of $\textbf{M}$ and $\textbf{W}_{in}$, we do not claim that our results are universal across all possible $\textbf{M}$ and $\textbf{W}_{in}$. At the same time the analysis tools we use throughout this paper provide significant insight towards how a RC solves the seeing double problem and are applicable to all random realisations of $\textbf{M}$ and $\textbf{W}_{in}$. Where appropriate, we highlight certain features of our results that appear across several random realisations of $\textbf{M}$ and $\textbf{W}_{in}$ and other features that may be heavily dependent on the particular $\textbf{M}$ and $\textbf{W}_{in}$ that is used to generate our results.

\subsection{Exploring Multifunctionality in \texorpdfstring{$\left( x_{cen}, \, \rho \right)$-plane}{TEXT}}\label{sec:SD_RegionsMF}

We now explore the long-term dynamics of the closed-loop RC (described in Eq.\,\eqref{eq:PredRes}) when initialised with $\hat{\boldsymbol{r}}(0)=\boldsymbol{r}_{\left(\pzc_{A}\right)}(t_{\text{listen}})$ and $\hat{\boldsymbol{r}}(0)=\boldsymbol{r}_{\left(\pzc_{B}\right)}(t_{\text{listen}})$, as viewed from $\mathbb{P}$, for a given $x_{cen} \in \left[ -10, 10 \right]$ and $\rho \in \left[ 0.1, 2.5 \right]$.

\subsubsection{Prediction Analysis}\label{ssec:PredAnal}

After the training a number of outcomes are possible. If the training is successful then $\hat{\pzc}_{A} \approx \pzc_{A}$ and $\hat{\pzc}_{B} \approx \pzc_{B}$ and the closed-loop RC's predicted trajectory on either $\pzc_{A}$ or $\pzc_{B}$ 
will remain on $\hat{\pzc}_{A}$ or $\hat{\pzc}_{B}$ in the long-term ($t \to \infty$). However, given the nature of the seeing double problem, it is also possible for the predicted trajectory to switch from $\hat{\pzc}_{A}$ to $\hat{\pzc}_{B}$ or vice-versa if the RC fails to lift the corresponding driving input data, $\boldsymbol{u}_{\left(\pzc_{A}\right)}(t)$ and $\boldsymbol{u}_{\left(\pzc_{B}\right)}(t)$, into separate regions of $\mathbb{S}$. Furthermore, if the RC fails to reconstruct $\pzc_{A}$ or $\pzc_{B}$ then the dynamics of Eq.\,\eqref{eq:PredRes} can, for example, decay to a fixed point, some other limit cycle, or never settle down to periodic behaviour and display evidence of the RC entering a chaotic regime.

An error metric known as the `roundness' is used to assess the accuracy of the RC's reconstruction of a given orbit . The roundness of a given cycle, $\pzc$, is denoted by $\delta \left( \pzc \right)$ and is found by calculating the difference between the radii of the maximum and minimum circles needed to enclose and inscribe the cycle. For simplicity, the centres of these maximum and minimum circles are taken to be located at the centre of the given target orbit. It was found empirically that if the relative roundness of a given cycle, $\delta_{rel} \left( \pzc \right) = \delta \left( \pzc \right) / b$, was $< \delta_{rel}^{+} = 0.25$, then the RC is in general able to successfully reconstruct the desired orbit with a high level of accuracy. This error metric is typically used in the production of circular shaped materials, like plastic tubes or copper pipe, and the formula to calculate the roundness is found in most machining manuals in line with the International Organisation for Standardisation (ISO).

The RC's reconstruction of a given orbit 
is characterised as the particular attractor that the state of the closed-loop RC is approaching for $t_{\text{predict}} - t^{*} \leq t \leq t_{\text{predict}}$ where $t_{\text{predict}}=600$ and $t^{*}=40$ as viewed from $\mathbb{P}$. The result of this is displayed in Figs.\,\ref{fig:Lerr_oppdir}-\ref{fig:Rerr_oppdir} where each point in the $\left( x_{cen}, \, \rho \right)$-plane is coloured according to the specifications in Table\,\ref{tab:ColorScheme}

\begin{table}
    \centering
    \begin{tabular}{l|l}
    \hline \hline
        Blue & Correct cycle is reconstructed and $\delta_{rel} \left( \pzc \right) < \delta_{rel}^{+}$. \\
        \hline
        Yellow & Prediction switches to the other cycle. \\
        \hline
        Magenta & Prediction does not settle down to periodic motion\\ & and label this as some attractor, $\pzs$. \\
        \hline
        Black & Prediction tends to some other limit cycle \\ & or $\delta_{rel} \left( \pzc \right) \geq \delta_{rel}^{+}$, this attractor is labelled as, LC.\\
        \hline
        Green & Prediction decays toward some fixed point labelled as FP.\\
        \hline \hline
    \end{tabular}
    \caption{Colour scheme for prediction analysis in Fig.\,\ref{fig:LR_oppdir}}
    \label{tab:ColorScheme}
\end{table}

In order to deduce the regions of the $\left( x_{cen}, \, \rho \right)$-plane in which multifunctionality was achieved, an additional picture is provided in Fig.\,\ref{fig:SD_MF_regions} where the common blue regions of Figs.\,\ref{fig:Lerr_oppdir}-\ref{fig:Rerr_oppdir} are identified and the larger of the $\delta_{rel} \left( \pzc \right)$ values when assessing the accuracy of the closed-loop RC's reconstruction of $\pzc_{A}$ and $\pzc_{B}$, defined as $\delta_{rel}^{max} = \max \left( \delta_{rel}\left( \hat{\pzc}_{A} \right), \delta_{rel}\left( \hat{\pzc}_{B} \right) \right)$, is plotted for each point in the $\left( x_{cen}, \, \rho \right)$-plane. 

\subsubsection{Regions of multifunctionality}

\begin{figure*}
    \centering
    \begin{subfigure}{0.48\textwidth}
        \centering
        \includegraphics[width=0.98\textwidth]{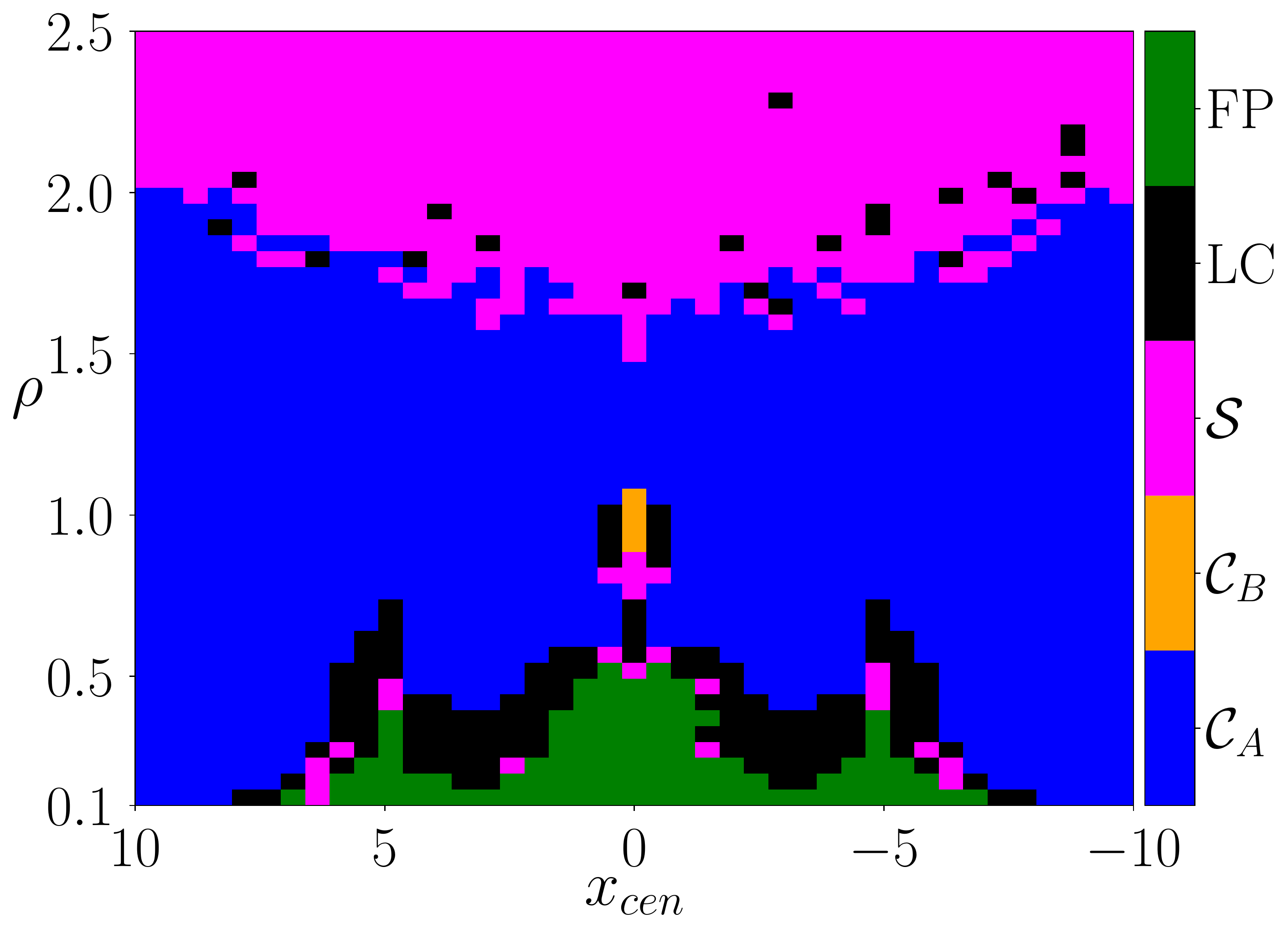}
        \caption{$\pzc_{A}$ prediction}
        \label{fig:Lerr_oppdir}
    \end{subfigure}
    \begin{subfigure}{0.48\textwidth}
        \centering
        \includegraphics[width=0.98\textwidth]{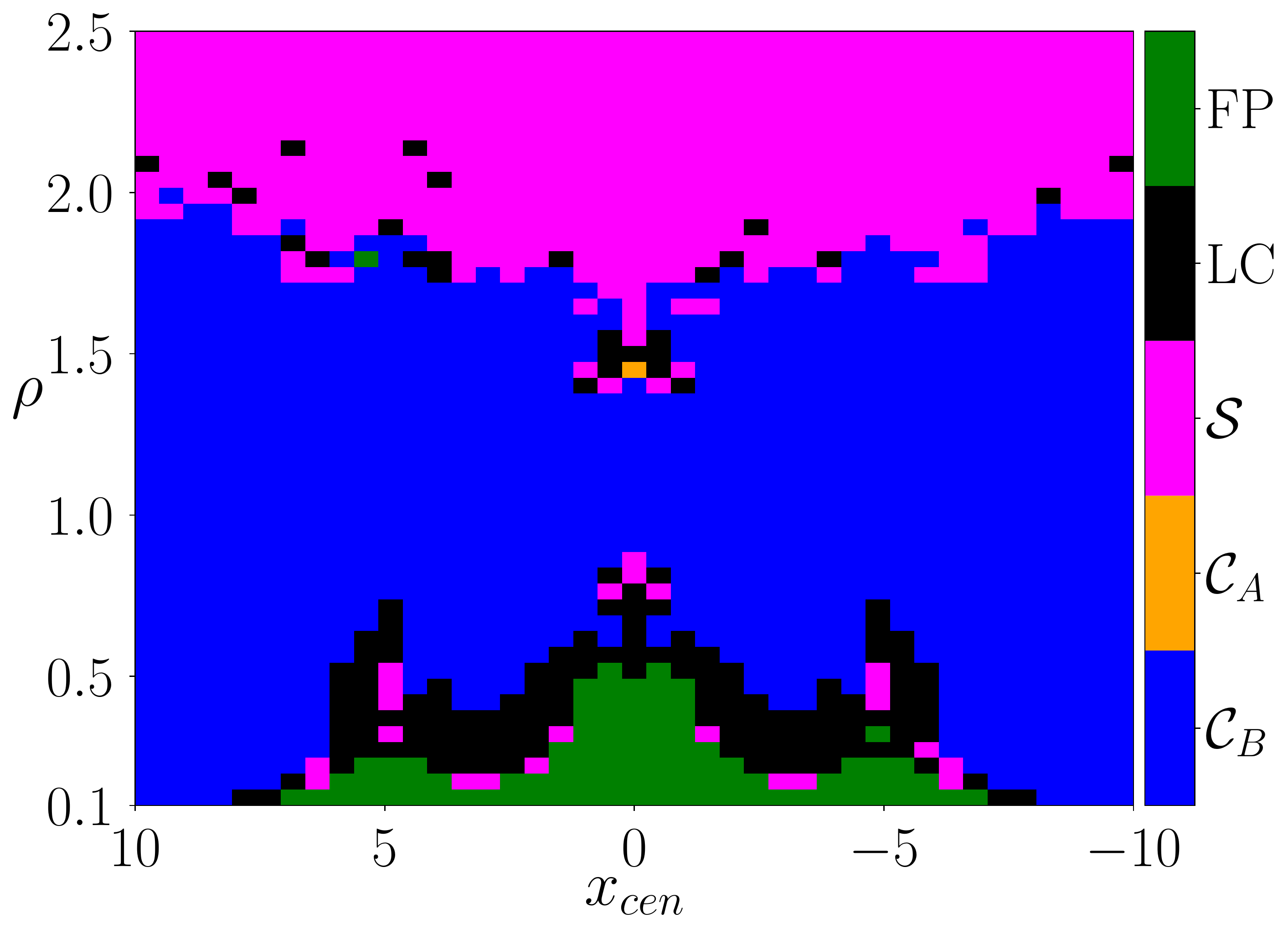}
        \caption{$\pzc_{B}$ prediction}
        \label{fig:Rerr_oppdir}
    \end{subfigure}
    \caption{Characterising the prediction of (a) $\pzc_{A}$ and (b) $\pzc_{B}$ in the $\left( x_{cen}, \, \rho \right)$-plane. Colouring scheme is as outlined in Table\,\ref{tab:ColorScheme}.}
    \label{fig:LR_oppdir}
\end{figure*}

\begin{figure}
    \centering
    \includegraphics[width=0.48\textwidth]{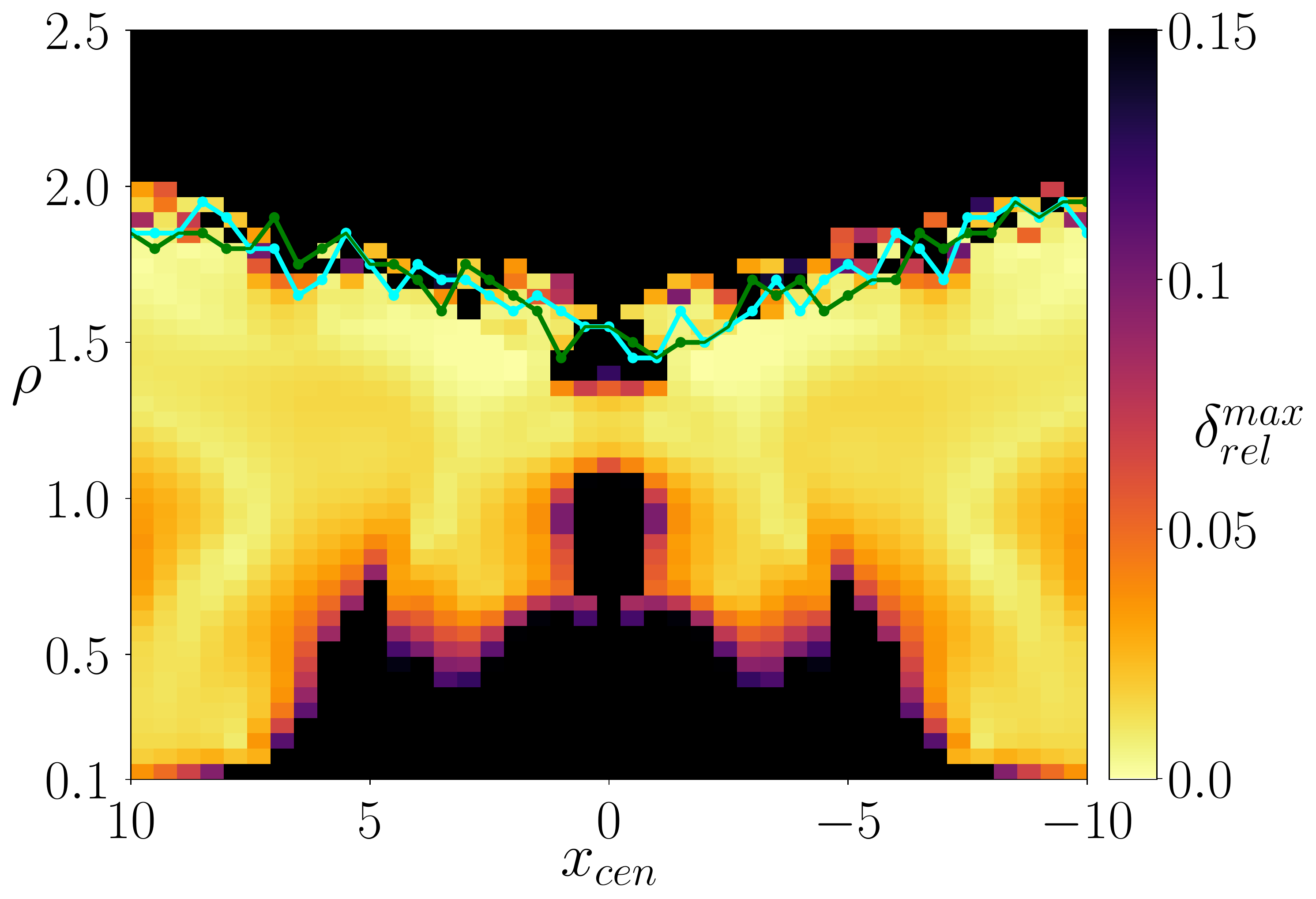}
    \caption{Regions of multifunctionality in $\left( x_{cen}, \, \rho \right)$-plane with error $\delta_{rel}^{max}$ for $\pzc_{A}$ and $\pzc_{B}$ rotating in opposite directions. Green and cyan curves highlight where $\lambda_{\text{MAX}}\,>\,\lambda_{C}$ for the first time for $\hat{\pzc}_{A}$ and $\hat{\pzc}_{B}$ respectively.}
    \label{fig:SD_MF_regions}
\end{figure}

The relationship between overlapping training data, multifunctionality, and $\rho$ becomes evident in Figs.\,\ref{fig:LR_oppdir}-\ref{fig:SD_MF_regions} where a `Goldilocks effect' is found to occur. By this it is meant that if $\rho$ is too small or too large then the closed-loop RC fails to achieve multifunctionality. When the orbits are furthest apart we see that the closed-loop RC is multifunctional even at small values of $\rho$. However, as the orbits are moved closer together and begin to overlap (i.e. when $\pzc_{A} \cap \, \pzc_{B} \neq \emptyset$), then the closed-loop RC requires a much larger $\rho$ in order to distinguish between the cycles and prevent the predicted trajectories from decaying to a fixed point or some other limit cycle. At the same time, if $\rho$ is too large the closed-loop RC fails to reconstruct either orbit regardless the investigated values of $x_{cen}$ and displays evidence of chaos.

We comment that from further inspection there are some subtle differences in the results for different $\textbf{M}$ and $\textbf{W}_{in}$ however this Goldilocks effect is a consistent feature.

Furthermore, it is important to note that a near symmetric set of results occur about $x_{cen} = 0$ in all Figs.\,\ref{fig:LR_oppdir}-\ref{fig:SD_MF_regions}, some insight to this is provided in Appendix\,\ref{ssec:SD_symm}.

\subsubsection{$\pzc_{A}$ and $\pzc_{B}$ rotating in same direction}\label{sssec:RotateSameDir}

\begin{figure}
    \centering
    \includegraphics[width=0.48\textwidth]{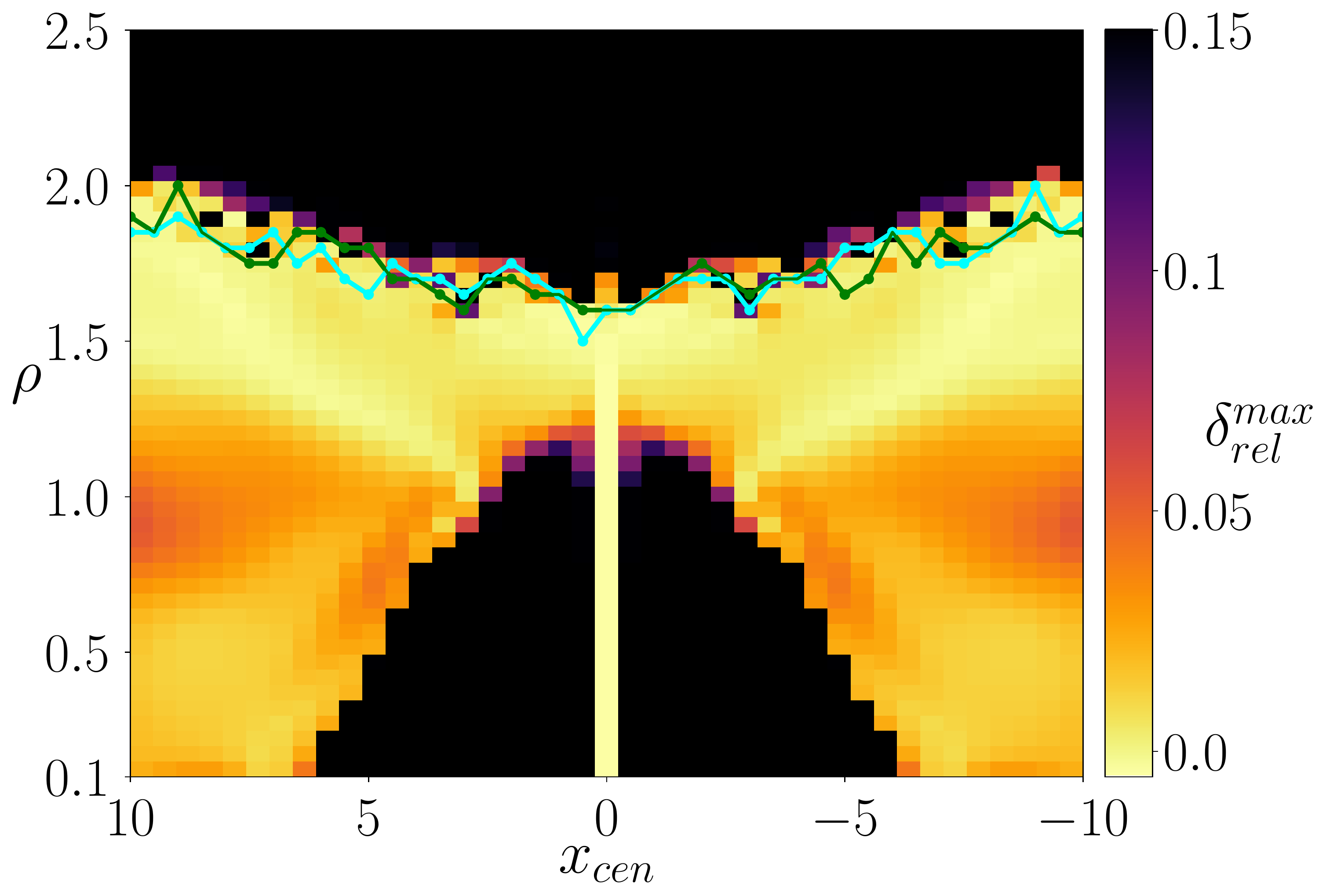}
    \caption{Regions of multifunctionality in $\left( x_{cen}, \, \rho \right)$-plane with error $\delta_{rel}^{max}$ for $\pzc_{A}$ and $\pzc_{B}$ rotating in the same direction. Green and cyan curves highlight where $\lambda_{\text{MAX}}\,>\,\lambda_{C}$ for the first time for $\hat{\pzc}_{A}$ and $\hat{\pzc}_{B}$ respectively.}
    \label{fig:SD_MF_regions_samedir}
\end{figure}

To provide a comparison, the same RC (constructed via the parameters in Table\,\ref{tab:RC_SDparams}) is trained to reconstruct a coexistence of $\pzc_{A}$ and $\pzc_{B}$ except, in this case, both cycles are rotating in the same direction. This is done by setting $b_{x}=b_{y}=5$ for both $\pzc_{A}$ and $\pzc_{B}$ in Eq.\,\eqref{eq:InputSys}. We examine the regions of multifunctionality in the $\left( x_{cen}, \, \rho \right)$-plane and these are computed using the same procedure that produced Fig.\,\ref{fig:SD_MF_regions}. 

The result of the above is presented in Fig.\,\ref{fig:SD_MF_regions_samedir} where a similar phenomenon to Fig.\,\ref{fig:SD_MF_regions} is shown to occur, once the overlap is introduced then multifunctionality is achieved only when $\rho$ is sufficiently large. However, the main difference between Fig.\,\ref{fig:SD_MF_regions_samedir} and Fig.\,\ref{fig:SD_MF_regions} is the trivial case for $x_{cen}=0$ as $\pzc_{A}$ and $\pzc_{B}$ are the exact same and the closed-loop RC has no issue reconstructing the dynamics for $\rho \in \left[0.1, 1.7 \right]$. For $\rho > 2.1$, the RC fails to become multifunctional $\forall x_{cen} \in \left[-10,10\right]$.

\subsection{Transition to Chaotic Dynamics}\label{ssec:EoCDyn_full}

\subsubsection{Dynamics before Chaos}\label{ssec:beforeEoC}

In Figs.\,\ref{fig:SD_MF_regions} and \ref{fig:SD_MF_regions_samedir} we find that for a given value of $x_{cen}$ and $\rho$, the smallest error, i.e. the smallest $\delta_{rel}^{max}$ value, occurs in many instances before the closed-loop RC's behaviour becomes chaotic for large $\rho$. This transition to chaos is quantified in Figs.\,\ref{fig:SD_MF_regions} and \ref{fig:SD_MF_regions_samedir} by the green and cyan curves which show the location where the largest Lyapunov exponent (LLE), denoted by $\lambda_{\text{MAX}}$, first becomes larger than $\lambda_{C} = 0.01$ when increasing $\rho$ for a given $x_{cen}$. We do this in order to indicate, within some measurement error, that $\lambda_{\text{MAX}}$ is positive.

Fig.\,\ref{fig:SD_MF_regions_samedir} shows that this transition to chaotic dynamics appears in the same vicinity of $\rho$ values as it does in Fig.\,\ref{fig:SD_MF_regions} for $|x_{cen}| \gtrsim 1.5$. As the same RC design is used in both cases, this suggests that the location of this transition may also arise as an intrinsic property of the RC rather than solely being related to the nature of the task itself. 

While the above comment is only speculation, some progress has been made in answering questions of this nature. For instance, \textcite{JiangLai19_RCspectral} conducted many extensive experiments which provide significant insight towards assessing the role of $\rho$ in the performance of a RC when training 100 random realisations of this RC in attractor reconstruction problems. The central result from Jiang and Lai's experiments is that, for a given RC, there exists an interval of $\rho$ values where optimal or near-optimal performance is achieved. A similar phenomenon is seen to occur in Figs.\,\ref{fig:SD_MF_regions} and \ref{fig:SD_MF_regions_samedir} which compliments the work of Jiang and Lai and also offers additional information as Figs.\,\ref{fig:SD_MF_regions} and \ref{fig:SD_MF_regions_samedir} both show that as the difficulty of the task is varied (by changing $x_{cen}$) then the width of these intervals change accordingly.

We remark that in studies involving the training of input-driven RCs (after the training the readout function does not replace the driving input like in Eq.\,\eqref{eq:PredRes}), a similar phenomenon known as the `edge of chaos' or the `edge of stability' and has been the focus of many investigations which suggest that RCs achieve optimal computational capacity at the edge of chaos \cite{bertschinger04CompEoC,boedecker12informationEoC,teuscher2022revisitingEoC}. 
However, \textcite{Carroll20EoC} presents several examples of RCs which perform better on certain problems just prior to crossing the edge of chaos.

\subsubsection{Dynamics after Chaos}\label{ssec:beyondEoC}

\begin{figure}
    \centering
    \includegraphics[width=0.479\textwidth]{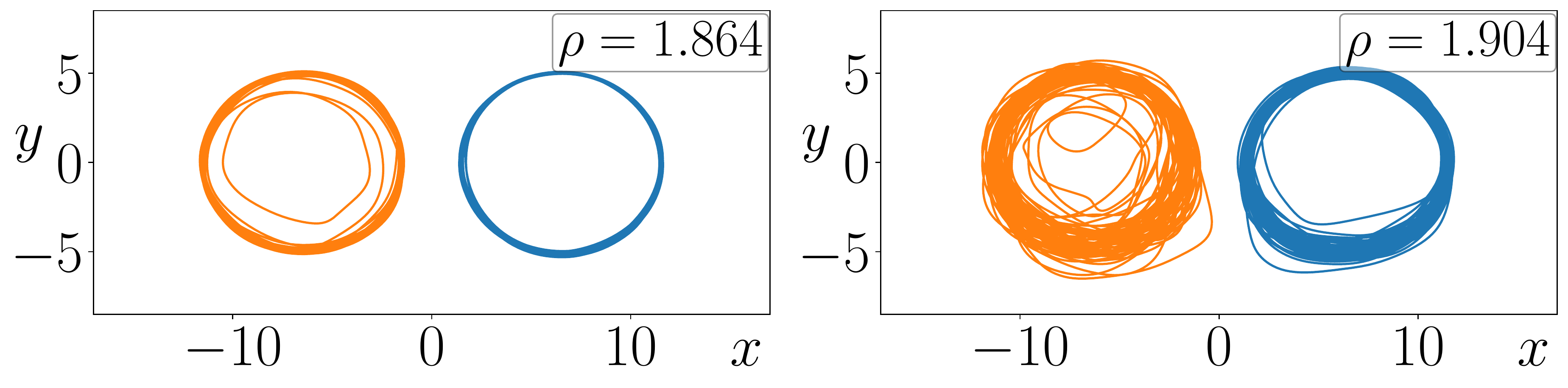}
    \includegraphics[width=0.479\textwidth]{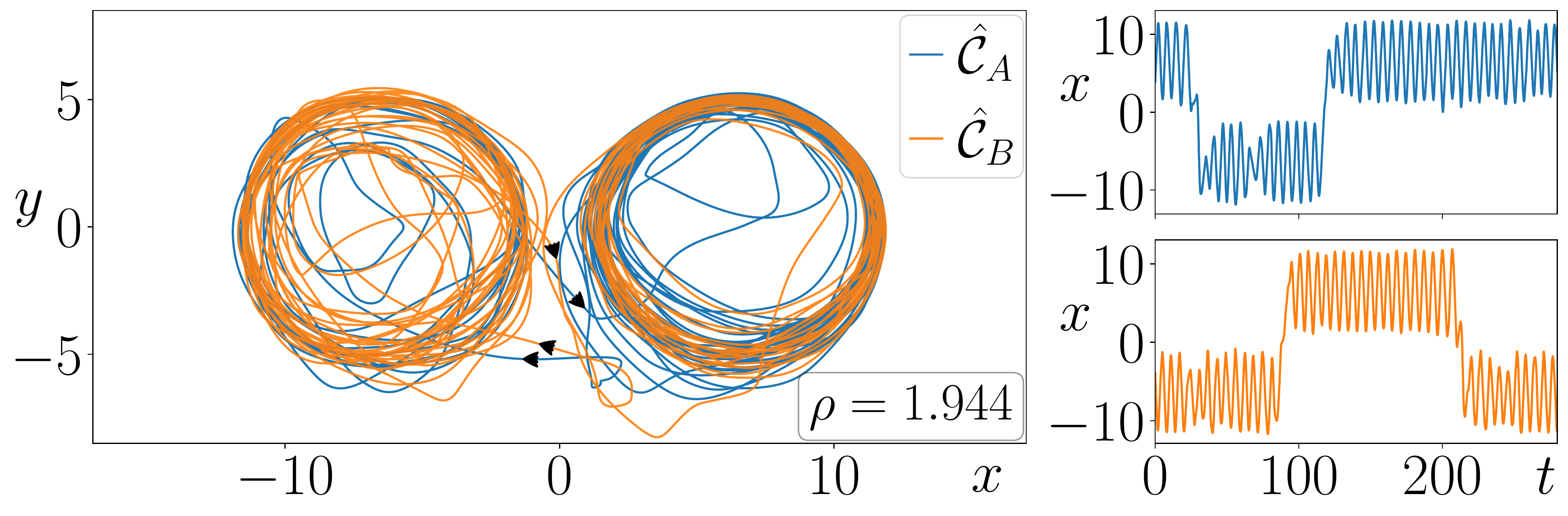}
    \includegraphics[width=0.479\textwidth]{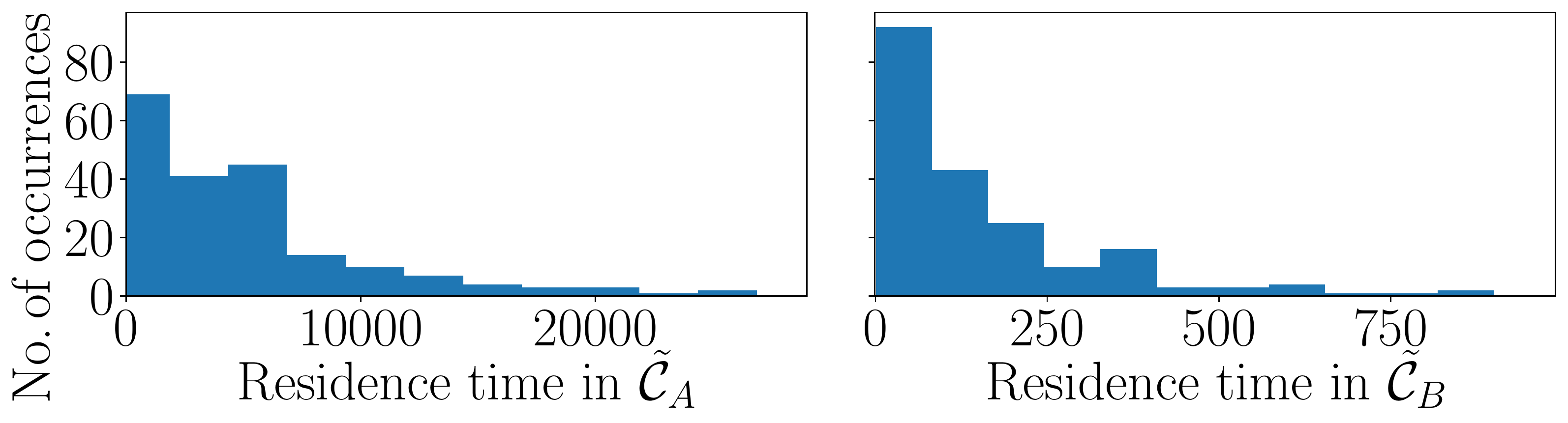}
    \caption{Closed-loop RC's reconstruction of $\pzc_{A}$ and $\pzc_{B}$ when $x_{cen}=6.5$ for $\rho = 1.864$ (top left), $1.904$ (top right), $1.944$ (middle left). Middle right: Reconstructed $x$ variable time traces highlight the nature of the RC's itinerant dynamics for $\rho = 1.944$. Distribution of residence times the RC spends in either $\tilde{\pzc}_{A}$ (bottom left) or $\tilde{\pzc}_{B}$ (bottom right) when computing a trajectory of the closed-loop RC up to $t=1000000$.}
    \label{fig:CI_example}
\end{figure}

In some cases, after the transition to chaotic dynamics, an interesting sequence of events is found to occur when increasing $\rho$ past the point of where multifunctionality is no longer achieved. As an example, Fig.\,\ref{fig:CI_example} illustrates the changes in the dynamics of Eq.\,\eqref{eq:PredRes} as viewed from $\mathbb{P}$ as $\rho$ is increased from $1.864$ to $1.944$ in the case where $x_{cen}=6.5$. 
In relation to Figs.\,\ref{fig:Lerr_oppdir} and \ref{fig:Rerr_oppdir}, the long-term dynamics of the closed-loop RC when trained at these $\rho$ values and initialised from $\hat{\boldsymbol{r}}_{\left(\pzc_{A}\right)}$ or $\hat{\boldsymbol{r}}_{\left(\pzc_{B}\right)}$ are colour-coded as magenta, indicating that the state of Eq.\,\eqref{eq:PredRes} does not settle down to periodic motion.

Fig.\,\ref{fig:CI_example} shows that for $\rho=1.864$ and $1.904$, the dynamics of both $\hat{\pzc}_{A}$ and $\hat{\pzc}_{B}$ begin to occupy larger and larger regions of $\mathbb{P}$. This behaviour continues up to a point where the closed-loop RC loses multistability and the state of Eq.\,\eqref{eq:PredRes} begins to switch from one orbit to the other and back again, shown here for $\rho = 1.944$.

While further analysis is required in order to determine the route of these switching dynamics, it is possible that the closed-loop RC has entered into the dynamical regime known as `chaotic itinerancy'. Chaotic itinerancy\cite{tsuda1987CImemory,ikeda1989itinerancy,kaneko1990CI}, describes the dynamical phenomenon of an autonomous switching process whereby the state of a given autonomous dynamical system switches between several `attractor ruins' or `quasi-attractors'. These were previously coexisting attractors that are now connected to form one global attractor. The quasi-attractors retain much of their original features except trajectories on these quasi-attractors leak into each other. In the case of Fig.\,\ref{fig:CI_example} for $\rho=1.944$, $\hat{\pzc}_{A}$ and $\hat{\pzc}_{B}$ exhibit these quasi-attractor properties. The state of the closed-loop RC wanders between visiting different regions of $\mathbb{P}$ associated with attractor ruins of the previously coexisting $\hat{\pzc}_{A}$ and $\hat{\pzc}_{B}$. The corresponding time traces of the reconstructed $x$ variable (as shown in the middle right in Fig.\,\ref{fig:CI_example}) provides a clearer picture of the switching dynamics between attractor ruins as emphasised by the arrows in the associated $\mathbb{P}$ (middle left panel in Fig.\,\ref{fig:CI_example}). After simulating the dynamics of Eq.\,\eqref{eq:PredRes} up to $t=1000000$ in this metastable regime from $\hat{\boldsymbol{r}}\left( 0 \right) = \boldsymbol{r}_{(\pzc_{A})}\left( t_{train} \right)$, we find that there are over $200$ switchings between the two quasi-attractors. In the bottom panels we also plot the distribution of time spent in each quasi-attractor, labelled here as $\tilde{\pzc}_{A}$ and $\tilde{\pzc}_{B}$. These plots show that the state of the closed-loop RC spends a far greater amount of time in $\tilde{\pzc}_{A}$ as opposed to $\tilde{\pzc}_{B}$. In future work we aim to explore the relationship between $\rho$ and the resulting distribution of residence times in each quasi-attractor.

It is important to stress here that there is no external input used to generate the switching dynamics shown in Fig.\,\ref{fig:CI_example} as it is an intrinsic property of the closed-loop RC after it loses multifunctionality and subsequently multistability. This differs from previous cases where RCs were trained to switch between learning different chaotic and periodic attractors based on an external input like, for example, in the work of Inoue \textit{et al.} in \cite{Nakajima20designing_itinerancy} and Lu and Bassett in \cite{LuBassett20_switching_learning}. In these scenarios, $\textbf{W}_{out}$ is updated online in response to a given driving input from an individual attractor and only one attractor is shown to exist at each point in time.

\subsection{Achieving multifunctionality when \texorpdfstring{$x_{cen}=0$}{TEXT}}\label{ssec:RCdynamics_xcen0}

\subsubsection{Dynamics in \texorpdfstring{$\mathbb{P}$}{TEXT}}\label{sssec:DynamicsInP}

From a dynamical systems perspective it is remarkable that the closed-loop RC reconstructs a coexistence of $\pzc_{A}$ and $\pzc_{B}$ when these share similar regions of state space. However, the most striking result from Fig.\,\ref{fig:SD_MF_regions} is that the closed-loop RC reconstructs a coexistence of these orbits even when they are completely overlapping, i.e., when $x_{cen}=0$. In this extreme scenario, multifunctionality is achieved for a small range of $\rho$ values, specifically for $1.1 \leq \rho \leq 1.4$. The most accurate reconstruction of both orbits is achieved when $\rho = 1.25$. 

\begin{figure*}
    \centering
    \includegraphics[width=0.99\textwidth]{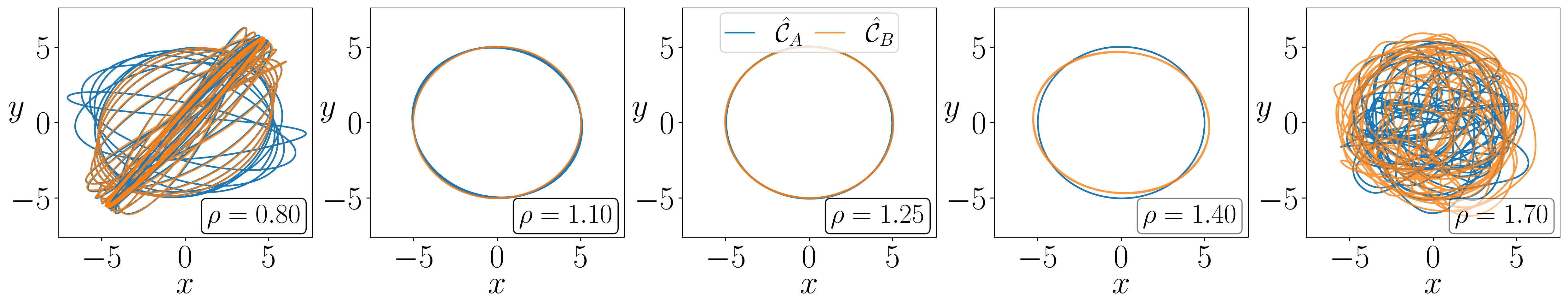}
    \caption{RC's reconstruction of $\pzc_{A}$ and $\pzc_{B}$ when $x_{cen}=0$ for $\rho = 0.80$, $1.10$, $1.25$, $1.40$, and $1.70$.}
    \label{fig:Circles_rho_08_17}
\end{figure*}

Fig.\,\ref{fig:Circles_rho_08_17} highlights how values of $\delta_{rel}^{max}$ differ in regions where multifunctionality was achieved, for instance, $\delta_{rel}^{max} \approx 0.062, 0.019, 0.125$ for $\rho=1.1, \,1.25, \,1.4$ respectively. These account for the slightly off-circular reconstructions of $\pzc_{A}$ for $\rho = 1.1$ and $\pzc_{B}$ for $\rho = 1.4$ as shown in Fig.\,\ref{fig:Circles_rho_08_17}.

Examples where the closed-loop RC fails to achieve multifunctionality just outside the small window of $\rho$ values identified in Fig.\,\ref{fig:SD_MF_regions} are also illustrated in Fig.\,\ref{fig:Circles_rho_08_17} for $\rho = 1.7$ and $0.8$. Here the chaotic-like behaviour of the closed-loop RC is seen for $\rho = 1.7$. What is important to note in this example is that shape swept out by the closed-loop RC's trajectories when initialised from both $\hat{\boldsymbol{r}}(0)=\boldsymbol{r}_{\left(\pzc_{A}\right)}(t_{\text{listen}})$ and $\hat{\boldsymbol{r}}(0)=\boldsymbol{r}_{\left(\pzc_{B}\right)}(t_{\text{listen}})$ is bounded to regions nearby the original circular orbits, however the direction of rotation is not preserved.

Fig.\,\ref{fig:Circles_rho_08_17} also shows that when the reconstruction of both $\pzc_{A}$ and $\pzc_{B}$ fails for $\rho = 0.8$ then the state of Eq.\,\eqref{eq:PredRes} approaches the same limit cycle which bears no resemblance to either $\pzc_{A}$ and $\pzc_{B}$. The role this limit cycle plays in how the closed-loop RC solves the seeing double problem is assessed in Sec.\,\ref{sec:SDdynamicsanal}.

Much of the following results in this paper focus on exploring the dynamics of the closed-loop RC in this `simplest hardest' case, `simplest' in terms of the dynamics being limit cycles and `hardest' as the training data is entirely overlapping. We comment that, from analysis not shown in the present paper, the small range of $\rho$ values where multifunctionality is achieved is relatively consistent across several random realisations of $\textbf{M}$ and $\textbf{W}_{in}$.

\subsubsection{Dynamics in \texorpdfstring{$\mathbb{S}$}{TEXT}}\label{sssec:DynamicsInS}

\begin{figure}
    \centering
    \includegraphics[width=0.45\textwidth]{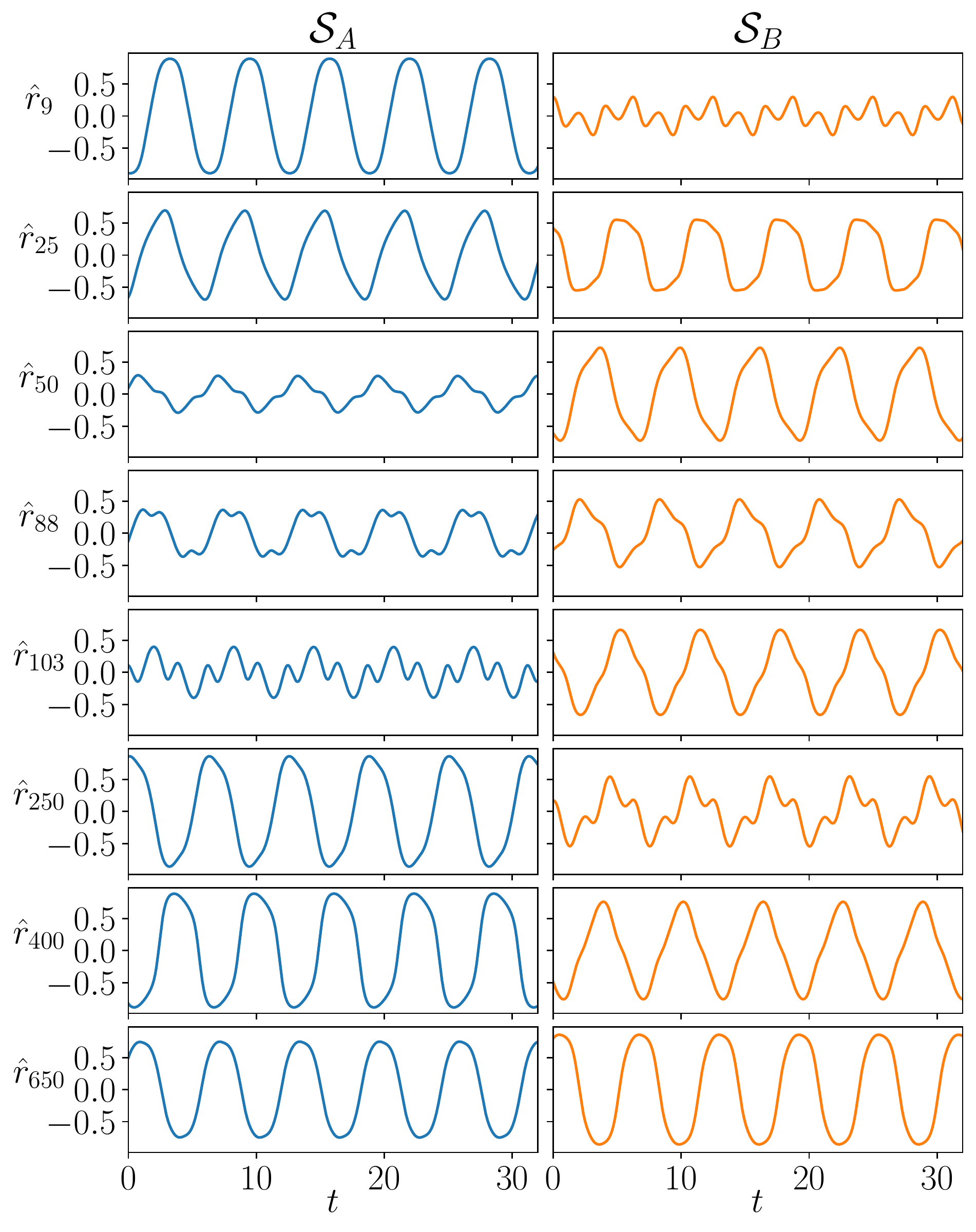}
    \caption{Illustrating how the dynamics of individual neurons can differ depending on whether the multifunctional RC approaches either $\pzs_{A}$ or $\pzs_{B}$ in order to reconstruct $\pzc_{A}$ or $\pzc_{B}$ when $x_{cen}=0$ for $\rho = 1.25$.}
    \label{fig:InternalDynamics_rho125}
\end{figure}

Before discussing these further results we wish to provide some insight towards the internal dynamics of the multifunctional RC and how the behaviour of the individual neurons differ when reconstructing either $\pzc_{A}$ or $\pzc_{B}$ in the case where $x_{cen}=0$. Fig.\,\ref{fig:InternalDynamics_rho125} depicts a representative example of how the dynamics of the same set of neurons can significantly differ in $\mathbb{S}$ as the state of the multifunctional RC approaches either $\pzs_{A}$ or $\pzs_{B}$ while the only difference between the projected dynamics (the dynamics of the reconstructed attractors, $\hat{\pzc}_{A}$ and $\hat{\pzc}_{B}$), is the direction of rotation. To place a stronger emphasis on the fact that we are now examining the dynamics of $\pzs_{A}$ and $\pzs_{B}$, in this section we refer to the corresponding states of Eq.\,\eqref{eq:PredRes} as $\hat{\boldsymbol{r}}_{\left(\pzs_{A}\right)}(t)$ and $\hat{\boldsymbol{r}}_{\left(\pzs_{B}\right)}(t)$. Fig.\,\ref{fig:InternalDynamics_rho125} highlights how the state of the multifunctional RC's neurons $9, 25, 50, 88, 103, 250, 400$, and $650$ evolve over time (starting from when Eq.\,\eqref{eq:PredRes} is initialised with $\hat{\boldsymbol{r}}\left( 0 \right) = \hat{\boldsymbol{r}}_{\left(\pzs_{A}\right)}(0) = \boldsymbol{r}_{(\pzc_{A})}\left( t_{train} \right)$ or $\hat{\boldsymbol{r}}\left( 0 \right) = \hat{\boldsymbol{r}}_{\left(\pzs_{B}\right)}(0) = \boldsymbol{r}_{(\pzc_{B})}\left( t_{train} \right)$) for the case of $\rho=1.25$ whose dynamics in $\mathbb{P}$ are shown in Fig.\,\ref{fig:Circles_rho_08_17}. It is found that many neurons behave similarly to each other when the state of the multifunctional RC approaches either $\pzs_{A}$ or $\pzs_{B}$, see the dynamics of $\hat{r}_{650}$ in Fig.\,\ref{fig:InternalDynamics_rho125} as an example. At the same time, Fig.\,\ref{fig:InternalDynamics_rho125} shows that the behaviour of many other neurons may bare no resemblance to one another. While all neurons exhibit periodic dynamics some neurons oscillate with multiple local maxima, however, what is consistent across all neurons is that they each have the same period of rotation.

\begin{figure*}
    \centering
    \includegraphics[width=0.85\textwidth]{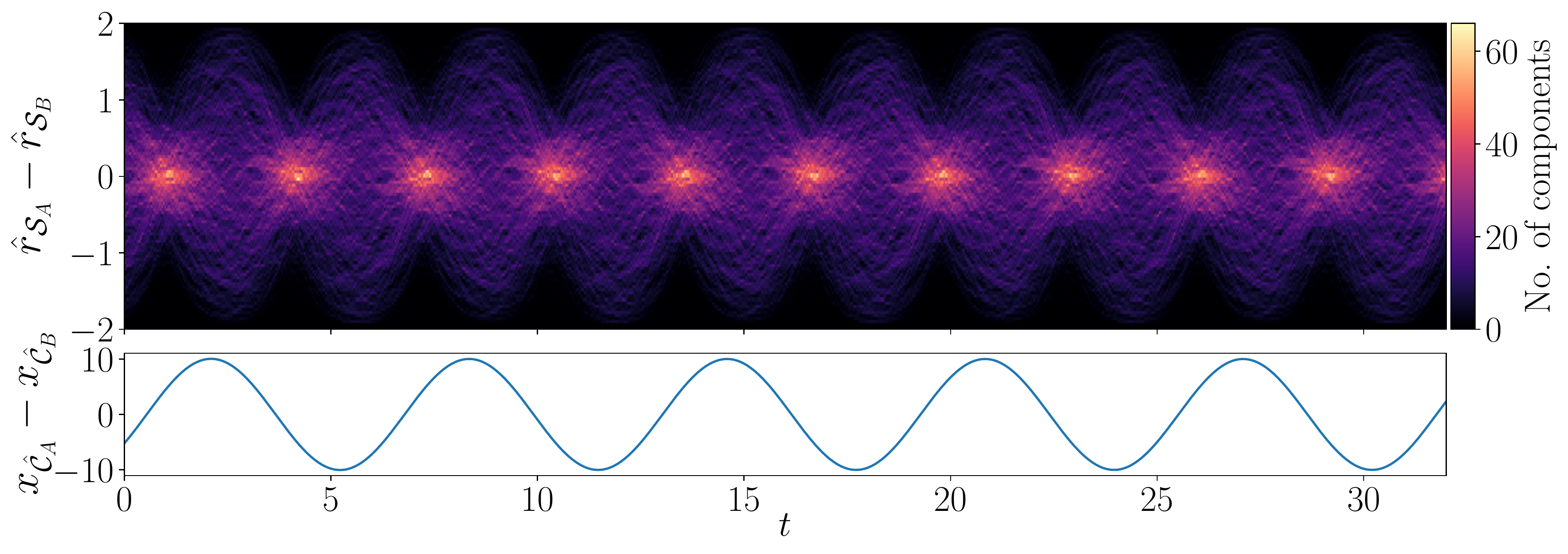}
    \caption{For $x_{cen}=0$ and $\rho = 1.25$, top panel illustrates the difference between the state of the multifunctional RC when approaching $\pzs_{A}$ or $\pzs_{B}$, lower panel depicts the difference between the respective $x$-variables when reconstructing $\pzc_{A}$ or $\pzc_{B}$.}
    \label{fig:InternalDynamicsAllNeurons_rho125}
\end{figure*}

In order to provide further insight to the dynamics of the individual neurons we now continue our study of the multifunctional RC in the above case where $\rho=1.25$ and $x_{cen}=0$. We now examine how the dynamics of the these neurons differ when the state of the multifunctional RC approaches the same point on each reconstructed attractor in $\mathbb{P}$, i.e. when $x_{\hat{\pzc}_{A}}(t)-x_{\hat{\pzc}_{B}}(t)=0$ where $x_{\hat{\pzc}_{A}}(t)$ and $x_{\hat{\pzc}_{A}}(t)$ denote the location of the $x$-variable in $\mathbb{P}$ at a given time $t$ (starting from when Eq.\,\eqref{eq:PredRes} is initialised) when reconstructing $\pzc_{A}$ and similarly for $\pzc_{B}$. Note that the associated $y$-variables are always equal for all values of $t$, this is why we need only look at the $x$-variable. Using the same notation convention as for $x_{\hat{\pzc}_{A}}$ and $x_{\hat{\pzc}_{B}}$, in the top panel of Fig.\,\ref{fig:InternalDynamicsAllNeurons_rho125} we plot how the difference between $\hat{r}_{\pzs_{A}}$ and $\hat{r}_{\pzs_{B}}$ evolves over time. To convey how these $N=1000$ neurons behave, we do this by generating a histogram at each time-step of the simulation which describes how many components of $\hat{r}_{\pzs_{A}}(t) - \hat{r}_{\pzs_{B}}(t)$ have a particular value at a given time $t$. These values are confined to be within $-2$ and $2$ (which is the maximum and minimum allowed values due to the $\tanh$ activation function in the RC design) in intervals of $0.02$ (which is the bin width of the histograms and is chosen to show sufficient details of the dynamics at an appropriate resolution). In the lower panel we plot the evolution of the corresponding difference between $x_{\hat{\pzc}_{A}}$ and $x_{\hat{\pzc}_{B}}$.

From Fig.\,\ref{fig:InternalDynamicsAllNeurons_rho125} it is evident that at times when $x_{\hat{\pzc}_{A}}(t) - x_{\hat{\pzc}_{B}}(t) \approx 0$, even though many neurons fire identically the crucial aspect for multifunctionality is that not all neurons fire identically as evidenced by the relatively wide distribution of the values of the components of $\hat{r}_{\pzs_{A}}(t) - \hat{r}_{\pzs_{B}}(t)$ at these particular times. We also see here that the greater the distance between points on each reconstructed attractor, i.e. when $|x_{\hat{\pzc}_{A}} - x_{\hat{\pzc}_{B}}|$ grows larger, the greater the difference between $\hat{r}_{\pzs_{A}}$ and $\hat{r}_{\pzs_{B}}$ (the wider the distribution of the firing patterns). 

\subsection{Projected basins of attraction for \texorpdfstring{$x_{cen}=0$}{TEXT}}\label{sec:SD_BoAs}

The focus of this section is to provide some insight towards the nature of the basins of attraction created by the closed-loop RC in order to solve the seeing double problem in the extreme case of $x_{cen}=0$.

\subsubsection{Method of generating the basins}\label{sssec:GenBasins}

Considering the dimension of the closed-loop RC ($N=1000$), to determine the precise location of the basins of attraction in $\mathbb{S}$ for both $\pzs_{A}$ and $\pzs_{B}$ is highly expensive. Instead, to get some insight about the structure of these basins, we examine how the state of the closed-loop RC evolves from the open-loop RC's representation of different points in $\mathbb{R}^{D}$. From this it is possible to see how $\mathbb{P}$ is split to accommodate a coexistence of $\hat{\pzc}_{A}$ and $\hat{\pzc}_{B}$ when $\pzc_{A}$ and $\pzc_{B}$ are completely overlapping.

To be more specific, for a given $\rho$ and $x_{cen}=0$, we drive Eq.\,\eqref{eq:ListenRes} with a point in $\mathbb{R}^{D}$ from $t=0$ to $t=t_{\text{listen}}$. After this time, the state of Eq.\,\eqref{eq:ListenRes} converges to a representation of the point in $\mathbb{R}^{N}$. This response vector is labelled as $\boldsymbol{r}_{D} = \boldsymbol{r}(t_{\text{listen}})$. The corresponding closed-loop RC is then initialised with $\boldsymbol{\hat{r}}(0) = \boldsymbol{r}_{D}$ and the long-term behaviour of Eq.\,\eqref{eq:PredRes} is characterised from this initial condition.

\subsubsection{Basins for \texorpdfstring{$\rho = 1.25$}{TEXT}}\label{sssec:rho125}

\begin{figure}
    \centering
    \includegraphics[width=0.4\textwidth]{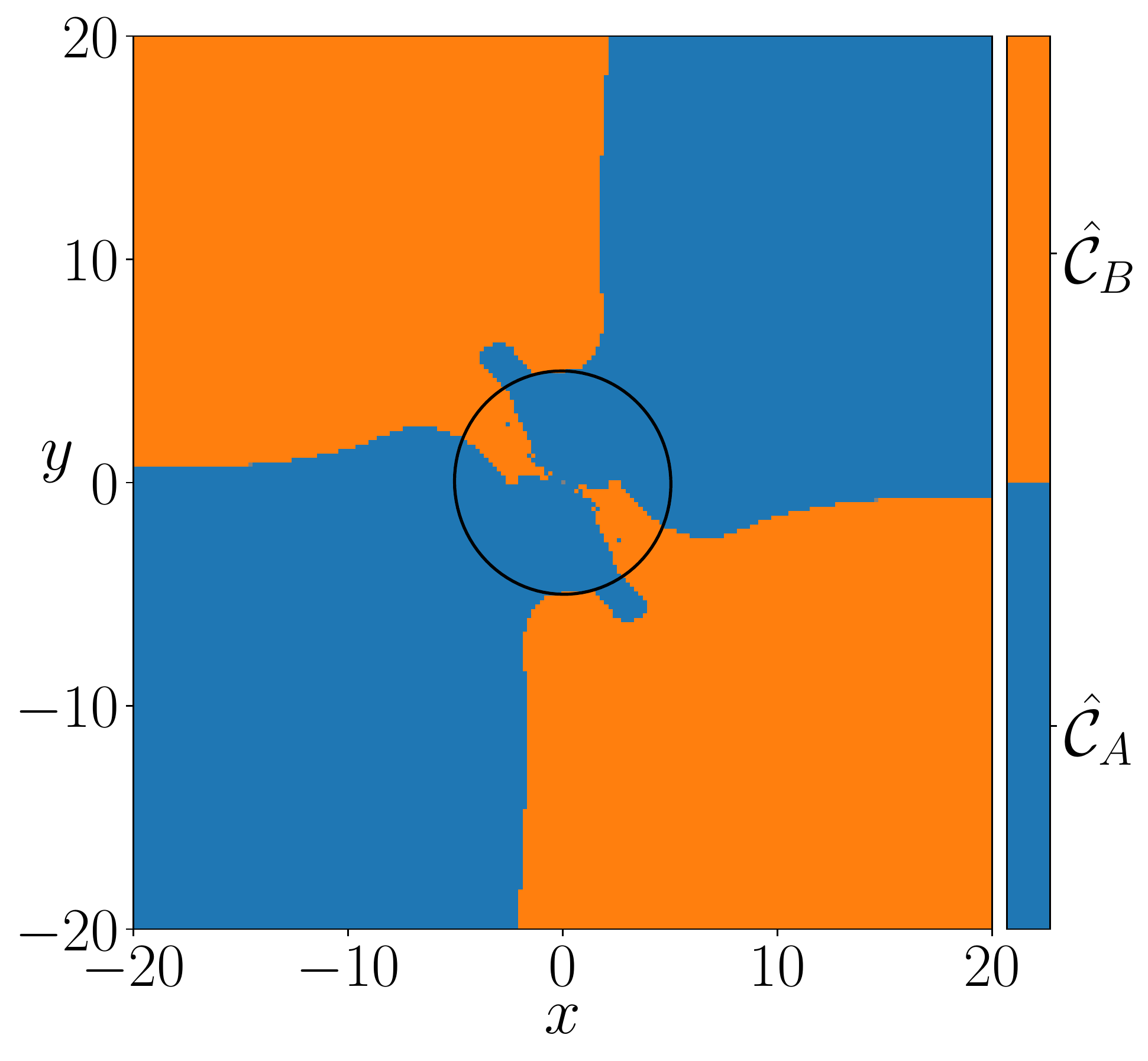}
    \caption{Basins of attraction for $\hat{\pzc}_{A}$ and $\hat{\pzc}_{B}$
    as viewed from $\mathbb{P}$ for $\rho = 1.25$. The respective attractors are plotted in black.} 
    \label{fig:BoA_rho_125}
\end{figure}

The projected basins of attraction for $\rho = 1.25$, which is the case where multifunctionality was achieved with the best accuracy, is shown in Fig.\,\ref{fig:BoA_rho_125}. Each point in Fig.\,\ref{fig:BoA_rho_125} describes which attractor Eq.\,\eqref{eq:PredRes} approaches in $\mathbb{P}$ for a given point in $\mathbb{R}^{D}$. If the state of Eq.\,\eqref{eq:PredRes} approaches $\hat{\pzc}_{A}$ then the point is coloured blue and similarly for $\hat{\pzc}_{B}$ the point is coloured orange. There is one exception, when Eq.\,\eqref{eq:ListenRes} is driven with $\left( 0, 0 \right)$ its state remains at the origin and so too does Eq.\,\eqref{eq:PredRes}, this is coloured grey. The reconstructed attractors, $\hat{\pzc}_{A}$ and $\hat{\pzc}_{B}$, are both plotted in black.

Fig.\,\ref{fig:BoA_rho_125} shows how $\hat{\pzc}_{A}$ and $\hat{\pzc}_{B}$ share $\mathbb{P}$ and that the number of points which converge to either attractor is not evenly distributed. In particular, most of the points closest to the location of the reconstructed attractors in $\mathbb{P}$ tend towards $\hat{\pzc}_{A}$. Furthermore, there is no single curve which entirely divides the boundary between the basins of $\hat{\pzc}_{A}$ and $\hat{\pzc}_{B}$ in Fig.\,\ref{fig:BoA_rho_125}, at the same time, a more regular or complex pattern may occur in $\mathbb{S}$. Nonetheless, the method used to produce Fig.\,\ref{fig:BoA_rho_125} provides an indication that the basin of attraction for the attractor $\pzs_{A}$ in $\mathbb{S}$, which when projected to $\mathbb{P}$ describes $\hat{\pzc}_{A}$, is of a larger volume than the basin for $\pzs_{B}$ in $\mathbb{S}$. We comment that this result may not be the case across all random realisations of $\textbf{M}$ and $\textbf{W}_{in}$.

\subsubsection{Basins for \texorpdfstring{$\rho = 0.9, 0.8, 0.4, 0.3$}{TEXT}}\label{sssec:rho09_03}

In order to gain a broader understanding of the closed-loop RC's dynamics when it fails to distinguish between $\pzc_{A}$ and $\pzc_{B}$, we now study the different basins of attraction which exist before the closed-loop RC exhibits multifunctionality. The same procedure which generated the basins for $\rho = 1.25$ is now used to find the projected basins of attraction in $\mathbb{P}$ when setting $\rho = 0.9, 0.8, 0.4$ and, $0.3$ in the case of $x_{cen}=0$.

\begin{figure*}
    \centering
    \begin{subfigure}{0.246\textwidth}
        \centering
        \includegraphics[width=\textwidth]{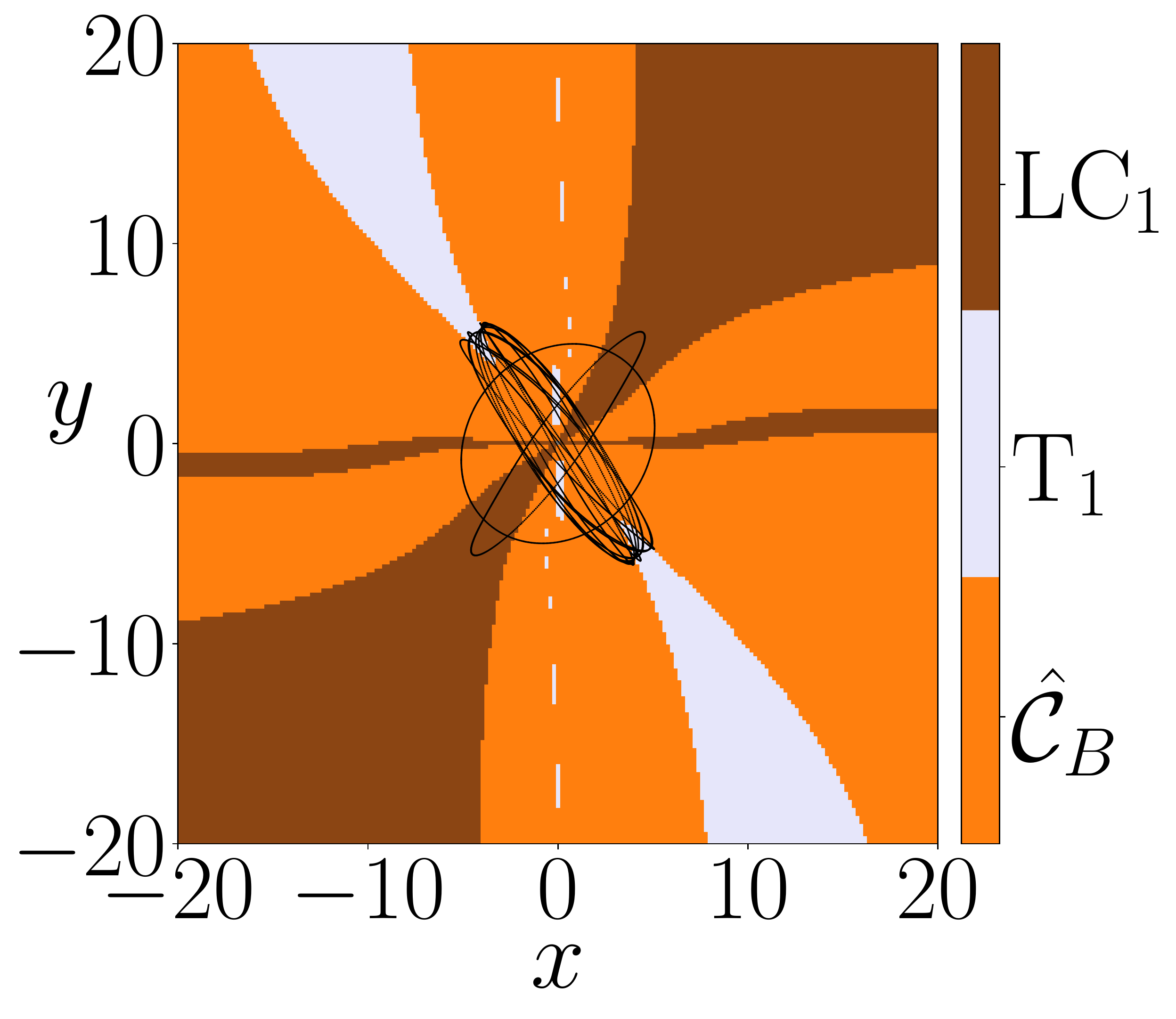}
        \caption{$\rho = 0.9$}
        \label{fig:BoA_rho_09}
    \end{subfigure}
    \begin{subfigure}{0.246\textwidth}
        \centering
        \includegraphics[width=\textwidth]{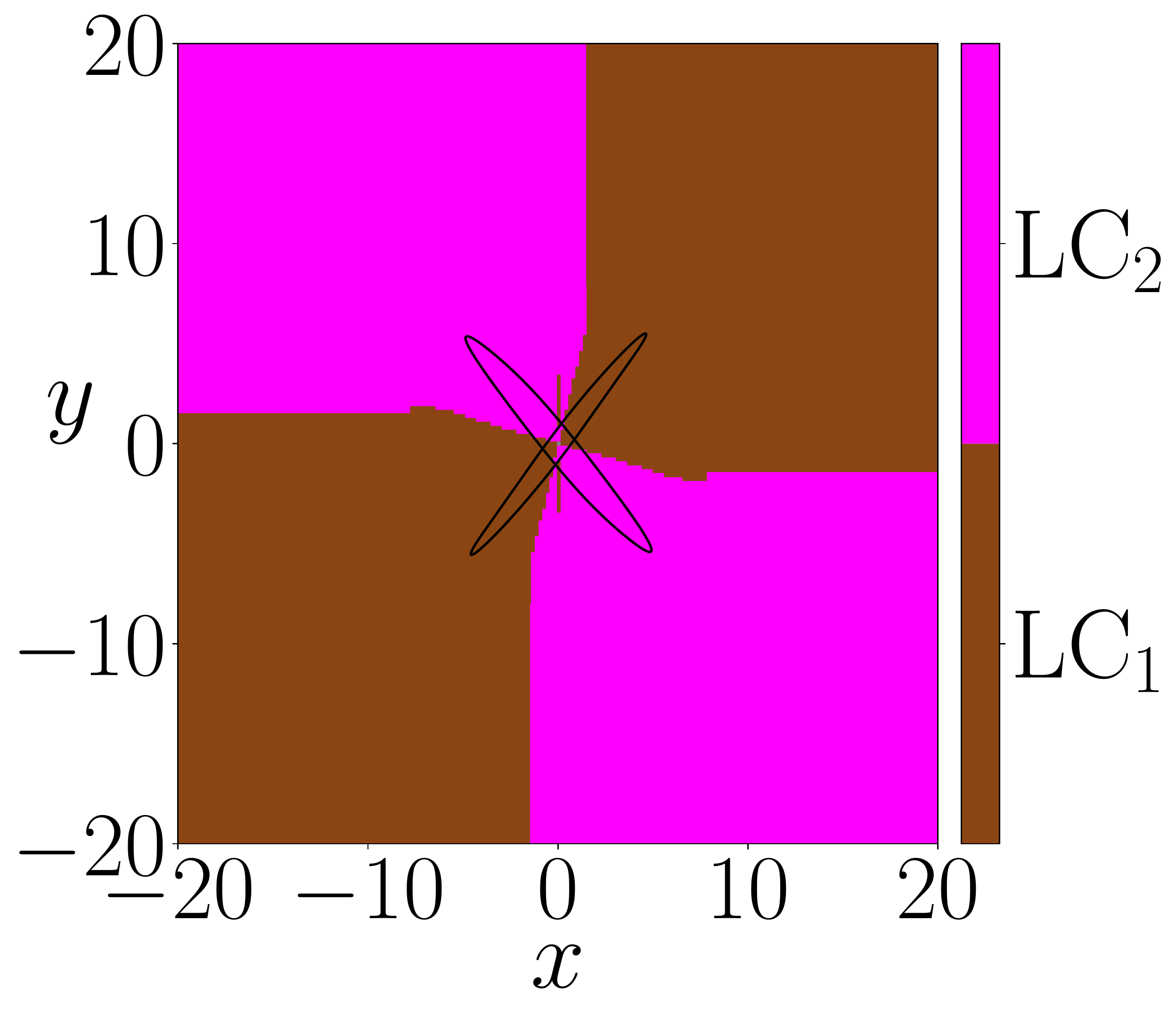}
        \caption{$\rho = 0.8$}
        \label{fig:BoA_rho_08}
    \end{subfigure}
    \begin{subfigure}{0.246\textwidth}
        \centering
        \includegraphics[width=\textwidth]{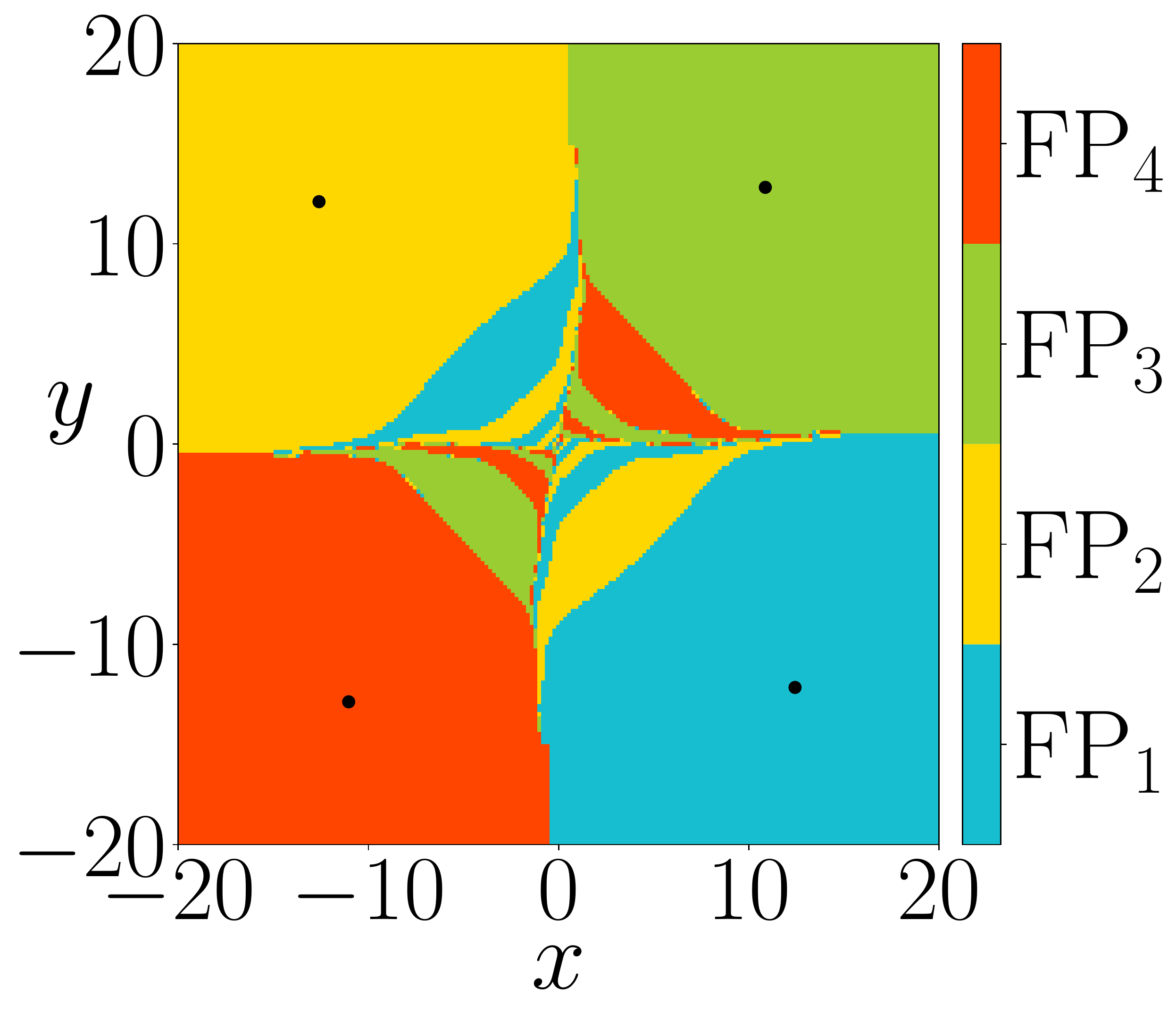}
        \caption{$\rho = 0.4$}
        \label{fig:BoA_rho_04}
    \end{subfigure}
    \begin{subfigure}{0.246\textwidth}
        \centering
        \includegraphics[width=\textwidth]{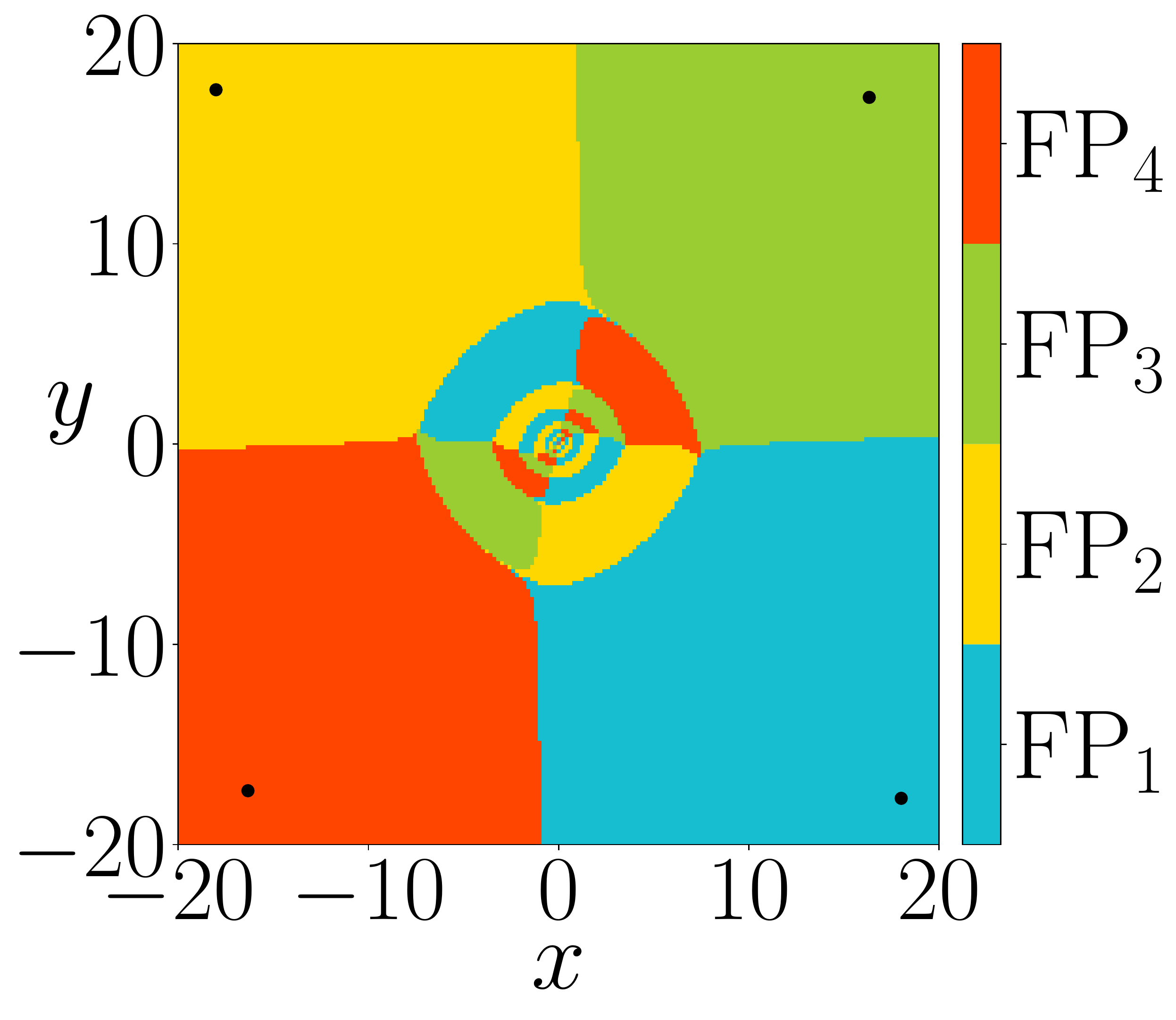}
        \caption{$\rho = 0.3$}
        \label{fig:BoA_rho_03}
    \end{subfigure}
    \caption{Basins of attraction for $\hat{\pzc}_{B}$, LC$_{1}$, T$_{1}$, LC$_{2}$, FP$_{1}$, FP$_{2}$, FP$_{3}$, and FP$_{4}$
    as viewed from $\mathbb{P}$ when setting $x_{cen}=0$ and (a) $\rho=0.9$, (b) $\rho=0.8$, (c) $\rho=0.4$, and (d) $\rho=0.3$. The respective attractors are plotted in black.}
    \label{fig:BoAs}
\end{figure*}

Prior to achieving multifunctionality, Fig.\,\ref{fig:LR_oppdir} shows that the closed-loop RC is able to reconstruct $\pzc_{B}$ before $\pzc_{A}$. For $\rho = 0.9$, Fig.\,\ref{fig:BoA_rho_09} illustrates that in addition to $\hat{\pzc}_{B}$ there are other basins of attraction for a limit cycle labelled as LC$_{1}$ and a quasi-periodic torus labelled as T$_{1}$. So while the aim is to reconstruct a coexistence of $\pzc_{A}$ and $\pzc_{B}$, there are other attractors that can be approached from an abundance of initial conditions. These untrained attractors, are attractors that exist in $\mathbb{P}$ but were not present during the training, these attractors were also found in \textcite{flynn2021multifunctionality}.

In the case where neither $\pzc_{A}$ or $\pzc_{B}$ can be reconstructed, for example, when $\rho = 0.8$ as seen in Figs.\,\ref{fig:LR_oppdir} and \ref{fig:Circles_rho_08_17}, Fig.\,\ref{fig:BoA_rho_08} shows that there are basins of attraction for two different limit cycles, LC$_{1}$ from Fig.\,\ref{fig:BoA_rho_09} and another limit cycle denoted as LC$_{2}$. So at one point before achieving multifunctionality there is a coexistence of two different limit cycles in $\mathbb{P}$.

When setting $\rho = 0.4$, Fig.\,\ref{fig:BoA_rho_04} reveals that four fixed points FP$_{1}$, FP$_{2}$, FP$_{3}$, and FP$_{4}$ are the only stable attractors which populate $\mathbb{P}$. Far away from the origin it can be seen that $\mathbb{P}$ is split into nearly equal quadrants. However, closer to the origin there is a greater interaction between the basins resulting in the emergence of a fractal pattern. A similar result is found when setting $\rho = 0.3$ in Fig.\,\ref{fig:BoA_rho_03} with fractal basins boundaries appearing closer to the origin. For further reading on fractal basin boundaries and the difficulties these cause in predicting the long-term behaviour of a dynamical system see: \citet{grebogi1983fractal,grebogi1983finalstatefractal} and \textcite{mcdonald1985fractal}.

The symmetrical properties of the basins discussed in the present section are analysed in Appendix.\,\ref{sssec:SymmUAs}.


\section{Seeing Double Dynamics}\label{sec:SDdynamicsanal}

Sec.\,\ref{sec:SD_results} outlines where multifunctionality is achieved in the $\left( x_{cen}, \, \rho \right)$-plane by assessing the long-term dynamics of the closed-loop RC given by Eq.\,\eqref{eq:PredRes}. However, from Figs.\,\ref{fig:Lerr_oppdir}, \,\ref{fig:Rerr_oppdir}, and \,\ref{fig:BoAs}, it appears that the closed-loop RC must overcome the influence of several untrained attractors before either $\pzc_{A}$ or $\pzc_{B}$ can be reconstructed. From our previous work in \textcite{flynn2021multifunctionality}, it is clear that there are many advantages to treating Eq.\,\eqref{eq:PredRes} like any other dynamical system as by exploring the closed-loop RC's dynamics through a bifurcation analysis further advancements on, for instance, how the RC learns to solve a given problem can be made. In this section, the interplay between the untrained attractors and the reconstructed attractors is investigated in greater detail and the particular bifurcations through which the RC learns how to solve the seeing double problem are identified.

\subsection[Bifurcation Analysis for \texorpdfstring{$x_{cen}=0$}{TEXT}]{Bifurcation Analysis for \texorpdfstring{$x_{cen}=0$}{TEXT}}\label{sec:Bif_xcen0rho}

In Sec.\,\ref{sec:SD_BoAs} multiple fixed points, limit cycles, and a torus were found to coexist in $\mathbb{P}$ at different $\rho$ values. However, according to the results shown in Fig.\,\ref{fig:BoA_rho_125}, these attractors seemingly vanish from $\mathbb{P}$ as multifunctionality is achieved. Moreover, for a small change in $\rho$ there is also a small change in the location of these fixed points and the nature of these limit cycles. This intriguing behaviour warrants further investigation in order to establish the suggested link between these untrained attractors and the reconstructed attractors. 

\begin{figure}
    \centering
    \includegraphics[angle=0,width=0.48\textwidth]{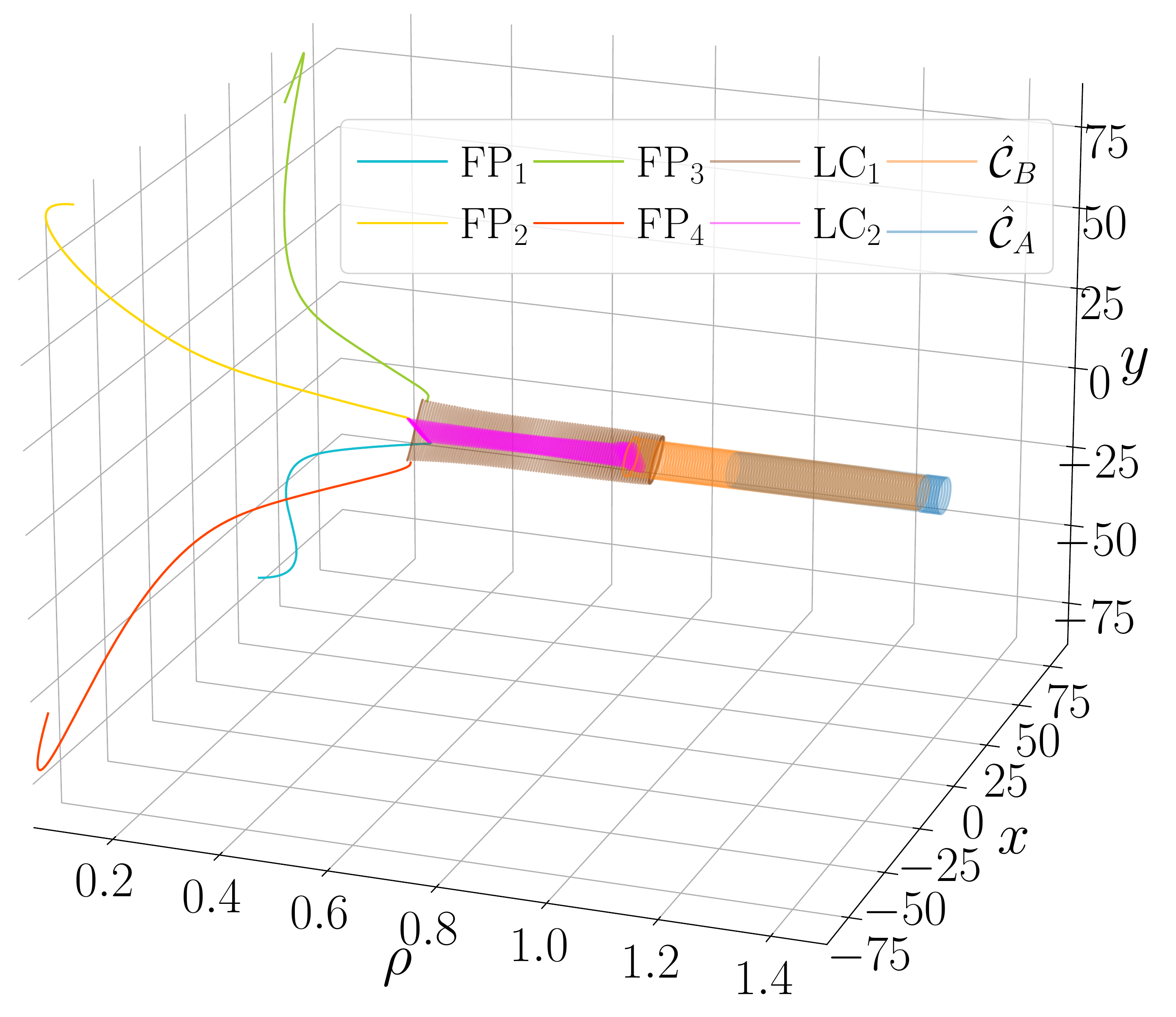}
    \caption{Stable attractors in $\mathbb{P}$ for $\rho \in \left[0.01, 1.47\right]$.}
    \label{fig:rho_bif}
\end{figure}

Taking the above statements into account, the changes in the dynamics of the fixed points, FP$_{1}$, FP$_{2}$, FP$_{3}$, and FP$_{4}$, shown in Figs.\,\ref{fig:BoA_rho_04}-\ref{fig:BoA_rho_03} are now tracked with respect to changes in $\rho$. The same tracking procedure in \textcite{flynn2021multifunctionality} is used here which consists of the following: the changes in the dynamics of a given attractor with respect to, for instance, $\rho$ is tracked by repeating the process of incrementally changing $\rho$, retraining and initialising the state of the closed-loop RC with the attractor corresponding to the previous $\rho$. If there is a change in the dynamics of the closed-loop RC, for instance, a bifurcation from fixed point to limit cycle or a period-doubling (P-D) bifurcation, then the dynamics of the resulting stable attractor is tracked. Fig.\,\ref{fig:rho_bif} shows the results of this.

While the particular bifurcation structure in Fig.\,\ref{fig:rho_bif} may be heavily dependent on the choice of $\textbf{M}$ and $\textbf{W}_{in}$, overall, Fig.\,\ref{fig:rho_bif} provides a clear picture on the challenges the closed-loop RC overcomes in order to achieve multifunctionality. For instance, Fig.\,\ref{fig:rho_bif} shows that for small $\rho$ the closed-loop RC is unable to create any time-varying dynamics and instead produces two pairs of antisymmetric fixed points. As $\rho$ is increased, Fig.\,\ref{fig:rho_bif} establishes that these branches of mirrored fixed points are drawn closer together and two limit cycles are born in quick succession of one another. These are the same limit cycles, LC$_{1}$ and LC$_{2}$, that were found previously in Fig.\,\ref{fig:BoA_rho_08}. 

The role played by these fixed points and limit cycles is akin to the formation of a memory in the brain. As the closed-loop RC learns the correct memory, in this case the correct attractors, along the way its dynamics begins to bear a stronger resemblance to $\hat{\pzc}_{A}$ and $\hat{\pzc}_{B}$ as the influence of the training data grows until attractor reconstruction is achieved. Once $\rho$ is sufficiently large, Fig.\,\ref{fig:rho_bif} shows that the closed-loop RC is first able to reconstruct $\pzc_{B}$ and subsequently for larger $\rho$, the closed-loop RC achieves multifunctionality. A more comprehensive breakdown on each of the major events described above is outlined in the remainder of the present Sec.\,\ref{sec:SDdynamicsanal}.

\subsubsection{Bifurcations of FP$_{1}$ and FP$_{2}$ to LC$_{1}$, and of FP$_{3}$ and FP$_{4}$ to LC$_{2}$}\label{sssec:FP1FP2LC1Bif}

\begin{figure}
    \centering
    \includegraphics[width=0.48\textwidth]{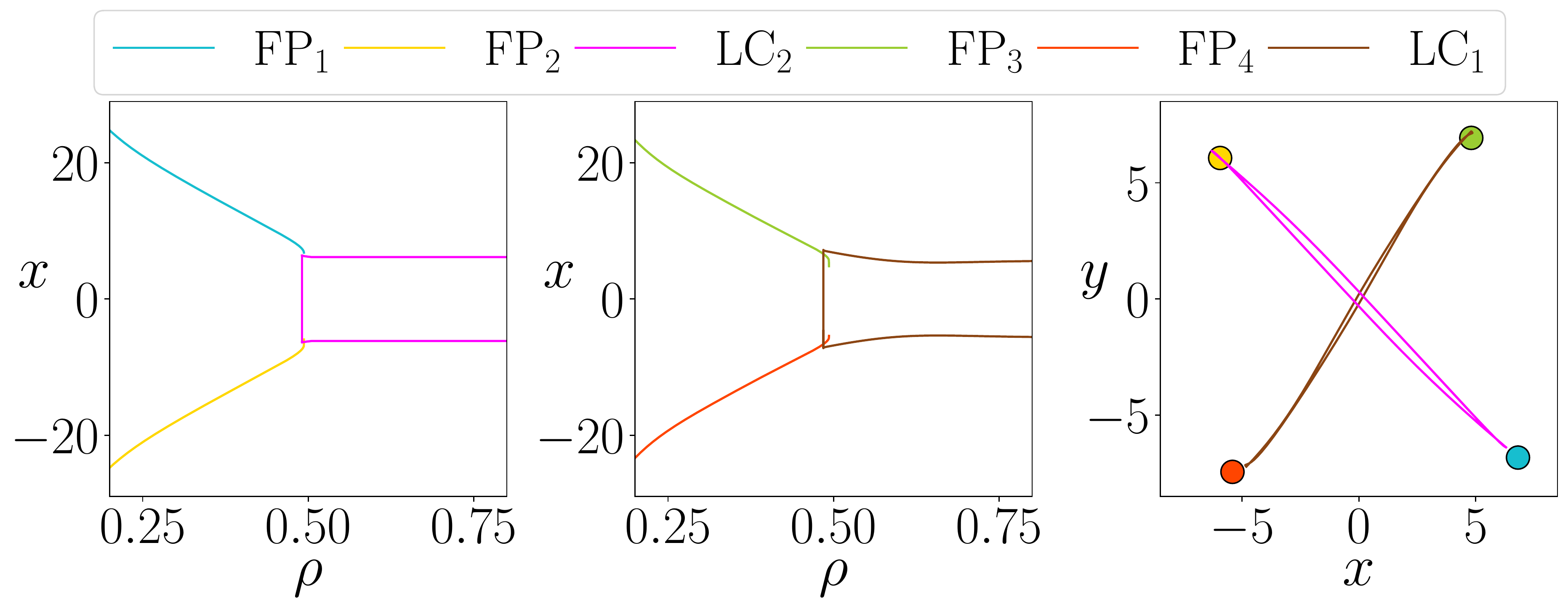}
    \caption{Left and middle panels illustrate the bifurcations from FP$_{1}$ and FP$_{2}$ to LC$_{1}$, and from FP$_{3}$ and FP$_{4}$ to LC$_{2}$, for $\rho \in \left[0.2, 0.8 \right]$ in the $\left( \rho, x \right)$-plane. Right panel shows the location of the fixed points and limit cycles in $\mathbb{P}$ just before they can no longer be tracked.}
    \label{fig:SD_FP_LC_bif}
\end{figure}

Fig.\,\ref{fig:SD_FP_LC_bif} provides a more detailed picture on the bifurcations from FP$_{1}$ and FP$_{2}$ to LC$_{2}$ at $\rho = 0.4935$ in the left hand panel and from FP$_{3}$ and FP$_{4}$ to LC$_{1}$ at $\rho = 0.4931$ in the middle panel. What is plotted here are the location of the fixed points and the local maxima and minima of the limit cycles in the $\left( \rho , x \right)$-plane. When tracking the evolution of these limit cycles for decreasing $\rho$, it is found that LC$_{2}$ can no longer be tracked for $\rho < 0.4905$ and similarly for LC$_{1}$ for $\rho < 0.4847$.

The panel on the right of Fig.\,\ref{fig:SD_FP_LC_bif} displays the location of the fixed points and limit cycles just before and after the respective bifurcations take place on each fixed point and limit cycle. It also appears here that the location of the fixed points just prior to the bifurcation are related to the amplitude of oscillation for both limits cycles after the bifurcation. 

\begin{figure}
    \centering
    \includegraphics[width=0.48\textwidth]{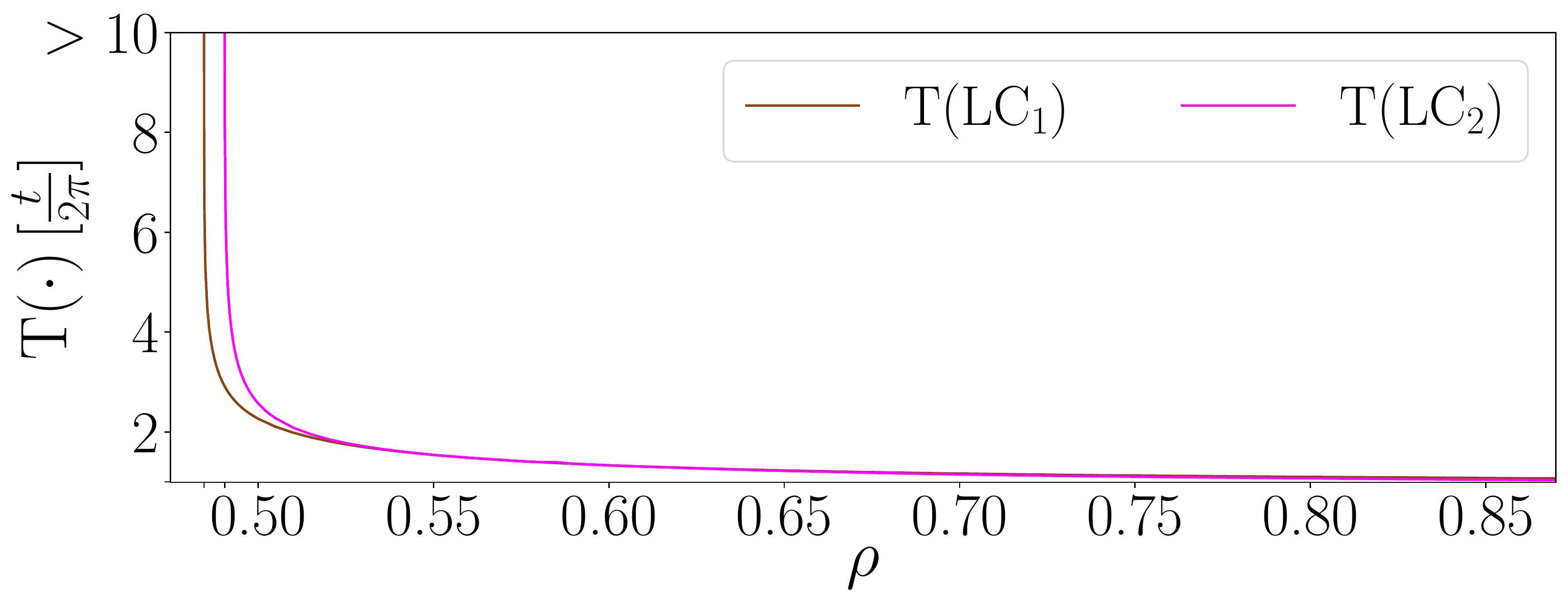}
    \caption{Period, T$\left( \cdot \right)$, of LC$_{1}$ and LC$_{2}$ in units of $\text{t}/2 \pi$ vs. $\rho$.}
    \label{fig:rho_period}
\end{figure}

Fig\,\ref{fig:rho_period} shows that as $\rho$ decreases and approaches the point at which LC$_{1}$ and LC$_{2}$ can no longer be tracked, the period, T$\left( \cdot \right)$, of both limit cycles tends to infinity. T$\left( \cdot \right)$ is plotted here in units of $t/2 \pi$ to provide a comparison to the period of $\pzc_{A}$ and $\pzc_{B}$ which correspond to T$\left( \cdot \right) = 1$.

Fig\,\ref{fig:rho_period} also shows that as $\rho$ increases, T(LC$_{1}$) and T(LC$_{2}$) converge to $1$. This provides further indication that while there does not appear to be a smooth transition from either LC$_{1}$ or LC$_{2}$ to $\hat{\pzc}_{A}$ or $\hat{\pzc}_{B}$, these untrained attractors bear a stronger resemblance to the desired dynamics as $\rho$ increases.

Taking all of the above information on LC$_{1}$ and LC$_{2}$ into account, there are two potential bifurcations which describe how these transitions from two fixed points to one limit cycle takes place. On one hand, since there appears to be a small range of $\rho$ values where FP$_{1}$, FP$_{2}$, and LC$_{2}$ coexist for $\rho \in \left(0.4905, 0.4935\right)$ and similarly where FP$_{3}$, FP$_{4}$, and LC$_{1}$ coexist for $\rho \in \left(0.4847, 0.4931\right)$, then a homoclinic bifurcation may be responsible for the birth of these stable limit cycles. On the other hand, since these windows of multistability appear for only a small range of $\rho$ values, their existence may be an artefact of the tracking procedure used to compute the evolution of these fixed points and limit cycles with regard to changes in $\rho$. If there is no window of multistability then a saddle-node infinite period (SNIPER) bifurcation may instead be responsible for the transition from antisymmetric pairs of fixed points to the corresponding limit cycle.

The above quandary identifies a pitfall of this tracking procedure as it is unable to determine when, for instance, an attractor becomes unstable. Instead this approach relies on a suitably chosen step size of the bifurcation parameter, in this case $\rho$, to account for when a bifurcation takes place. If the step size is sufficiently small then this method can, for example, continue to track unstable solutions of the given dynamical system for a certain range of the bifurcation parameter values.
Despite these issues, the results presented in this section provide significant insight towards the role of untrained attractors, how this closed-loop RC solves the seeing double problem and achieves multifunctionality. 

\subsubsection{Bifurcations of LC$_{1}$ and LC$_{2}$ and the birth of $\hat{\pzc}_{B}$}\label{sssec:Bif_LC1LC2}

\begin{figure}
    \centering
    \includegraphics[width=0.48\textwidth]{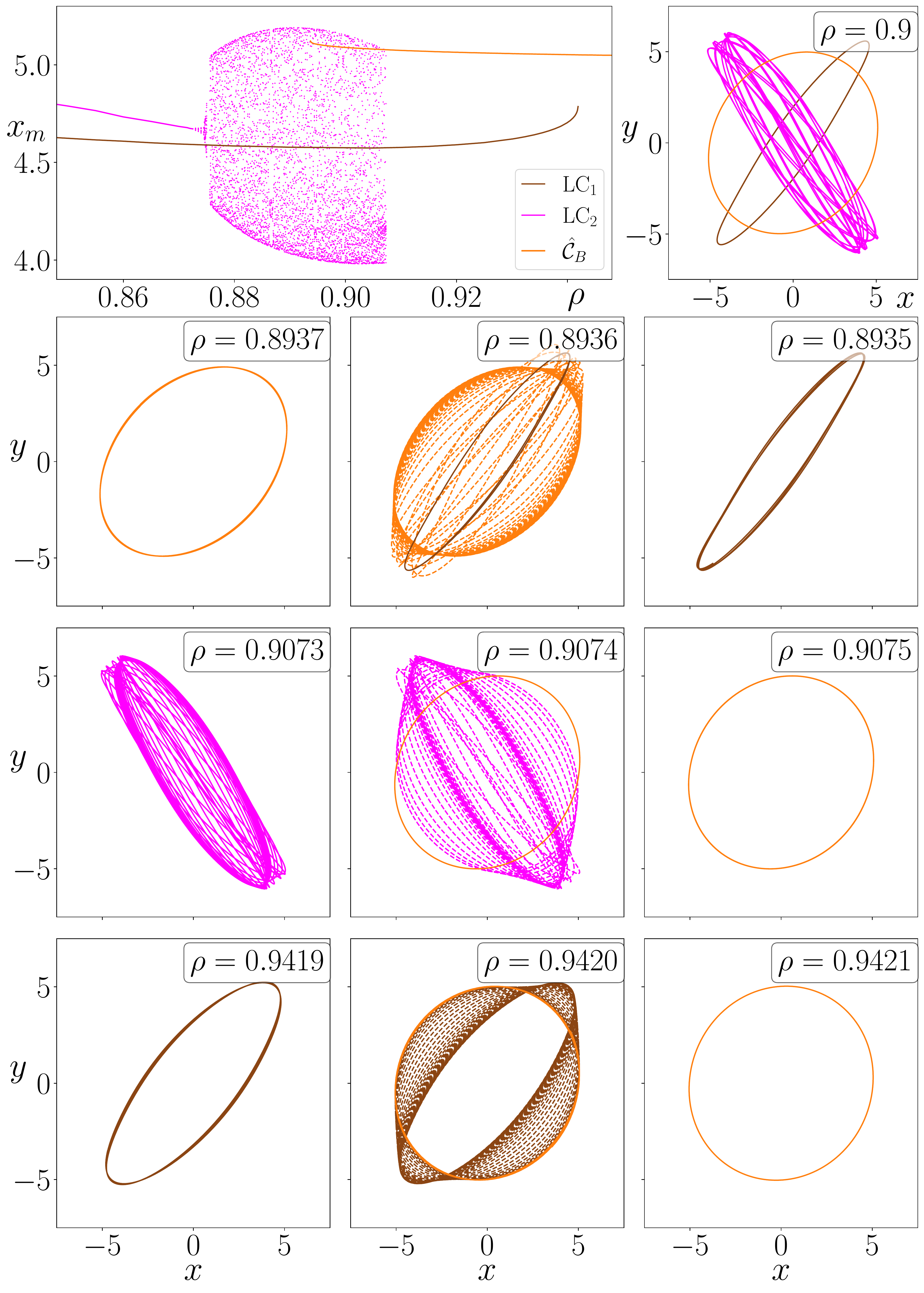}
    \caption{Top left panel shows how the local x maxima of LC$_{1}$, LC$_{2}$, and $\hat{\pzc}_{B}$ change for changes in $\rho$. Top right panel shows the attractors in $\mathbb{P}$ for $\rho = 0.9$. Second, third, and fourth rows show snapshots of $\hat{\pzc}_{B}$, LC$_{2}$, and LC$_{1}$ before they can no longer be tracked.}
    \label{fig:SD_LC1LC2_bif_2}
\end{figure}

Fig.\,\ref{fig:rho_bif} also illustrates how the limit cycles, LC$_{1}$ and LC$_{2}$, evolve in $\mathbb{P}$ with regard to changes in $\rho$. An extensive account of this behaviour and the events which result in the death of these stable limit cycles and the creation of the first reconstructed cycle, $\hat{\pzc}_{B}$, is provided in Fig.\,\ref{fig:SD_LC1LC2_bif_2}.

The evolution of the local maxima of the closed-loop RC's stable solutions in the projected $x$ coordinate (denoted as $x_{m}$) for $\rho \in \left( 0.848,0.948 \right)$ is plotted in the top left panel of Fig.\,\ref{fig:SD_LC1LC2_bif_2}. Here it can be seen that the quasi-periodic torus, T$_{1}$, shown earlier in Fig.\,\ref{fig:BoA_rho_09} comes into existence due to LC$_{2}$ undergoing a torus bifurcation at $\rho = 0.8746$. In light of this, the top right panel provides an updated picture on the stable attractors which, in Fig.\,\ref{fig:BoA_rho_09}, were shown to coexist in $\mathbb{P}$ for $\rho = 0.9$. The colour scheme in the top right panel is consistent with the bifurcation diagram shown in the top left panel instead of assigning the light grey colour to the torus as in Fig.\,\ref{fig:BoA_rho_09}.

The second row of images in Fig.\,\ref{fig:SD_LC1LC2_bif_2} provides additional insight towards the sequence of events which result in the creation/destruction of $\hat{\pzc}_{B}$ as $\rho$ increases/decreases. Moving from left to right, these images show the result of tracking the evolution of $\hat{\pzc}_{B}$ as $\rho$ decreases. The dashed trajectory indicates that $\hat{\pzc}_{B}$ is found to either become unstable or no longer exist between $\rho=0.8937$ and $0.8936$ and the state of the closed-loop RC is on a transient to LC$_{1}$. The closed-loop RC then remains on LC$_{1}$ when $\rho$ is decreased further as illustrated by the image on the right for $\rho = 0.8935$.

The third row of images in Fig.\,\ref{fig:SD_LC1LC2_bif_2} depicts the previously mentioned torus bifurcation and the subsequent death of this torus as $\rho$ increases. By increasing $\rho$ from $0.9073$ to $0.9074$, the torus can no longer be tracked. Following this bifurcation the state of the closed-loop RC departs from the unstable torus (as illustrated by the dashed line) and approaches $\hat{\pzc}_{B}$ and, according to the tracking algorithm, remains on $\hat{\pzc}_{B}$ when $\rho$ is increased as shown in the right panel for $\rho = 0.9075$.

The bottom row of images in Fig.\,\ref{fig:SD_LC1LC2_bif_2} show that by increasing $\rho$ from $0.9419$ to $0.9420$ LC$_{1}$ can no longer be tracked. As $\hat{\pzc}_{B}$ is the only remaining stable solution of the closed-loop RC for $\rho = 0.9420$, the state of the closed-loop RC has no choice but to embark on a transient from the previously stable LC$_{1}$ to $\hat{\pzc}_{B}$ (plotted as the dashed line) and remains on $\hat{\pzc}_{B}$ as $\rho$ is increased to $0.9421$.

\subsubsection{The rise and fall of multifunctionality}\label{sssec:Bif_CbCa}

\begin{figure}
    \centering
    \includegraphics[width=0.48\textwidth]{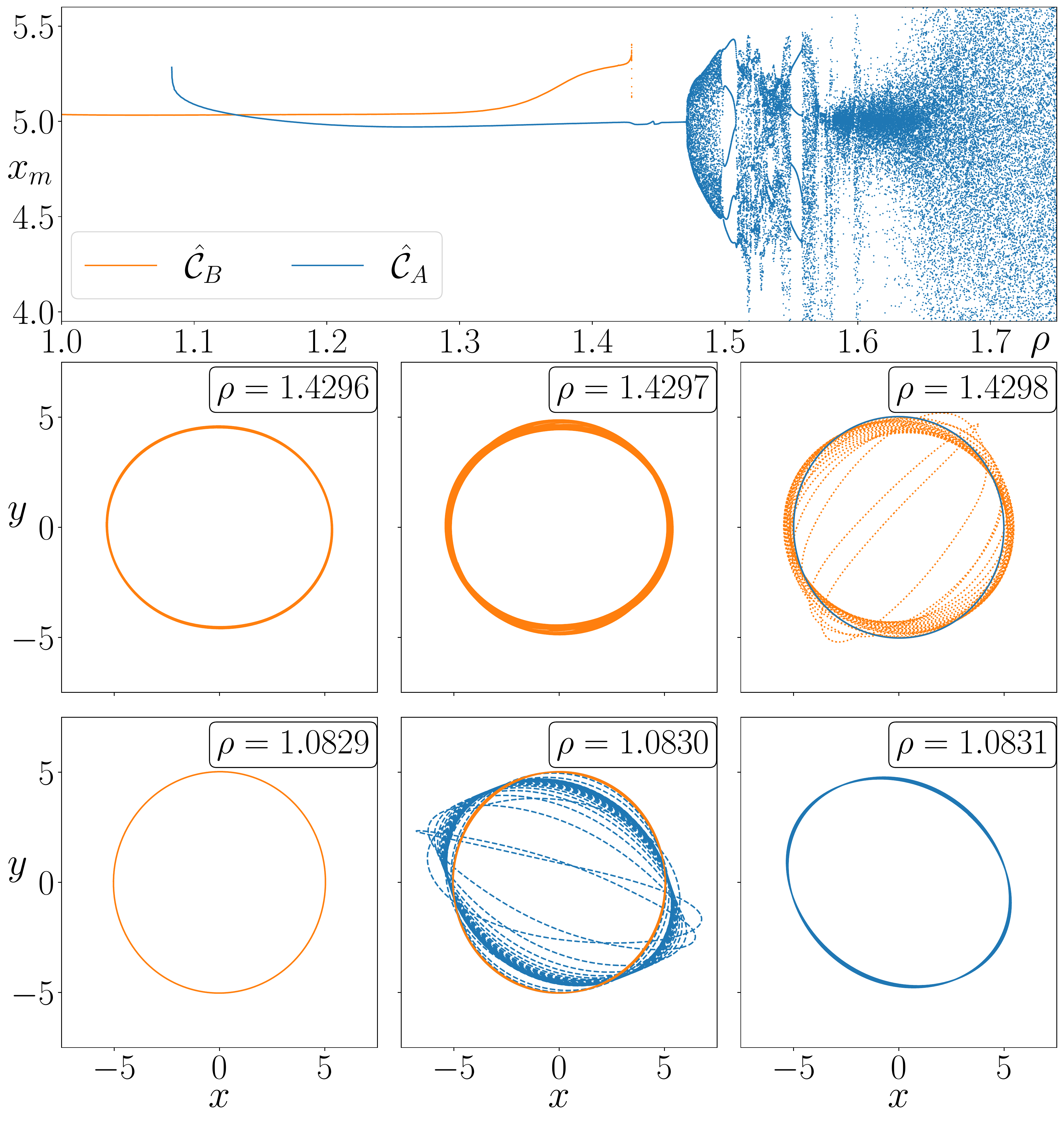}
    \caption{Top panel: shows how the local $x$ maxima of $\hat{\pzc}_{B}$ and $\hat{\pzc}_{A}$ change for changes in $\rho$. Second and third row show snapshots of $\hat{\pzc}_{B}$ and $\hat{\pzc}_{A}$ before they can no longer be tracked.}
    \label{fig:SD_CbCa_bif_2}
\end{figure}

While Fig.\,\ref{fig:rho_bif} highlights the range of $\rho$ values where multifunctionality is achieved, Fig.\,\ref{fig:SD_CbCa_bif_2} expands on this and shows how the rise and fall of multifunctionality begins with the arrival of $\hat{\pzc}_{A}$ at $\rho = 1.0831$ and death of $\hat{\pzc}_{B}$ at $\rho = 1.4297$.

The top panel of Fig.\,\ref{fig:SD_CbCa_bif_2} picks up where Fig.\,\ref{fig:SD_LC1LC2_bif_2} leaves off by continuing to track the changes in the dynamics of $\hat{\pzc}_{B}$ with regard to changes in $\rho$. The second row of images in Fig.\,\ref{fig:SD_CbCa_bif_2} highlights the circumstances in which $\hat{\pzc}_{B}$ can no longer be tracked. By increasing $\rho$ from $1.4296$ to $1.4297$ it is shown here that $\hat{\pzc}_{B}$ no longer has a fixed amplitude of oscillation. Following this, by increasing $\rho$ to $1.4298$, $\hat{\pzc}_{B}$ is no longer a stable limit cycle and the state of the closed-loop RC is instead on a transient to $\hat{\pzc}_{A}$ and remains on $\hat{\pzc}_{A}$ for increasing $\rho$.

The third row of images in Fig.\,\ref{fig:SD_CbCa_bif_2} depicts how $\hat{\pzc}_{A}$ is born. Shown from right to left is the result of tracking the evolution of $\hat{\pzc}_{A}$ for decreasing $\rho$. Here it can be seen that between $\rho = 1.0831$ and $1.0830$, $\hat{\pzc}_{A}$ can no longer be tracked and the state of the closed-loop RC subsequently tends towards $\hat{\pzc}_{B}$ and remains on $\hat{\pzc}_{B}$ for smaller $\rho$ values.

While Fig.\,\ref{fig:Circles_rho_08_17} shows an instance in which the closed-loop RC fails to achieve multifunctionality for $\rho = 1.7$, the top panel of Fig.\,\ref{fig:SD_CbCa_bif_2} reveals the means in which this comes about. Following the death of $\hat{\pzc}_{B}$ at $\rho = 1.4298$, $\hat{\pzc}_{A}$ is the only stable solution present in $\mathbb{P}$. However, $\hat{\pzc}_{A}$ is found to undergo a torus bifurcation at $\rho = 1.4711$. Following this, several windows of periodic and chaotic behaviour appear and the closed-loop RC remains chaotic for $\rho > 1.6$. The first instance of chaos for increasing $\rho$ is in agreement with the corresponding Lyapunov analysis in Fig.\,\ref{fig:SD_MF_regions}. As $\rho$ is increased beyond this point, the range of local maxima in the $x$ variable, as plotted here, continues to grow. As a result, the dynamics of the closed-loop RC loses its resemblance to the original driving input.

Furthermore, the bifurcation diagram presented in Fig.\,\ref{fig:SD_CbCa_bif_2} provides additional information on the results shown previously in Fig.\,\ref{fig:Circles_rho_08_17}. Fig.\,\ref{fig:SD_CbCa_bif_2} shows that the trajectories taken by the closed-loop RC when initialised with either $\boldsymbol{r}_{\left(\pzc_{A}\right)}(t_{\text{listen}})$ or $\boldsymbol{r}_{\left(\pzc_{B}\right)}(t_{\text{listen}})$ for $\rho = 1.7$ are trajectories on the same chaotic attractor. However, we do not claim that the mechanisms responsible for the loss of multifunctionality are universal across all possible $\textbf{M}$ and $\textbf{W}_{in}$ as we find, from further analysis not presented in the present paper, examples where $\hat{\pzc}_{A}$ is first to disappear as $\rho$ increases.

\subsubsection{Floquet analysis}\label{sssec:SD_Floq_anal}

The results presented thus far in Sec.\,\ref{sec:SDdynamicsanal} illustrate some of the hurdles that the closed-loop RC overcomes in order to achieve multifunctionality. However, given the nature of the seeing double problem, it is also possible to provide some quantitative insight into how multifunctionality is achieved by conducting a Floquet analysis. This involves computing the eigenvalues of the `monodromy matrix', $\textbf{Q}$, to classify the stability of a particular attractor using a trajectory of one period on the attractor. $\textbf{Q}$ is the solution of,
\begin{align}
    \dot{\textbf{Q}}(t) = \textbf{J}(t) \textbf{Q}(t), \quad \textbf{Q}(0) = \boldsymbol{I}_{N}, \label{eq:FloqMonodromy}
\end{align}
after one period, $T$. $\textbf{J}(t)$ is the Jacobian matrix, which in this case, is the Jacobian of the closed-loop RC in Eq.\,\eqref{eq:PredRes} and $\textbf{J}(t)$ is computed at each point along a trajectory of one period on the attractor. The eigenvalues of $\textbf{Q}$, denoted by $\mu_{i}$, are called the Floquet multipliers. In the case of the problem at hand, by solving Eq.\,\ref{eq:FloqMonodromy} with input of one period from, for instance, $\hat{\pzc}_{A}$ then the largest multiplier, $\mu_{1}$, will have magnitude $= 1$ and all other $|\mu_{i}| < 1$. While this provides information on when $\hat{\pzc}_{A}$ and $\hat{\pzc}_{B}$ are stable at certain $\rho$ values, the particular bifurcation through which the attractors become stable/unstable can also be determined. Furthermore, dynamical features which the tracking algorithm (used to compute the bifurcation diagrams shown previously) fails to notice may be unearthed by conducting this Floquet analysis. For more information on Floquet theory, see \textcite{nayfeh2008FloquetRef}.

\begin{figure*}
    \centering
    \includegraphics[width=0.86\textwidth]{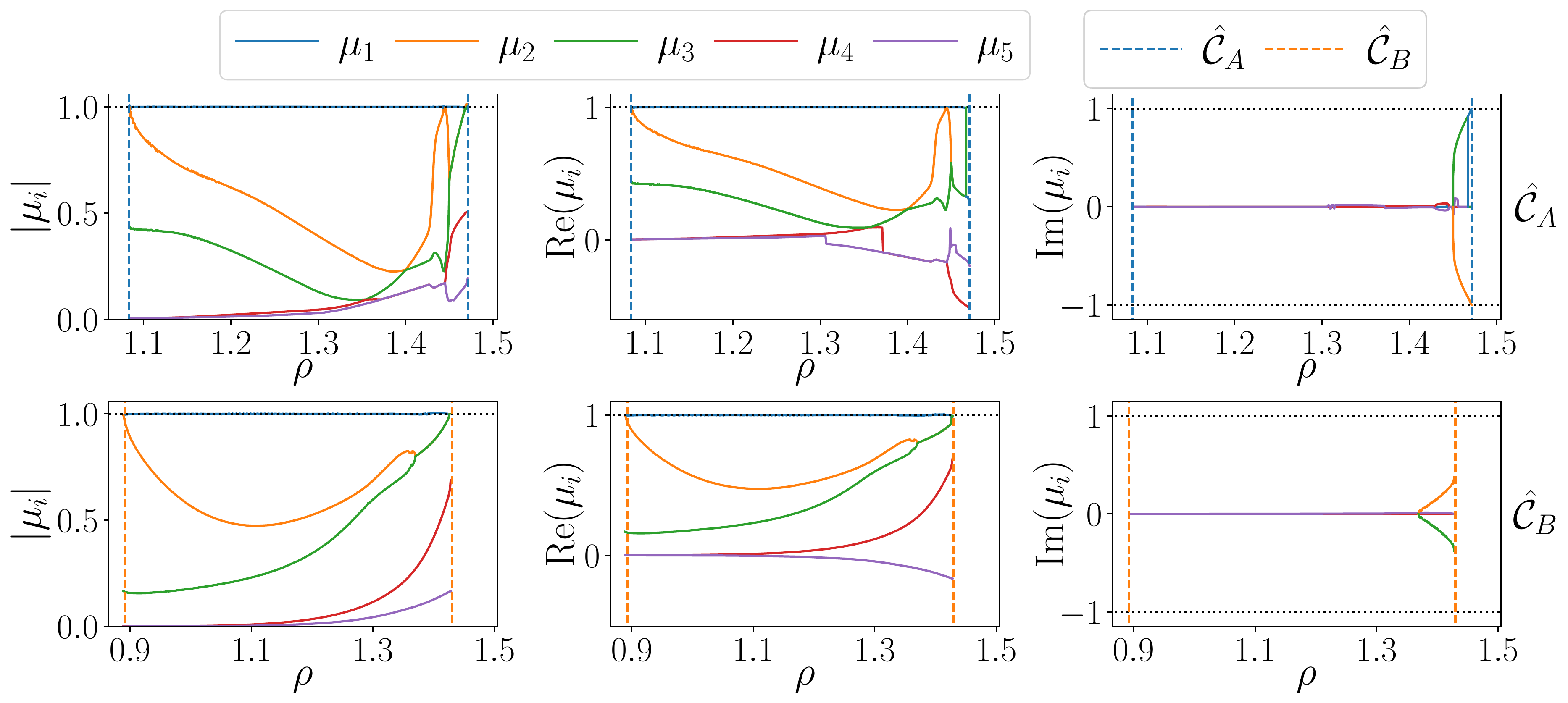}
    \caption{Magnitude, real part, and imaginary part of the first $i=5$ largest Floquet multipliers, $\mu_{i}$, for changes in $\rho$ for $\hat{\pzc}_{A}$ (top row) and $\hat{\pzc}_{B}$ (bottom row). Dashed vertical lines highlight the $\rho$ values where $\hat{\pzc}_{A}$ and $\hat{\pzc}_{B}$ were found to exist between.}
    \label{fig:FloquetCaCb_rho}
\end{figure*}

Taking this into account, Fig.\,\ref{fig:FloquetCaCb_rho} shows how the 5 largest Floquet multipliers evolve for changes in $\rho$ when solving Eq.\,\ref{eq:FloqMonodromy} with drive from $\hat{\pzc}_{A}$ and $\hat{\pzc}_{B}$ in the top and bottom panels respectively. In each case, $| \mu_{1} | = 1$ at all $\rho$ values where $\hat{\pzc}_{A}$ and $\hat{\pzc}_{B}$ were found to exist as stable attractors.

As $\rho$ decreases, Fig.\,\ref{fig:FloquetCaCb_rho} shows that Re$\left( \mu_{2} \right) \to 1$ and since all Im$\left( \mu_{i} \right) = 0$ it can be deduced that $\hat{\pzc}_{A}$ and $\hat{\pzc}_{B}$ are born through a saddle-node bifurcation. For $\hat{\pzc}_{B}$, Fig.\,\ref{fig:FloquetCaCb_rho} demonstrates that Re$\left( \mu_{2} \right)$ and Re$\left( \mu_{3} \right)$ become equal at $\rho = 1.3690$ and beyond this point the corresponding imaginary components of $\mu_{2}$ and $\mu_{3}$ grow according to Im$\left( \mu_{2} \right)$ $= -$ Im$\left( \mu_{3} \right)$ until $|\mu_{2}|$ and $|\mu_{3}|$ tend to $1$ resulting in $\hat{\pzc}_{B}$ becoming unstable at $\rho = 1.4927$. Based on this information and the results shown in Fig.\,\ref{fig:SD_CbCa_bif_2}, Fig.\,\ref{fig:FloquetCaCb_rho} indicates that $\hat{\pzc}_{B}$ becomes unstable through a subcritical torus bifurcation. The behaviour of the Floquet multipliers computed for $\hat{\pzc}_{A}$ in Fig.\,\ref{fig:FloquetCaCb_rho} establishes that a supercritical torus bifurcation is responsible for the $\hat{\pzc}_{A}$'s change in dynamics from stable limit cycle to stable torus at $\rho = 1.4711$ as illustrated in Fig.\,\ref{fig:SD_CbCa_bif_2}. Fig.\,\ref{fig:FloquetCaCb_rho} also reveals further information which Fig.\,\ref{fig:SD_CbCa_bif_2} fails to capture. For instance, a symmetry breaking bifurcation on $\hat{\pzc}_{A}$ is found to occur at $\rho = 1.4450$ as Re$\left( \mu_{2} \right)$ sharply increases to $1$ and quickly decreases from $1$.

\subsection{Bifurcation Analysis in the \texorpdfstring{$\left( x_{cen}, \, \rho \right)$}{TEXT}-plane}\label{ssec:Bif_xcenrho}

Sec.\,\ref{sec:Bif_xcen0rho} provides an extensive account of the particular bifurcations the closed-loop RC goes through in order to achieve multifunctionality when $x_{cen} = 0$. These results suggest that for increasing $\rho$, on one hand, the dynamics of FP$_{1}$, FP$_{2}$, FP$_{3}$, and FP$_{4}$ which lead to the creation of LC$_{1}$ and LC$_{2}$, symbolise a set of challenges for the closed-loop RC to overcome in order to achieve multifunctionality. On the other hand, as $\rho$ increases it can also be argued that the role these fixed points and limit cycles play is to sculpt $\mathbb{P}$ for the arrival of $\hat{\pzc}_{A}$ and $\hat{\pzc}_{B}$. Furthermore, given the ubiquity of the Goldilocks effect across multiple random realisations of $\textbf{M}$ and $\textbf{W}_{in}$ as mentioned in Sec.\,\ref{sec:SD_RegionsMF}, it is worthwhile to study the behaviour of these untrained attractors at different $x_{cen}$ values.



In this section many of the tools applied in Sec.\,\ref{sec:Bif_xcen0rho} are used to explore the dynamics of the closed-loop RC in the $\left( x_{cen},\, \rho \right)$-plane.
The results presented in this section reveal how the relationships between the fixed points and limit cycles mentioned above evolve with respect to changes in $x_{cen}$ nearby where these six attractors coexist for a small range of $\rho$ values when $x_{cen}=0$ in Figs.\,\ref{fig:rho_bif} and \ref{fig:SD_FP_LC_bif}. These results also provide insight on the different roads to $\pzc_{A}$ and $\pzc_{B}$'s reconstruction as conceivably carved out by these fixed points and limit cycles.  

\subsubsection{Tracking stable solutions in the \texorpdfstring{$\left(x_{cen}, \, \rho \right)$}{TEXT}-plane}\label{sssec:Trackingrhoxcen}

\begin{figure*}
    \centering
    \includegraphics[scale=0.34]{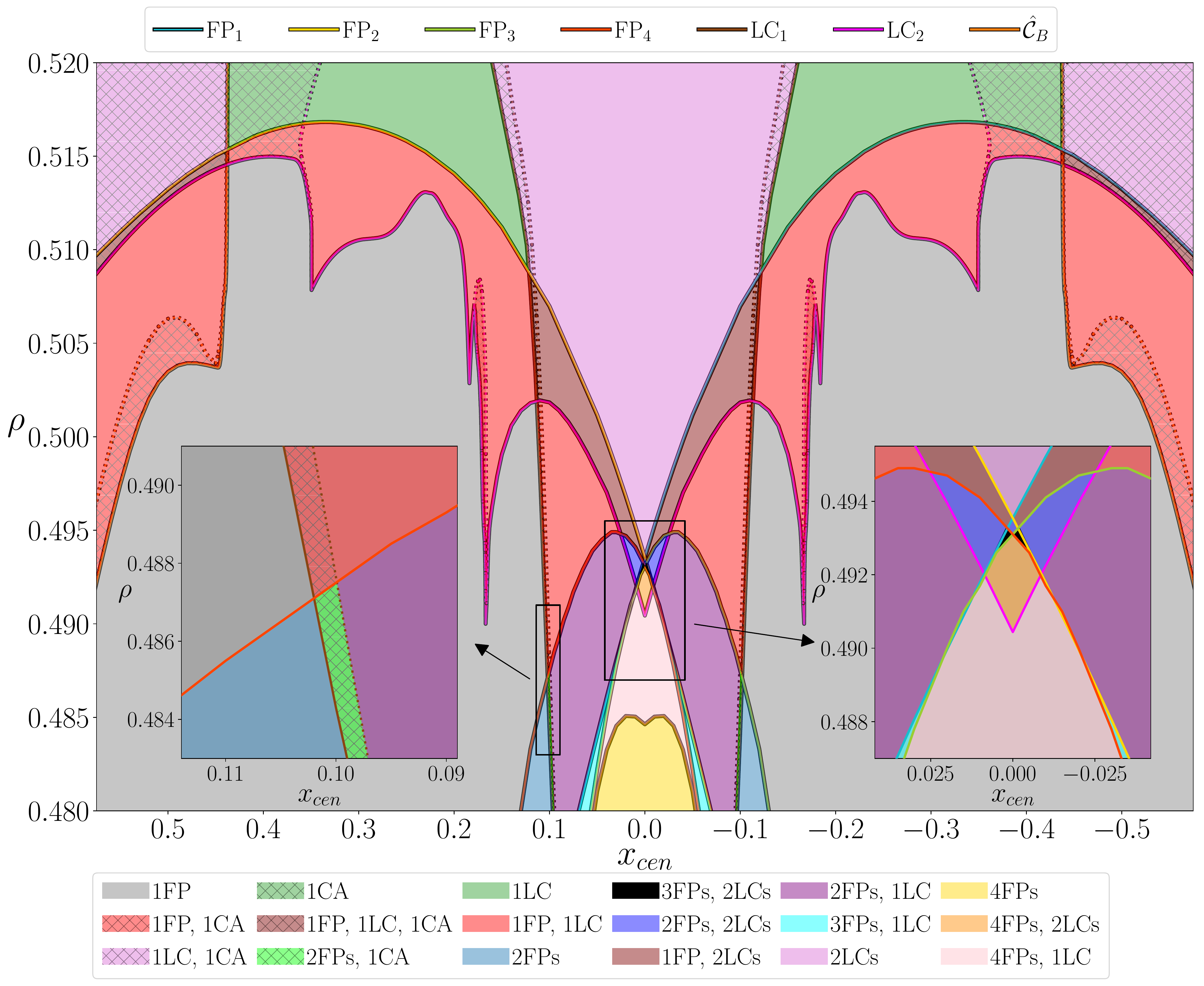}
    \caption{Result of tracking the evolution of stable attractors in $\mathbb{P}$ for $\rho \in \left[0.48,0.52\right]$ and $x_{cen} \in \left[ 0.575, -0.575\right]$. Curves represent the end point of tracking stable solutions. Filled regions describe the nature of the closed-loop RC's stable solutions.}
    \label{fig:rhoxcen_bif}
\end{figure*}

The result of tracking the changes in the dynamics of the closed-loop RC's stable solutions for $\rho \in \left[0.48,0.52\right]$ and $x_{cen} \in \left[ 0.575, -0.575\right]$ are presented in Fig.\,\ref{fig:rhoxcen_bif}. These results are generated by tracking the evolution of LC$_{1}$ and LC$_{2}$ for changes in $x_{cen}$ beginning nearby the bifurcations at $\rho = 0.4847$ and $0.4905$, which spawn LC$_{1}$ and LC$_{2}$ respectively when $x_{cen}=0$, as illustrated in Figs.\,\ref{fig:rho_bif} and \ref{fig:SD_FP_LC_bif}. Similarly, the evolution of FP$_{1}$, FP$_{2}$, FP$_{3}$, and FP$_{4}$ are also tracked with respect to changes in $x_{cen}$ nearby the bifurcations at $\rho = 0.4935$ and $0.4931$ that result in the death of these fixed points. Once these attractors can no longer be tracked for changes in $x_{cen}$, this location in the $\left( x_{cen}, \, \rho \right)$-plane is recorded and added to the corresponding curve in Fig.\,\ref{fig:rhoxcen_bif}. The same process is repeated for incremental changes in $\rho$ for $\rho \in \left[ 0.48, 0.52 \right]$. We did not find any additional attractors in the portion of the $\left( x_{cen}, \, \rho \right)$-plane shown in Fig.\,\ref{fig:rhoxcen_bif} that were unrelated to those studied previously in Sec.\,\ref{sec:Bif_xcen0rho}.

Each curve and filled region in Fig.\,\ref{fig:rhoxcen_bif} has a specific meaning which we outline below.

The labels and colour scheme of each curve in Fig.\,\ref{fig:rhoxcen_bif} are consistent with illustrations of the attractors studied in the previous section. While, for example, the stable period-1 limit cycle, LC$_{1}$, is found to undergo P-D bifurcations at several locations in the $\left( x_{cen},\, \rho \right)$-plane, the resulting stable period-2 limit cycle is also referred to as LC$_{1}$. This convention, which was also taken in relation to LC$_{2}$ and T$_{1}$ in Fig.\,\ref{fig:SD_LC1LC2_bif_2}, is applied even when, for instance, an infinite sequence of P-D bifurcations have taken place when tracking the dynamics of the original period-1 limit cycle. What is not shown nor highlighted along these curves are the continuation of the corresponding codim-1 bifurcation curves and location of codim-2 points associated with the dynamics of certain attractors in the $\left( x_{cen},\, \rho \right)$-plane. This decision was taken in order to avoid overcrowding Fig.\,\ref{fig:rhoxcen_bif} with information that may distract from its main message of tracking how the attractors mentioned above behave in the $\left( x_{cen},\, \rho \right)$-plane. Instead, a more detailed breakdown on some of the bifurcations shown in Fig.\,\ref{fig:rhoxcen_bif} is provided later in the present section.

The colour scheme of the filled regions in the portion of the $\left( x_{cen},\, \rho \right)$-plane shown in Fig.\,\ref{fig:rhoxcen_bif}, as defined by the area the above curves enclose, describe the nature of the closed-loop RC's stable solutions in $\mathbb{P}$ for a given $x_{cen}$ and $\rho$.

\subsubsection{General comments on RC's dynamics in the \texorpdfstring{$\left(x_{cen}, \, \rho \right)$}{TEXT}-plane}\label{sssec:GenCommentsrhoxcen}

The primary message of Fig.\,\ref{fig:rhoxcen_bif} continues the central narrative of this paper, the greater the overlap between $\pzc_{A}$ and $\pzc_{B}$, the closed-loop RC requires larger $\rho$ values to recall more information from the past in order to achieve multifunctionality and the greater the influence of these fixed points and limit cycles until reconstruction occurs. Fig.\,\ref{fig:rhoxcen_bif} further substantiates this statement through examining the behaviour of the stable solutions in $\mathbb{P}$ for the selected $x_{cen}$ and $\rho$ values. What this reveals is the wider role these attractors play in terms of how the closed-loop RC solves the seeing double problem.



One of the main features of Fig.\,\ref{fig:rhoxcen_bif} is the symmetry in the closed-loop RC's dynamics about $x_{cen} = 0$ and this is studied in much greater detail in Appendix\,\ref{sssec:SymmUAs}.

\subsubsection{Behaviour of FP\texorpdfstring{$_{3}$}{TEXT} and FP\texorpdfstring{$_{4}$}{TEXT}}\label{sssec:FP3FP4_rhoxcen}

Concentrating first on FP$_{3}$ and FP$_{4}$, Fig.\,\ref{fig:rhoxcen_bif} reveals that as the magnitude of $x_{cen}$ increases these fixed points play a diminishing role in $\mathbb{P}$. At $\rho = 0.48$, these fixed points cease to coexist for $|x_{cen}| > 0.055$ and both no longer exist for $|x_{cen}| > 0.13$. The inset plot on the right hand side of Fig.\,\ref{fig:rhoxcen_bif} highlights that for increasing $\rho$ along the FP$_{3}$ and FP$_{4}$ curves from $x_{cen} = 0.13$ and $-0.13$ respectively, the FP$_{3}$ and FP$_{4}$ curves become point-like at $x_{cen} \approx -0.035$ and $0.035$ respectively and can no longer be tracked for $\rho > 0.495$.

\subsubsection{Behaviour of FP\texorpdfstring{$_{1}$}{TEXT} and FP\texorpdfstring{$_{2}$}{TEXT}}\label{sssec:FP1FP2_rhoxcen}

Fig.\,\ref{fig:rhoxcen_bif} reveals that as the magnitude of $x_{cen}$ increases, FP$_{1}$ and FP$_{2}$ begin to play a more prominent role in the closed-loop RC's dynamics in comparison to FP$_{3}$ and FP$_{4}$ for the particular $\rho$ values explored here in the $\left( x_{cen},\, \rho \right)$-plane.

While Fig.\,\ref{fig:rhoxcen_bif} shows that the coexistence of FP$_{1}$ and FP$_{2}$ comes to an end at $\rho = 0.48$ for $|x_{cen}| > 0.07$, these fixed points exist far beyond the limits of the $x_{cen}$-axis displayed here. These fixed points occupy the largest combined area in the portion of the $\left( x_{cen},\, \rho \right)$-plane shown here. For increasing $\rho$ along the curves which enclose the regions of the $\left( x_{cen},\, \rho \right)$-plane where FP$_{1}$ and FP$_{2}$ are both stable, the same fate as FP$_{3}$ and FP$_{4}$ awaits FP$_{1}$ and FP$_{2}$ as these curves become point-like at $x_{cen} \approx 0.34$ and $-0.34$ respectively.

As a further comment, the inset plot on the right hand side of Fig.\,\ref{fig:rhoxcen_bif} provides additional insight to the bifurcation diagrams in Figs.\,\ref{fig:rho_bif} and \ref{fig:SD_FP_LC_bif} and draws attention to how the coexistence between FP$_{1}$ and FP$_{2}$, and FP$_{3}$ and FP$_{4}$ at $x_{cen}=0$ unfolds with respect to changes in $x_{cen}$.

\subsubsection{Behaviour of LC\texorpdfstring{$_{1}$}{TEXT}}\label{sssec:LC1_rhoxcen}

Focusing next on LC$_{1}$, Fig.\,\ref{fig:rhoxcen_bif} illustrates the limited role this limit cycle plays in the closed-loop RC's dynamics as the magnitude of $x_{cen}$ increases from $0$. For decreasing $\rho$, the range of $x_{cen}$ values where LC$_{1}$ is stable also decreases. By increasing $|x_{cen}|$ from $0$ for $\rho = 0.48$, a region in the $\left( x_{cen},\, \rho \right)$-plane where all four fixed points coexist, LC$_{1}$ is first encountered as a stable solution at $x_{cen} \approx \pm 0.05$. LC$_{1}$ briefly coexists with these four fixed points until FP$_{3}$ and FP$_{4}$ can no longer be tracked as previously mentioned at $x_{cen} = \mp 0.055$ and likewise for FP$_{1}$ and FP$_{2}$ at $x_{cen} = \pm 0.07$. This leaves LC$_{1}$ in coexistence with FP$_{2}$ and FP$_{3}$ for $x_{cen} \in \left[ 0.071, 0.097 \right]$ and in coexistence with FP$_{1}$ and FP$_{4}$ for $x_{cen} \in \left[ -0.071, -0.097 \right]$. LC$_{1}$ can no longer be tracked for $|x_{cen}| > 0.097$ at $\rho = 0.48$. 

As a further comment, as $\rho$ decreases, the brown curves on either side of $x_{cen} = 0$, which outline where LC$_{1}$ is stable, appear to be converging to a point along the FP$_{3}$ and FP$_{4}$ curves. This suggests that a codim-2 bifurcation takes place in the $\left( x_{cen},\, \rho \right)$-plane outside of the regions explored in Fig.\,\ref{fig:rhoxcen_bif}.

Prior to no longer being able to track LC$_{1}$ for increasing the magnitude of $x_{cen}$, Fig.\,\ref{fig:rhoxcen_bif} also illustrates that LC$_{1}$ begins to exhibit chaos along the dotted brown curve $\forall \, \rho$ values shown here. The extent of this chaos in the $\left( x_{cen},\, \rho \right)$-plane is described by the area enclosed by this dotted curve and the LC$_{1}$ curve, and is further emphasised by the hatched pattern. It is shown here that, on one hand, these windows of chaos widen as $\rho$ increases, on the other hand, the inset plot on the left hand side of Fig.\,\ref{fig:rhoxcen_bif} provides a more detailed picture on how, this small region of chaos narrows as $\rho$ decreases.

\begin{figure}
    \centering
    \includegraphics[width=0.48\textwidth]{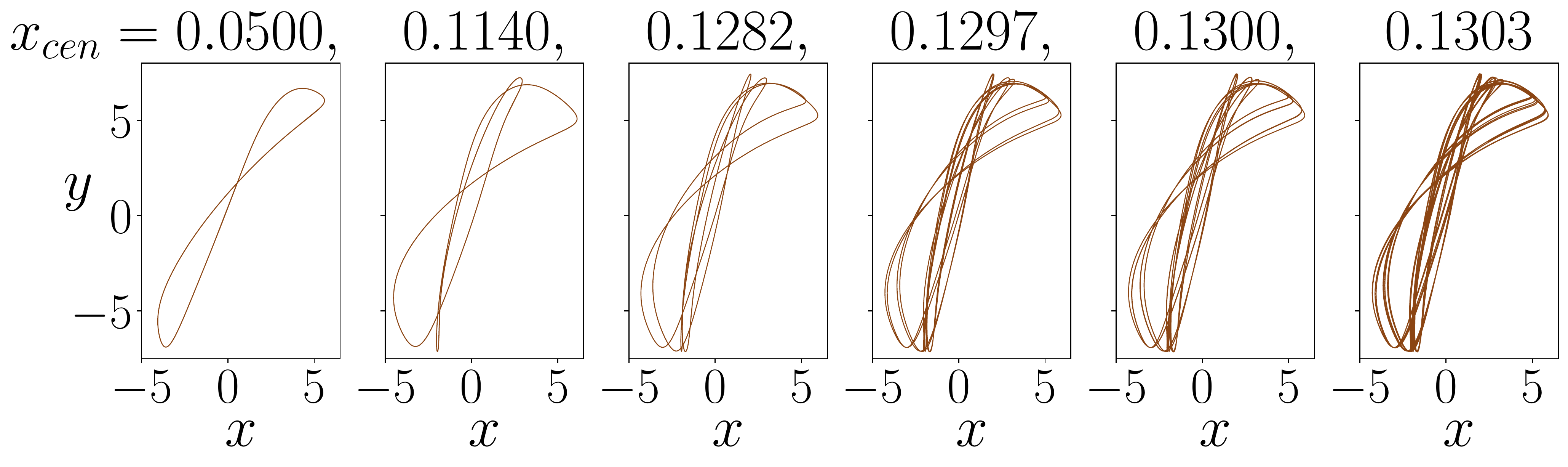}
    \caption{Evolution of LC$_{1}$ in $\mathbb{P}$ for changes in $x_{cen}$ with $\rho = 0.515$.}
    \label{fig:LC1_pdc}
\end{figure}

Fig.\,\ref{fig:LC1_pdc} sheds some light on LC$_{1}$'s change from periodic to chaotic dynamics.
For $\rho = 0.515$, it is shown here that as $x_{cen}$ increases, LC$_{1}$ undergoes a series of P-D bifurcations. An infinite sequence of these bifurcations occurs between $x_{cen} = 0.1301$ and $0.1302$ which gives rise to the chaotic dynamics shown in Fig.\,\ref{fig:LC1_pdc} for $x_{cen}=0.1303$.

\subsubsection{Behaviour of LC\texorpdfstring{$_{2}$}{TEXT}}\label{sssec:LC2_rhoxcen}

Fig.\,\ref{fig:rhoxcen_bif} lays bare the intricate relationship between LC$_{2}$, $\rho$, and $x_{cen}$. Several events occur along the LC$_{2}$ curve as the magnitude of $x_{cen}$ increases. The inset plot on the right hand side highlights the first of these incidents with the loss of its coexistence with the four fixed points and LC$_{1}$. This continues with larger $|x_{cen}|$ values and at $\rho = 0.502$, LC$_{2}$ is in coexistence with only FP$_{1}$ for $x_{cen} < -0.115$ and FP$_{2}$ for $x_{cen} > 0.115$ until $\pzc_{B}$ is reconstructed at larger $x_{cen}$ values.

It is also shown in Fig.\,\ref{fig:rhoxcen_bif} that nearby these points in the $\left( x_{cen},\, \rho \right)$-plane where LC$_{2}$ loses its coexistence with LC$_{1}$, the LC$_{2}$ curve itself reaches its first of several local maxima. In order to continue tracking LC$_{2}$ beyond this point, $\rho$ is required to decrease as $|x_{cen}|$ increases until the curve reaches its first of several local minima at $\left( x_{cen},\, \rho \right) = \left( \pm 0.16, 0.49 \right)$ where a P-D bifurcation occurs. Subsequent P-D bifurcations are encountered along the LC$_{2}$ curve for slightly increasing $|x_{cen}|$ and $\rho$. The hatched region of the $\left( x, \rho \right)$-plane enclosed by the LC$_{2}$ curve and the nearby dotted curve characterise the small window of chaotic solutions produced by an infinite sequence of these P-D bifurcations. 

\begin{figure}
    \centering
    \begin{subfigure}{0.48\textwidth}
    \includegraphics[width=\textwidth]{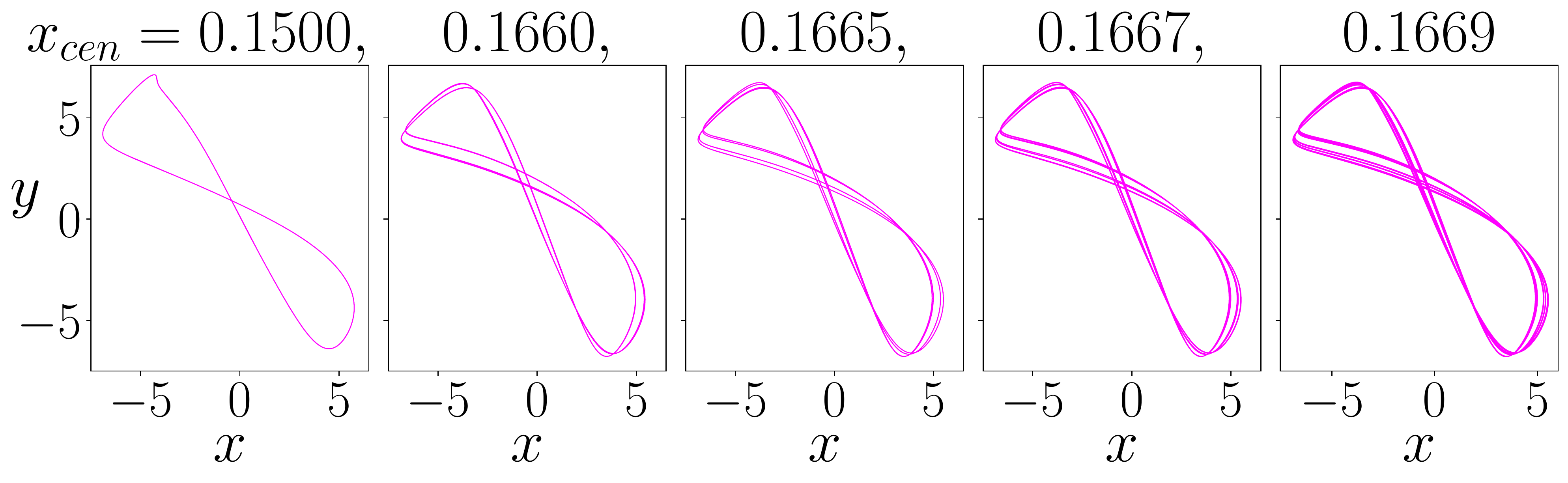}
    \caption{$\rho = 0.501$}
    \label{fig:LC2_rho501pdc}
    \end{subfigure}
    \hfill
    \begin{subfigure}{0.48\textwidth}
    \includegraphics[width=\textwidth]{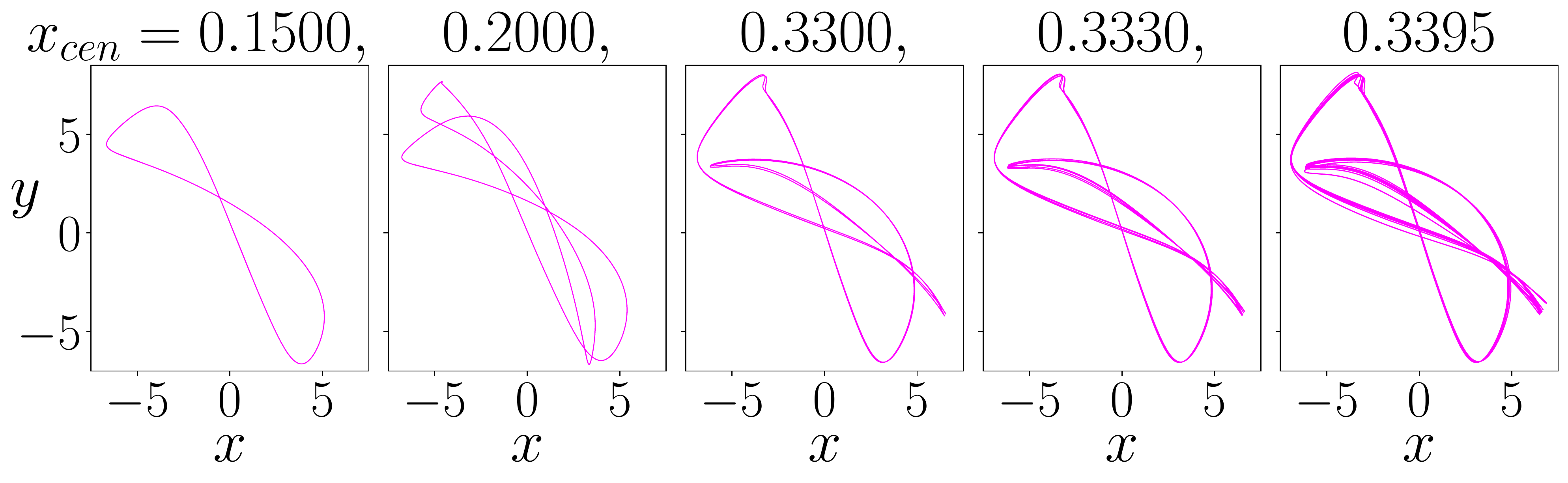}
    \caption{$\rho = 0.511$}
    \label{fig:LC2_rho511pdc}
    \end{subfigure}
    \caption{Evolution of LC$_{2}$ in $\mathbb{P}$ for changes in $x_{cen}$ with $\rho = 0.501$ in (a) and $0.511$ in (b).}
    \label{fig:LC2_rho501511pdc}
\end{figure}

Fig.\,\ref{fig:LC2_rho501pdc} provides an example of the closed-loop RC's dynamics nearby one of these P-D cascades for $\rho = 0.501$. It is shown here that when tracking the evolution of the original stable period-1 limit cycle, successive P-D bifurcations are encountered for increasing $x_{cen}$. Fig.\,\ref{fig:LC2_rho501pdc} shows examples of the corresponding period-1 dynamics at $x_{cen}=0.1500$, period-2 dynamics at $x_{cen}=0.1660$, and period-4 dynamics at $x_{cen}=0.1665$. In this scenario the period-8 limit cycle, whose dynamics are shown here for $x_{cen} = 0.1667$, undergoes an infinite sequence of P-D bifurcations for a small increase in $x_{cen}$. This produces the chaotic dynamics shown here for $x_{cen} = 0.1669$. This region of chaos in Fig.\,\ref{fig:rhoxcen_bif} ends once the LC$_{2}$ curve and dotted curve meet at $\left( x_{cen}, \, \rho \right) = \left( 0.178, 0.507 \right)$. For increasing $|x_{cen}|$ beyond this point, it is found that LC$_{2}$ is a period-2 limit cycle before it can no longer be tracked. However, subsequent P-D bifurcations occur along the LC$_{2}$ curve for increasing $|x_{cen}|$. 

Following an infinite sequence of these P-D bifurcations, what Fig.\,\ref{fig:rhoxcen_bif} also shows is that along the LC$_{2}$ curve, chaos reappears at $\left( x_{cen},\, \rho \right) = \left( \pm 0.35, 0.508 \right)$. The dotted magenta curve which stems from this point distinguishes between when LC$_{2}$ is periodic or chaotic as further emphasised by the hatched pattern. 
An example of this particular route to chaos is provided in Fig.\,\ref{fig:LC2_rho511pdc}. Furthermore, although it is not shown here, these P-D bifurcation do not only occur along the LC$_{2}$ curve, there is a locus which connects each given set of P-D bifurcations points in the $\left( x_{cen},\, \rho \right)$-plane in a similar way to the dotted magenta curve mentioned above. 

What Fig.\,\ref{fig:LC2_rho511pdc} depicts is the result of tracking the evolution of LC$_{2}$'s dynamics for changes in $x_{cen}$ when $\rho = 0.511$. The remnants of the previous small window of chaos, which occurs after the first local minima along the LC$_{2}$ curve, are visible here with a period-2 limit cycle found for $x_{cen}=0.2$. The major difference between this period-2 limit cycle and the example illustrated in Fig\,\ref{fig:LC2_rho501pdc} for $\rho = 0.501$ is the prominence of the sharp turn around the max $y$ value on the period-2 limit cycle in Fig.\,\ref{fig:LC2_rho511pdc}. This remains a distinctive feature of LC$_{2}$ when tracking its dynamics as $x_{cen}$ is further increased and persists after each of the subsequent P-D bifurcations that result in chaos. It is also worth noting that this kink along LC$_{2}$'s dynamics occurs nearby the location of FP$_{2}$ and FP$_{1}$ for the corresponding positive and negative $x_{cen}$ values respectively.

LC$_{2}$ remains chaotic as $|x_{cen}|$ increases from $0.35$ and, depending on the sign of $x_{cen}$, continues its coexistence with either FP$_{1}$ or FP$_{2}$ for the range of $x_{cen}$ values shown in Fig.\,\ref{fig:rhoxcen_bif}. Much like LC$_{1}$ with FP$_{3}$ and FP$_{4}$, Fig.\,\ref{fig:rhoxcen_bif} also suggests that the LC$_{1}$ curve may connect with the FP$_{1}$ and FP$_{2}$ curves beyond the range of $x_{cen}$ values shown here.

\begin{figure}
    \centering
    \begin{subfigure}{0.48\textwidth}
    \includegraphics[width=\textwidth]{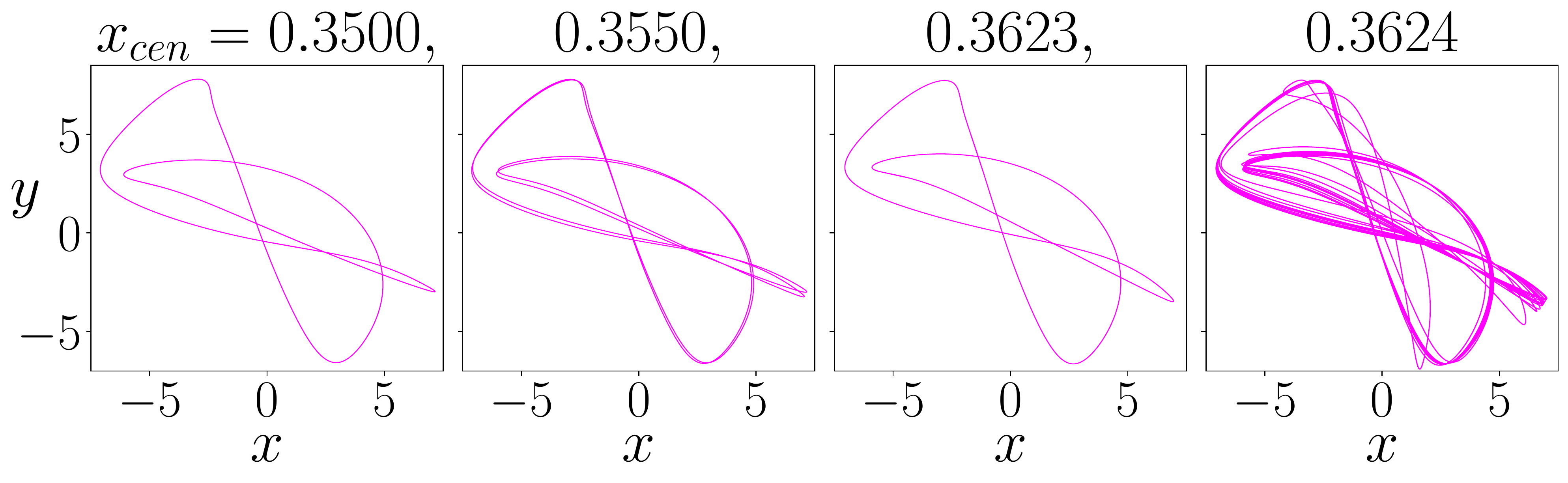}
    \caption{$\rho = 0.515$}
    \label{fig:LC2_rho515pdc}
    \end{subfigure}
    \hfill
    \begin{subfigure}{0.48\textwidth}
    \includegraphics[width=\textwidth]{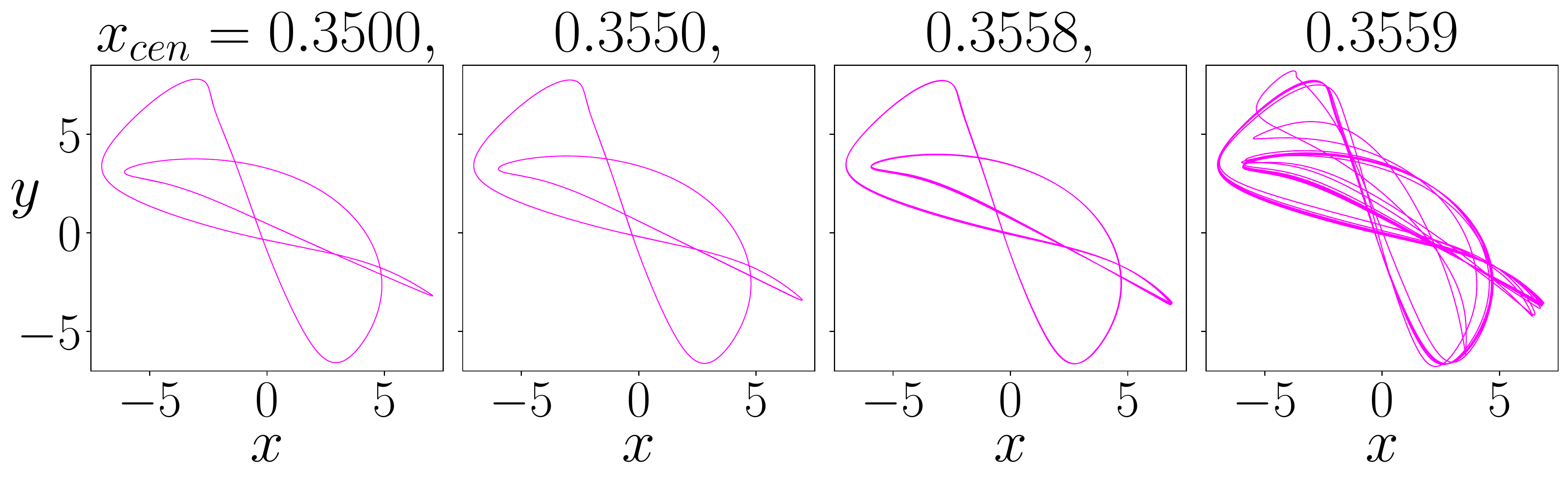}
    \caption{$\rho = 0.517$}
    \label{fig:LC2_rho517pdc}
    \end{subfigure}
    \caption{Evolution of LC$_{2}$ in $\mathbb{P}$ for changes in $x_{cen}$ with $\rho = 0.515$ in (a) and $\rho = 0.517$ in (b).}
    \label{fig:rho_515_517}
\end{figure}

As $\rho$ increases along this dotted curve of chaos, which begins at $\left( x_{cen},\, \rho \right) = \left( \pm 0.35, 0.508 \right)$, there is a point where these sequences of P-D bifurcations are no longer directly responsible for the onset of chaotic dynamics. For instance, Fig.\,\ref{fig:LC2_rho515pdc} shows that when tracking the evolution of LC$_{2}$ for $\rho = 0.515$, by increasing $x_{cen}$ from $0.3623$ to $0.3624$, there is a sudden change in dynamics from a period-2 limit cycle to a chaotic attractor. The same is shown Fig.\,\ref{fig:LC2_rho517pdc} for $\rho = 0.517$ and $x_{cen}$ increases from $0.3558$ to $0.3559$. The only difference between these examples of LC$_{2}$'s dynamics prior to the dawning of chaos is the lack of a period-4 limit cycle for $\rho = 0.517$.

\subsubsection{Behaviour of \texorpdfstring{$\hat{\pzc}_{B}$}{TEXT}}\label{sssec:Cb_rhoxcen}

Fig.\,\ref{fig:rhoxcen_bif} also illustrates regions in the $\left( x_{cen},\, \rho \right)$-plane where $\hat{\pzc}_{B}$ comes into existence. For much of its presence in Fig.\,\ref{fig:rhoxcen_bif}, $\hat{\pzc}_{B}$ is in coexistence with either, LC$_{2}$, or one of FP$_{2}$ or FP$_{1}$ depending on the sign of $x_{cen}$, or LC$_{2}$ and FP$_{2}$ or FP$_{1}$.

At the limits of the $x_{cen}$-axis in Fig.\,\ref{fig:rhoxcen_bif} (when $|x_{cen}| = 0.575$), the nature of $\hat{\pzc}_{B}$ differs greatly depending on $\rho$. In particular, the dotted orange curve describes the boundary between when $\hat{\pzc}_{B}$ exhibits periodic or chaotic behaviour, and it is within the hatched region that $\hat{\pzc}_{B}$ is chaotic. For decreasing $|x_{cen}|$, these two curves approach one another and the window of chaos becomes arbitrarily small at $|x_{cen}| = 0.45$. By slightly decreasing $|x_{cen}|$ from this point in the $\left( x_{cen},\, \rho \right)$-plane a sudden increase in $\rho$ is required in order to reconstruct $\hat{\pzc}_{B}$.

\begin{figure*}
    \centering
    \includegraphics[width=0.75\textwidth]{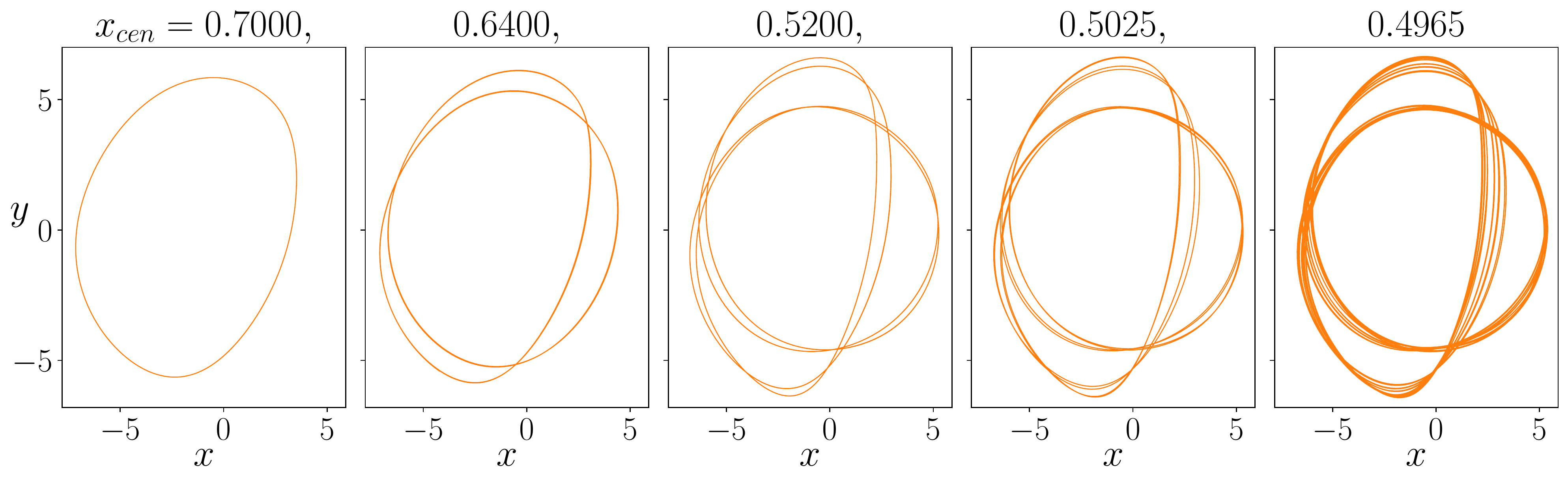}
    \caption{Evolution of $\hat{\pzc}_{B}$ in $\mathbb{P}$ for changes in $x_{cen}$ with $\rho = 0.504$.}
    \label{fig:Cb_pdc}
\end{figure*}

Fig.\,\ref{fig:Cb_pdc} illustrates the changes in $\hat{\pzc}_{B}$'s dynamics when tracking its evolution for decreasing $x_{cen}$ when $\rho = 0.504$. Here it can be seen that $\hat{\pzc}_{B}$ is a period-1 limit cycle for $x_{cen} = 0.7000$ and is also not perfectly circular. $\hat{\pzc}_{B}$ undergoes several P-D bifurcations with examples of period-2, period-4, and period-8 dynamics shown here for $x_{cen} = 0.6400, 0.5200,$ and $0.5205$. An infinite sequence of P-D bifurcations occurs as $x_{cen}$ is decreased further resulting in the chaotic dynamics like those shown here for $x_{cen} = 0.4965$. 

Above all, what Figs.\,\ref{fig:Cb_pdc} and \ref{fig:rhoxcen_bif} show is that while $\hat{\pzc}_{B}$ is stable its dynamics do not always resemble $\pzc_{B}$ even to the point when $\hat{\pzc}_{B}$ is not in any way circular, has a much larger period, or even exhibits chaos. At the same time, it is also worth noting that after a number of P-D bifurcations, for instance the period-4 limit cycle shown in Figs.\,\ref{fig:Cb_pdc}, portions of the trajectory on $\hat{\pzc}_{B}$ begin to resemble a more circular orbit with similar characteristics to $\pzc_{B}$. This suggests that these P-D bifurcations may play an important role in how the closed-loop RC improves its reconstruction of $\pzc_{B}$.

\section{\label{sec:DisConc}Conclusion}

The results presented here show that for the closed-loop RC described in Eq.\,\eqref{eq:PredRes} to achieve multifunctionality there is a crucial dependence on $\rho$, the spectral radius of the RC's internal connections. The intricacies of this relationship are revealed in Sec.\,\ref{sec:SD_RegionsMF} when the cycles $\pzc_{A}$ and $\pzc_{B}$ are moved closer together. As $\pzc_{A}$ and $\pzc_{B}$ begin to overlap, Fig.\,\ref{fig:SD_MF_regions} shows that the closed-loop RC requires a greater amount of memory in order to successfully distinguish between the orbits. A `Goldilocks' effect is also found whereby if $\rho$ is too small or too large then multifunctionality will not occur.

From our results we have obtained a greater understanding of multifunctionality in a RC. In the context of the formalism developed in Sec.\,\ref{sec:trainingMFRC}, we have shown that while, for instance, $\pza_{1} \cap \pza_{2} \neq \emptyset$ (there can be an overlap in the training data), for a closed-loop RC to achieve multifunctionality, we require that the open-loop RC responds uniquely to each of the driving inputs $\boldsymbol{u}_{\left(\pza_{1}\right)}(t)$ and $\boldsymbol{u}_{\left(\pza_{2}\right)}(t)$ so that at no time $t$ during the training $\boldsymbol{r}_{\left(\pza_{1}\right)}(t) = \boldsymbol{r}_{\left(\pza_{2}\right)}(t)$ even if $\boldsymbol{u}_{\left(\pza_{1}\right)}(t) = \boldsymbol{u}_{\left(\pza_{2}\right)}(t)$.

While Figs.\,\ref{fig:SD_MF_regions} and \ref{fig:SD_MF_regions_samedir} show that the closed-loop RC achieves its best performance just prior to exhibiting chaotic dynamics, Fig.\,\ref{fig:CI_example} provides some evidence that beyond this transition to chaos the closed-loop RC enters into the dynamical regime of chaotic itinerancy. From a different point of view, multifunctionality can be used as a precursor to generate a specifically designed chaotically itinerant orbit between several attractor ruins thereby providing a suitable means to realise Tsuda's thoughts on the role of `chaotic itinerancy' in cognitive neurodynamics in \textcite{tsuda2015chaotic}.
Tsuda puts forward chaotic itinerancy as a potential dynamical mechanism to describe how the brain transitions between performing different tasks or recounting different memories characterised by each attractor ruin. We intend to study this connection between multifunctionality and chaotic itinerancy and how can multifunctionality be exploited as a data-driven modelling tool in future work.

Sec.\,\ref{sec:SD_BoAs} focuses on how the closed-loop RC constructs basins of attraction from the perspective of $\mathbb{P}$ for $\hat{\pzc}_{A}$ and $\hat{\pzc}_{B}$ when they are entirely overlapping, i.e. when $x_{cen}=0$. The nature of $\mathbb{P}$ prior to the closed-loop RC achieving multifunctionality is also explored. Fig.\,\ref{fig:BoAs} unveils the coexistence of several fixed points with fractal basin boundaries, limit cycles and tori in $\mathbb{P}$, these are untrained attractors which were not involved in the training. This discovery provides the basis for conducting the bifurcation analyses carried out in Sec.\,\ref{sec:SDdynamicsanal} in order to determine the role these attractors play in how the closed-loop RC solves the seeing double problem. 

This bifurcation analysis is divided into two studies, the results of the first are outlined in Sec.\,\ref{sec:Bif_xcen0rho} which explores how changes in $\rho$ effects the closed-loop RC's dynamics in the case of $x_{cen} = 0$. The results of the second study are presented in Sec.\,\ref{ssec:Bif_xcenrho} where the closed-loop RC's dynamics in a portion of the $\left( x_{cen}, \, \rho \right)$-plane are assessed. Here the role of the untrained attractors in the learning dynamics of the closed-loop RC becomes apparent. The bifurcations which result in the rise and fall of multifunctionality are identified through the results of the Floquet analysis shown in Fig.\,\ref{fig:FloquetCaCb_rho}. The main take home message from both of these bifurcation analyses is the same, the more $\pzc_{A}$ and $\pzc_{B}$ overlap, the more memory the the closed-loop RC requires in order to exhibit multifunctionality and the greater the influence of these untrained attractors until reconstruction occurs.

Even though the bifurcation analysis tools used in this paper have provided significant insight into how the closed-loop RC learns, performs certain tasks, and interacts with untrained attractors, it is likely that only half the story is being told as the unstable solutions are not accounted for. 
As mentioned in \textcite{flynn2021multifunctionality}, these issues can be alleviated by adapting well-known continuation software like AUTO \cite{AUTO_Doedel1} to the case of reservoir computing.

While the bifurcation analysis in Sec.\,\ref{sec:SDdynamicsanal} is conducted for only one random realisation of $\textbf{M}$ and $\textbf{W}_{in}$, it is also possible that multifunctionality may be achieved in various other scenarios for different choices of $\textbf{M}$ and $\textbf{W}_{in}$. For instance, \textcite{morraflynn23_MF_fly} investigate the dynamics of the same RC design as in Eq.\,\ref{eq:ListenRes}, when trained to solve the seeing double problem in the case of $x_{cen}=0$, using the same technique that produced Fig.\,\ref{fig:rho_bif}. It is shown that for a much smaller number of neurons ($N=500$ as opposed to $1000$) there is a greater interaction and coexistence between the trained attractors and the untrained attractors in comparison to the results presented in Fig.\,\ref{fig:rho_bif}. This may indicate that as $N$ is increased, the less the untrained attractors interact and coexist with the trained attractors. We remark that even though these bifurcation diagrams are produced for smaller $N$ and different random realisations of $\textbf{M}$ and $\textbf{W}_{in}$, they share similar characteristics features in terms of there being a coexistence of four fixed points at small $\rho$, and that there is a bifurcation between these fixed points and two limit cycles. The results presented in \textcite{morraflynn23_MF_fly} also show that by constructing the weights of $\textbf{M}$ based on a region of a fruit-fly's connectome (where $N=426$) and repeating the same experiment, there is a slightly different bifurcation structure. What is interesting here, in comparison to the results presented in this paper, is that for different designs of $\textbf{M}$ and for $N=426$ as opposed to $N=1000$, multifunctionality is achieved for a much greater range of $\rho$ values. 

We comment that bifurcation analysis on ANNs similar to the closed-loop RC studied throughout this paper have so far mainly concentrated low dimensional ANNs such as the work of \textcite{doya1992bifurcations}, \textcite{beer1995DynamicsOfCTRNNs,beer2022codim2RNN}, and \textcite{ceni2020interpreting} where the results from a bifurcation analysis of a RC provides the insight needed to solve the `flip-flop problem'. By conducting a bifurcation analysis of low-dimensional RCs and studying how the dynamics of the RC changes as the dimension increases may provide further insight into the behaviour of higher dimensional RCs and establish how much multifunctionality a given RC can have.

We do not claim that the training method outlined in Sec.\,\ref{sec:trainingMFRC} is the optimal method for reconstructing a coexistence of any number of attractors. It is clear from the results presented in this paper that even for the closed-loop RC to reconstruct a coexistence of two relatively simple dynamical objects, $\pzc_{A}$ and $\pzc_{B}$, it is the relationship between them that is the more crucial aspect. By making further progress on establishing other limiting factors which impact the reconstruction capacity of a given closed-loop RC, it may be possible to extend the current record of the seven coexisting reconstructed chaotic attractors presented in \textcite{herteux2021MScthesis}.

In summary, this paper provides significant insight towards how a closed-loop RC achieves multifunctionality in a paradigmatic setting with a particular focus on exploring the dynamics of the closed-loop RC as it solves the seeing double problem.
However, 
there are still many unanswered questions concerning the reconstruction capacity of a given closed-loop RC and the limiting factors to the interesting dynamics a multifunctional RC can produce. While this paper improves the language surrounding the basic mechanisms needed for multifunctionality to occur from previous studies,
in future work we aim to extend the description of how RCs achieve attractor reconstruction using the formalism of generalised synchronisation, like in \textcite{LuHuntOtt18RC}, to the case of multifunctionality. Furthermore, we also aim to investigate the role of untrained attractors in how RCs achieve attractor reconstruction in more general cases. From an industry and applications perspective, the study of untrained attractors is an important issue for critical safety applications when, for instance, constructing digital twins of real-world systems as the operator needs to ensure that no unintentional dynamics can be triggered when noise is introduced. Moreover, the existence of untrained attractors raises questions about what behavioural assurances can be guaranteed for closed-loop RCs of this form.

\begin{acknowledgments}
This work was funded by the Irish Research Council Enterprise Partnership Scheme (Grant No. EPSPG/2017/301). We thank the reviewers of this paper for their very helpful comments. We would also like to thank Andrea Ceni, Lou Pecora, Tom Carroll, Ling-Wei Kong, Joschka Herteux, Daniel K\"{o}glmayr, Oliver Heilmann, Christoph R\"{a}th, Jacob Morra, and many fellow researchers from UCC, for their influential conversations and input when discussing the contents of this paper as it has evolved over the past few years.
\end{acknowledgments}

\appendix

\section{RC design and training parameters}\label{app:RCdesign}

In this paper, $\textbf{M}$ is constructed with an Erd\"{o}s-Renyi topology where each of the non-zero elements are then replaced with a random number between $-1$ and $1$, the matrix is subsequently scaled to a specific spectral radius, $\rho$. To elaborate, $\textbf{M}$ is designed such that each element is chosen independently to be nonzero with probability $P$ (i.e. sparsity $= P$ or degree $= N/P$) and these nonzero elements are chosen uniformly from $\left( -1, 1 \right)$. This random sparse matrix is rescaled such that the spectral radius, which is the maximum of the absolute values of its eigenvalues, is $\rho$. This topology is chosen in order to provide the RC with enough dynamical flexibility to solicit multistable dynamics in tasks requiring multifunctionality. Similarly, the input matrix, $\textbf{W}_{in}$, is designed such that each row has only one nonzero randomly assigned element, chosen uniformly from $\left( -1, 1 \right)$.

The numerical results presented throughout this paper are generated using the same $\textbf{M}$ and $\textbf{W}_{in}$. From further analysis
and results which are not shown here, we remark that there are relatively small quantitative changes to our results when using different initialisations of $\textbf{M}$ and $\textbf{W}_{in}$ however the main characteristics of our results remain similar. 

\begin{table}
    \centering
    \begin{tabular}{c c c c c c c c c }
    \hline \hline
        $N$ & $P$ & $\rho$ & $\sigma$ & $\gamma$ & $\beta$ & $t_{\text{listen}}$ & $t_{\text{train}}$ \\
        \hline
        $1000$ & $0.04$ & $\left[ 0.1, 2.5 \right]$ & $0.2$ & $5$ & $10^{-2}$ & $200$ & $400$ \\
        \hline \hline
    \end{tabular}
    \caption{RC design and training parameters used to generate the results displayed in Sec.\ref{sec:SD_results}.}
    \label{tab:RC_SDparams}
\end{table}

\section{Symmetry and Seeing Double}\label{sssec:SymmUAs}

It is important to note that in Figs.\,\ref{fig:BoA_rho_125} and \ref{fig:BoAs} there are certain symmetries present in the projected basins of attraction. 

In Figs.\,\ref{fig:BoA_rho_125},\,\ref{fig:BoA_rho_09}, and \ref{fig:BoA_rho_08}, the closed-loop RC approaches the same attractor whether it is initialised with its representation of either the point $\left( x', y' \right)$ or $\left( -x', -y' \right)$. However, there is a slightly different situation in Figs.\,\ref{fig:BoA_rho_04}-\ref{fig:BoA_rho_03}. Using the terminology developed by \textcite{herteux2020Symm}, when the closed-loop RC is initialised with its representation of the point $\left( x', y' \right)$ then $\hat{\boldsymbol{r}}(t)$ settles down to a particular fixed point, for example FP$_{1}$ or FP$_{3}$, if instead the closed-loop RC was initialised with its representation of $\left( -x', -y' \right)$ then $\hat{\boldsymbol{r}}(t)$ settles down to the `mirror-attractor' of that fixed point, which in this case is FP$_{2}$ or FP$_{4}$.

Despite breaking the symmetry in this closed-loop RC setup with the square terms in Eq.\,\eqref{eq:q_square} in the design of $\textbf{W}_{out}$, there are certain properties of the seeing double training data that allow for the emergence of mirror-attractors in a similar way to some of the results presented in \textcite{flynn2021symmetry} regarding the `symmetry kills the square' phenomenon. The manner in which these scenarios arise in the case of the seeing double problem are outlined in the remainder of this Appendix.

\subsection{Seeing Double and \texorpdfstring{$\textbf{W}^{(2)}_{out} \to \textbf{0}$}{TEXT}}\label{ssec:SD_symm}

In order for mirror-attractors to exist, the symmetry breaking terms in Eq.\,\eqref{eq:PredRes}, i.e., the $\boldsymbol{r}^{2}(t)$ term in Eq.\,\eqref{eq:q_square} and the corresponding $\textbf{W}_{out}^{(2)}$ term in Eq.\,\eqref{eq:BreakSym_r2}, cannot exist\cite{herteux2020Symm}.


In the case of the seeing double problem, there are two conditions which contribute to forcing $\textbf{W}^{(2)}_{out} \to \textbf{0}$ that need to be satisfied simultaneously. These are the anti-symmetric nature of the individual circular orbits, as described by Eqs.\,(28) and (29) in \textcite{flynn2021symmetry}, and the subsequent anti-symmetric nature of the combined training data set described in Eqs.\,(16) and (17) in \textcite{flynn2021symmetry} which determine $\textbf{W}_{out}$.

The first of the conditions identified above, which relates to Eqs.\,(28) and (29) in \textcite{flynn2021symmetry}, holds when the RC is trained to reconstruct only $\pzc_{A}$ or $\pzc_{B}$ for $x_{cen}=0$ as,
\begin{align}
    \hspace{-.2cm}\boldsymbol{u}_{\left(\pzc_{A}\right)}(t) \hspace{-.045cm} = \hspace{-.045cm} - \boldsymbol{u}_{\left(\pzc_{A}\right)}(t + \frac{T}{2}) \,\,\text{$\&$} \,\, \boldsymbol{u}_{\left(\pzc_{B}\right)}(t) \hspace{-.045cm} = \hspace{-.045cm} - \boldsymbol{u}_{\left(\pzc_{B}\right)}(t + \frac{T}{2}), \label{eq:uAB_Tsym}
\end{align}
is satisfied $\forall \, t$ so long as the open-loop RC responds to this driving input such that,
\begin{align}
    \hspace{-.2cm}\boldsymbol{r}_{\left(\pzc_{A}\right)}(t) \hspace{-.05cm} = \hspace{-.05cm} - \boldsymbol{r}_{\left(\pzc_{A}\right)}(t + \frac{T}{2}) \,\, \text{$\&$} \,\, \boldsymbol{r}_{\left(\pzc_{B}\right)}(t) \hspace{-.05cm} = \hspace{-.05cm} - \boldsymbol{r}_{\left(\pzc_{B}\right)}(t + \frac{T}{2}), \label{eq:rAB_Tsym}
\end{align}
where $T$ is the period of oscillation on each circular orbit. 

However, complications arise with the second condition that relates to Eqs.\,(16) and (17) in \textcite{flynn2021symmetry} and concerns the structure of the resulting $\textbf{X}_{C}$ and $\textbf{Y}_{C}$ matrices used in the case of multifunctionality. On one hand, this condition holds $\forall \, x_{cen}$ when $\pzc_{A}$ and $\pzc_{B}$ are rotating in the same direction as the RC is effectively being trained to reconstruct mirror-attractors. To be more specific, the input data, $\boldsymbol{u}(t)$, from, for instance, $\pzc_{A}$ when centered at $\left( x_{cen}, \, 0 \right)$, is equal to $-\boldsymbol{u}(t)$ when $\pzc_{A}$ is centered at $\left( -x_{cen}, \, 0 \right)$, denoted as $\pzc_{A}^{-}$. As the RC is trained to reconstruct a coexistence of $\pzc_{A}$ and $\pzc_{B}$ for a given $x_{cen}$, and $\pzc_{B} = \pzc_{A}^{-} = -\pzc_{A}$ then $\boldsymbol{u}_{\left(\pzc_{B}\right)}(t) = -\boldsymbol{u}_{\left(\pzc_{A}\right)}(t)$ and
\begin{align}
    \hspace{-0.27cm} \left( \hspace{-0.14cm} \begin{array}{c} \boldsymbol{r}_{\left(\pzc_{B}\right)}\hspace{-0.03cm}(t) \\ \boldsymbol{r}_{\left(\pzc_{B}\right)}^{2}\hspace{-0.03cm}(t) \end{array} \hspace{-0.14cm} \right) \hspace{-0.1cm} = \hspace{-0.1cm} \left( \hspace{-0.14cm} \begin{array}{c} \boldsymbol{r}_{\left(\pzc_{A}^{-}\right)}\hspace{-0.03cm}(t) \\ \boldsymbol{r}_{\left(\pzc_{A}^{-}\right)}^{2}\hspace{-0.03cm}(t) \end{array} 
\hspace{-0.14cm} \right) \hspace{-0.1cm} = \hspace{-0.1cm} \left( \hspace{-0.14cm} \begin{array}{c} - \boldsymbol{r}_{\left(\pzc_{A}\right)}\hspace{-0.03cm}(t) \\ \left(- \boldsymbol{r}_{\left(\pzc_{A}\right)}\hspace{-0.03cm}(t) \right)^{2} \end{array} \hspace{-0.14cm} \right) \hspace{-0.1cm} = \hspace{-0.1cm} \left(\hspace{-0.14cm} \begin{array}{c} - \boldsymbol{r}_{\left(\pzc_{A}\right)}\hspace{-0.03cm}(t) \\ \boldsymbol{r}_{\left(\pzc_{A}\right)}^{2}\hspace{-0.03cm}(t) \end{array} \hspace{-0.14cm} \right) \hspace{-0.1cm}.\label{eq:rpmCA}
\end{align}
Therefore, the resulting $\textbf{X}_{C}$ and $\textbf{Y}_{C}$ training data matrices are of the form described in Eqs.\,(16) and (17) in \textcite{flynn2021symmetry} and as a consequence, $\textbf{W}_{out}^{(2)} = \textbf{0}$.

On the other hand, when $\pzc_{A}$ and $\pzc_{B}$ are rotating in opposite directions there is no mirror symmetry present in the combined training data sets, it is only when $x_{cen}=0$ that `symmetry kills the square' as the combined training data matrices from both $\pzc_{A}$ and $\pzc_{B}$ satisfy necessary conditions. As a further comment, if instead the RC was trained to reconstruct a coexistence of $\pzc_{A}$ and $\pzc_{B}$ when these are of different radii and $x_{cen}=0$, so long as Eqs.\,(28) and (29) in \textcite{flynn2021symmetry} are satisfied for both $\boldsymbol{u}_{\left(\pzc_{A}\right)}(t)$ and $\boldsymbol{u}_{\left(\pzc_{B}\right)}(t)$ then $\textbf{W}_{out}^{(2)} = \textbf{0}$ as the RC is trained to provide a coexistence of attractors which are all anti-symmetric to themselves. In return, this places no constraint on the behaviour of the corresponding response attractors, $\pzs_{A}$ and $\pzs_{B}$, in $\mathbb{S}$ as these do not necessarily need to be mirror-attractors of one another as already shown in Fig.\,\ref{fig:BoA_rho_125}.

\subsection{\texorpdfstring{$\rho$ vs. $\textbf{W}_{out}^{(2)}$}{TEXT}}

\begin{figure*}
    \centering
    \begin{subfigure}{0.33\textwidth}
        \centering
        \includegraphics[width=\textwidth]{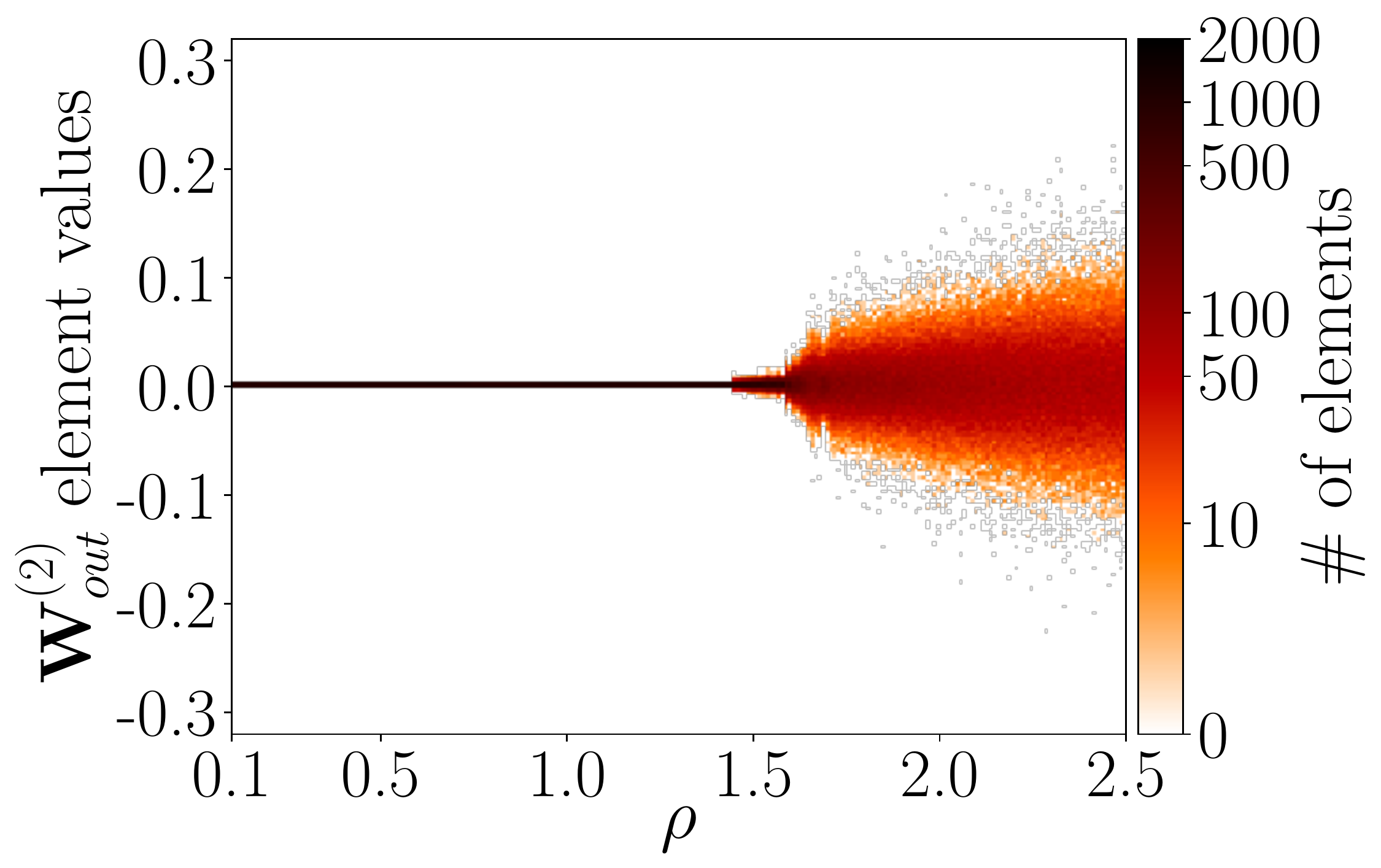}
        \caption{$x_{cen}=0$}
        \label{fig:SDrho_wout2xcen0}
    \end{subfigure}
        \begin{subfigure}{0.33\textwidth}
        \centering
        \includegraphics[width=\textwidth]{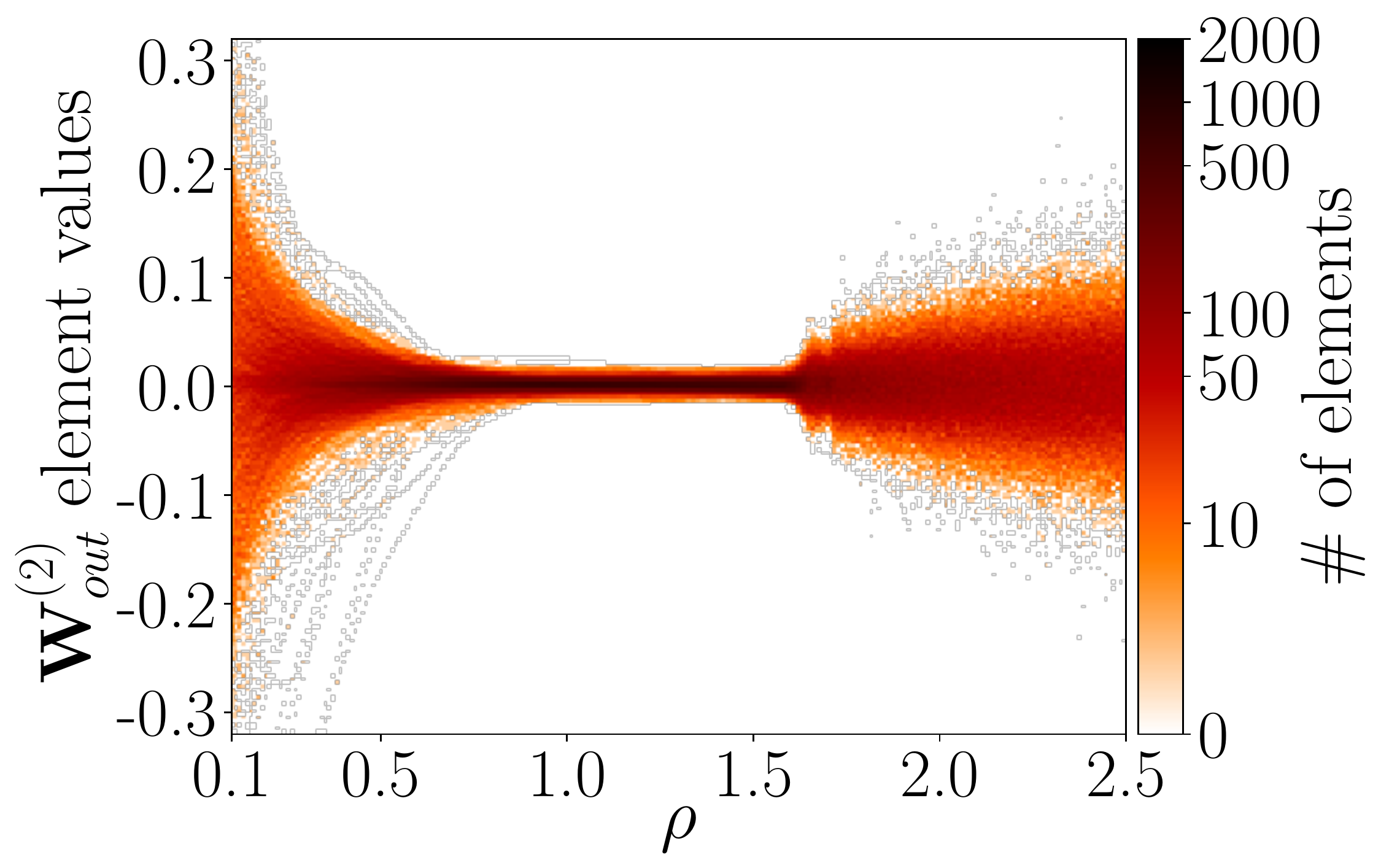}
        \caption{$x_{cen}=-1$}
        \label{fig:SDrho_wout2xcen1}
    \end{subfigure}
    \begin{subfigure}{0.33\textwidth}
        \centering
        \includegraphics[width=\textwidth]{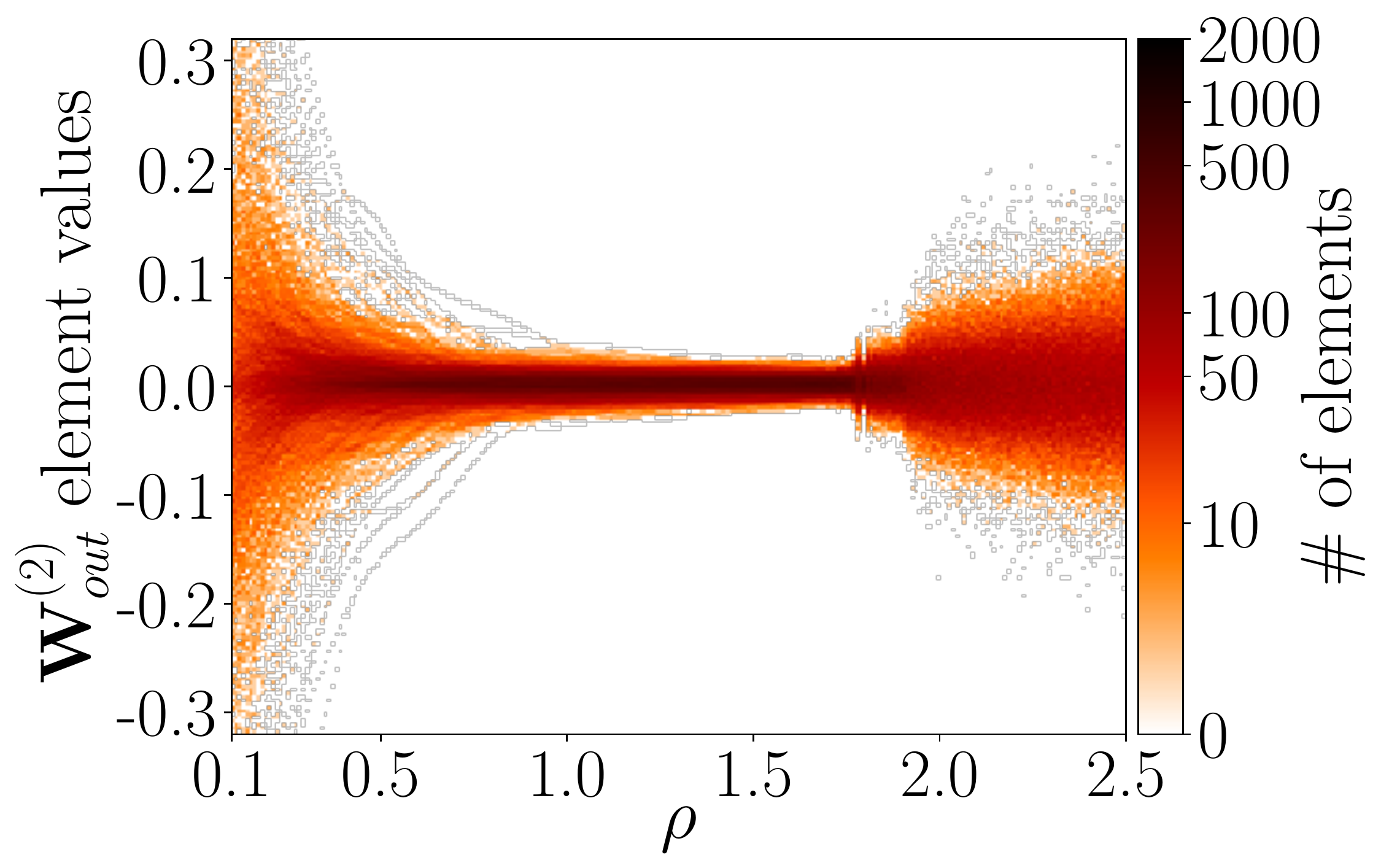}
        \caption{$x_{cen}=-5$}
        \label{fig:SDrho_wout2xcen5}
    \end{subfigure}
    \caption{Histograms of the $\textbf{W}^{(2)}_{out}$ elements for $\rho \in \left[ 0.1, 2.5 \right]$ when in (a) $x_{cen}=0$, in (b) $x_{cen}=-1$, and in (c) $x_{cen}=-5$.}
    \label{fig:SDrho_wout2}
\end{figure*}

Given the insight gained about the role of $\rho$ in the behaviour of $\textbf{W}^{(2)}_{out}$ from Fig.\,(7) in \textcite{flynn2021symmetry}, the same procedure is applied here for the present case of the seeing double problem when $x_{cen} = 0$ and the result of this is illustrated in Fig.\,\ref{fig:SDrho_wout2xcen0}. There is a threshold value of $\rho$ associated with the symmetry-breaking bifurcation discussed in Sec.\,\ref{sssec:SD_Floq_anal} which corresponds with $\textbf{W}^{(2)}_{out}$ coming into existence in Fig.\,\ref{fig:SDrho_wout2xcen0}. Since $\textbf{W}^{(2)}_{out} = \textbf{0}$ for all $\rho$ values less than this critical value, each solution of the closed-loop RC is required to either have a mirror-attractor or is anti-symmetric to itself, whether or not multifunctionality is achieved. While the behaviour of $\textbf{W}_{out}^{(2)}$ in Fig.\,\ref{fig:SDrho_wout2xcen0} is qualitatively similar to Fig.\,(7) in \textcite{flynn2021symmetry} there are some differences, mainly that $\textbf{W}_{out}^{(2)}$ reappears at smaller $\rho$ and its elements grow much faster and larger in magnitude over a similar interval of $\rho$. 

Nevertheless, the existence of mirror-attractors in Figs.\,\ref{fig:BoA_rho_04} and \ref{fig:BoA_rho_03}, which led towards this discussion in the first place, and that $\textbf{W}_{out}^{(2)}=\textbf{0}$ for the corresponding $\rho$ values shown in Fig.\,\ref{fig:SDrho_wout2xcen0}, together show that while both driving inputs possess the symmetry described in Eq.\,\eqref{eq:uAB_Tsym} and the open-loop RC may respond to these driving inputs according to Eq.\,\eqref{eq:rAB_Tsym}, it is still possible that $\textbf{W}_{out}^{(2)}=\textbf{0}$ even if the resulting attractors in $\mathbb{P}$ are not periodic but are instead mirrored fixed points.


When $x_{cen} \neq 0$ the symmetry in the training data is broken, which is much like shifting the location of the mirror-attractors as in \textcite{flynn2021symmetry}. The effect this has on the elements of $\textbf{W}_{out}^{(2)}$ for $\rho \in \left[ 0.1, 2.5 \right]$ is shown in Figs.\,\ref{fig:SDrho_wout2xcen1}-\ref{fig:SDrho_wout2xcen5}. Similar to the results depicted in Fig.\,2 in \textcite{flynn2021symmetry}, the central message of Figs.\,\ref{fig:SDrho_wout2xcen1}-\ref{fig:SDrho_wout2xcen5} is that as the magnitude of $x_{cen}$ increases so too does the magnitude of $\textbf{W}_{out}^{(2)}$'s elements.


\subsection{\texorpdfstring{$\left( \textbf{W}_{out}^{(1)} \right)_{x_{cen}} = \left( \textbf{W}_{out}^{(1)} \right)_{-x_{cen}}$ and $\left( \textbf{W}_{out}^{(2)} \right)_{x_{cen}} = - \left( \textbf{W}_{out}^{(2)} \right)_{-x_{cen}}$}{TEXT}}\label{sssec:W1W2pmxcen}

While it is shown in Fig.\,\ref{fig:SDrho_wout2xcen0} that $\textbf{W}_{out}^{(2)} = \textbf{0}$ for $x_{cen}=0$, there is another `symmetry kills the square' phenomenon and related to the structure of the seeing double training data that is generated as follows.

Building on the statements made earlier which led towards the relationship established in Eq.\,\eqref{eq:rpmCA}, $\boldsymbol{u}_{\left(\pzc_{A}\right)}(t) = -\boldsymbol{u}_{\left(\pzc_{A}^{-}\right)}(t)$ and the resulting input training data matrix, $\textbf{Y}$, for a given $x_{cen}$ and $-x_{cen}$ are related in the following way,
\begin{align}
    \textbf{Y}_{x_{cen}} = - \textbf{Y}_{-x_{cen}}.
\end{align}
Therefore, in order for the training to be successful and $\hat{\pzc_{A}} = - \hat{\pzc_{A}}^{-}$, 
then the following is also required,
\begin{align}
    \left( \textbf{W}_{out} \right)_{x_{cen}} \left( \hspace{-.12cm} \begin{array}{c} \boldsymbol{r}_{\left(\pzc_{A}\right)}(t) \\ \boldsymbol{r}_{\left(\pzc_{A}\right)}^{2}(t) \end{array} \hspace{-.12cm} \right) = - \left( \textbf{W}_{out} \right)_{-x_{cen}} \left( \hspace{-.12cm} \begin{array}{c} \boldsymbol{r}_{\left(\pzc_{A}^{-}\right)}(t) \\ \boldsymbol{r}_{\left(\pzc_{A}^{-}\right)}^{2}(t) \end{array} \hspace{-.12cm} \right),\label{eq:Wout_pmxcen}
\end{align}
where $\left( \textbf{W}_{out} \right)_{x_{cen}}$ is the resulting $\textbf{W}_{out}$ matrix trained using data from $\pzc_{A}$ and $\pzc_{B}^{-}$ when centered at $\left( x_{cen}, \, 0 \right)$ and $\left( -x_{cen}, \, 0 \right)$ respectively for a given $x_{cen}$, and similarly for $\left( \textbf{W}_{out} \right)_{-x_{cen}}$ with $\pzc_{A}^{-}$ and $\pzc_{B}$ for the corresponding $-x_{cen}$.

Using the information from Eq.\,\eqref{eq:rpmCA}, this allows Eq.\,\eqref{eq:Wout_pmxcen} to be rewritten as,
\begin{align}
    \hspace{-.25cm}\left( \textbf{W}_{out} \right)_{x_{cen}} \left( \hspace{-.12cm} \begin{array}{c} - \boldsymbol{r}_{\left(\pzc_{A}^{-}\right)}(t) \\ \boldsymbol{r}_{\left(\pzc_{A}^{-}\right)}^{2}(t) \end{array} \hspace{-.12cm} \right) \hspace{-.1cm}= - \left( \textbf{W}_{out} \right)_{-x_{cen}} \left( \hspace{-.12cm} \begin{array}{c} \boldsymbol{r}_{\left(\pzc_{A}^{-}\right)}(t) \\ \boldsymbol{r}_{\left(\pzc_{A}^{-}\right)}^{2}(t) \end{array} \hspace{-.12cm} \right),
\end{align}
from which it follows that,
\begin{align}
    &- \left( \textbf{W}_{out}^{(1)} \right)_{x_{cen}} \boldsymbol{r}_{\left(\pzc_{A}^{-}\right)}(t) + \left( \textbf{W}_{out}^{(2)} \right)_{x_{cen}} \boldsymbol{r}_{\left(\pzc_{A}^{-}\right)}^{2}(t) \nonumber \\
    & = - \left( \textbf{W}_{out}^{(1)} \right)_{-x_{cen}}  \boldsymbol{r}_{\left(\pzc_{A}^{-}\right)}(t) - \left( \textbf{W}_{out}^{(2)} \right)_{-x_{cen}} \boldsymbol{r}_{\left(\pzc_{A}^{-}\right)}^{2}(t),\\
    &\implies \left( \left( \textbf{W}_{out}^{(2)} \right)_{x_{cen}} + \left( \textbf{W}_{out}^{(2)} \right)_{-x_{cen}} \right) \boldsymbol{r}_{\left(\pzc_{A}^{-}\right)}^{2}(t) \nonumber \\
    &\quad \quad = \left( \left( \textbf{W}_{out}^{(1)} \right)_{x_{cen}} - \left( \textbf{W}_{out}^{(1)} \right)_{-x_{cen}} \right) \boldsymbol{r}_{\left(\pzc_{A}^{-}\right)}(t),\\
    &\therefore \left( \textbf{W}_{out}^{(1)} \right)_{x_{cen}} = \left( \textbf{W}_{out}^{(1)} \right)_{-x_{cen}}, \nonumber \\
    & \quad \quad \text{and} \,\, \left( \textbf{W}_{out}^{(2)} \right)_{x_{cen}} = - \left( \textbf{W}_{out}^{(2)} \right)_{-x_{cen}}.\label{eq:W1W2pmxcen}
\end{align}

\begin{figure*}
    \begin{subfigure}{0.245\textwidth}
        \centering
        \includegraphics[width=\textwidth]{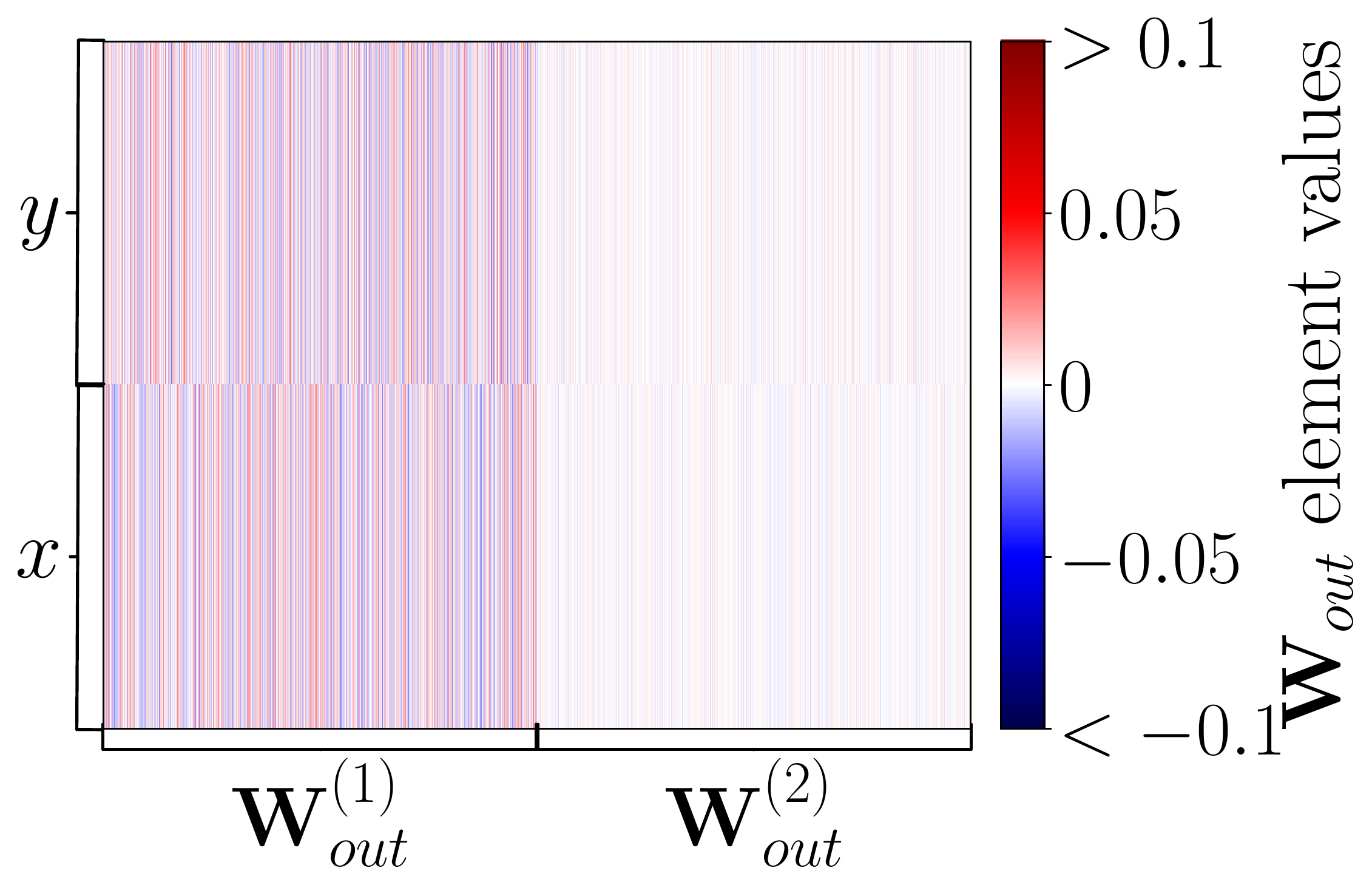}
        \caption{$\left( \textbf{W}_{out} \right)_{x_{cen}=3}$}
        \label{fig:SDWout_xcen3}
    \end{subfigure}
        \begin{subfigure}{0.245\textwidth}
        \centering
        \includegraphics[width=\textwidth]{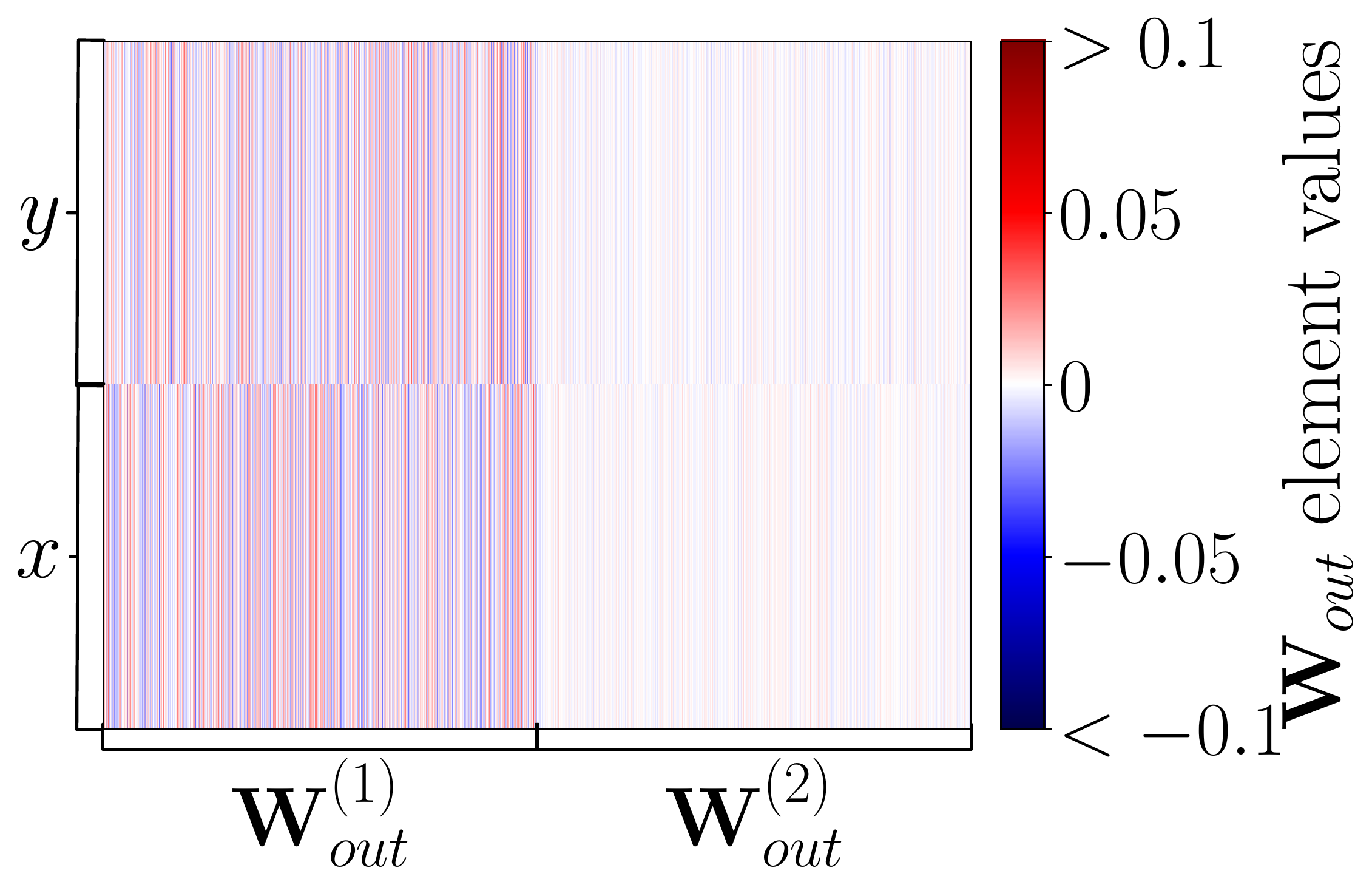}
        \caption{$\left( \textbf{W}_{out} \right)_{x_{cen}=-3}$}
        \label{fig:SDWout_xcenm3}
    \end{subfigure}
    \begin{subfigure}{0.245\textwidth}
        \centering
        \includegraphics[width=\textwidth]{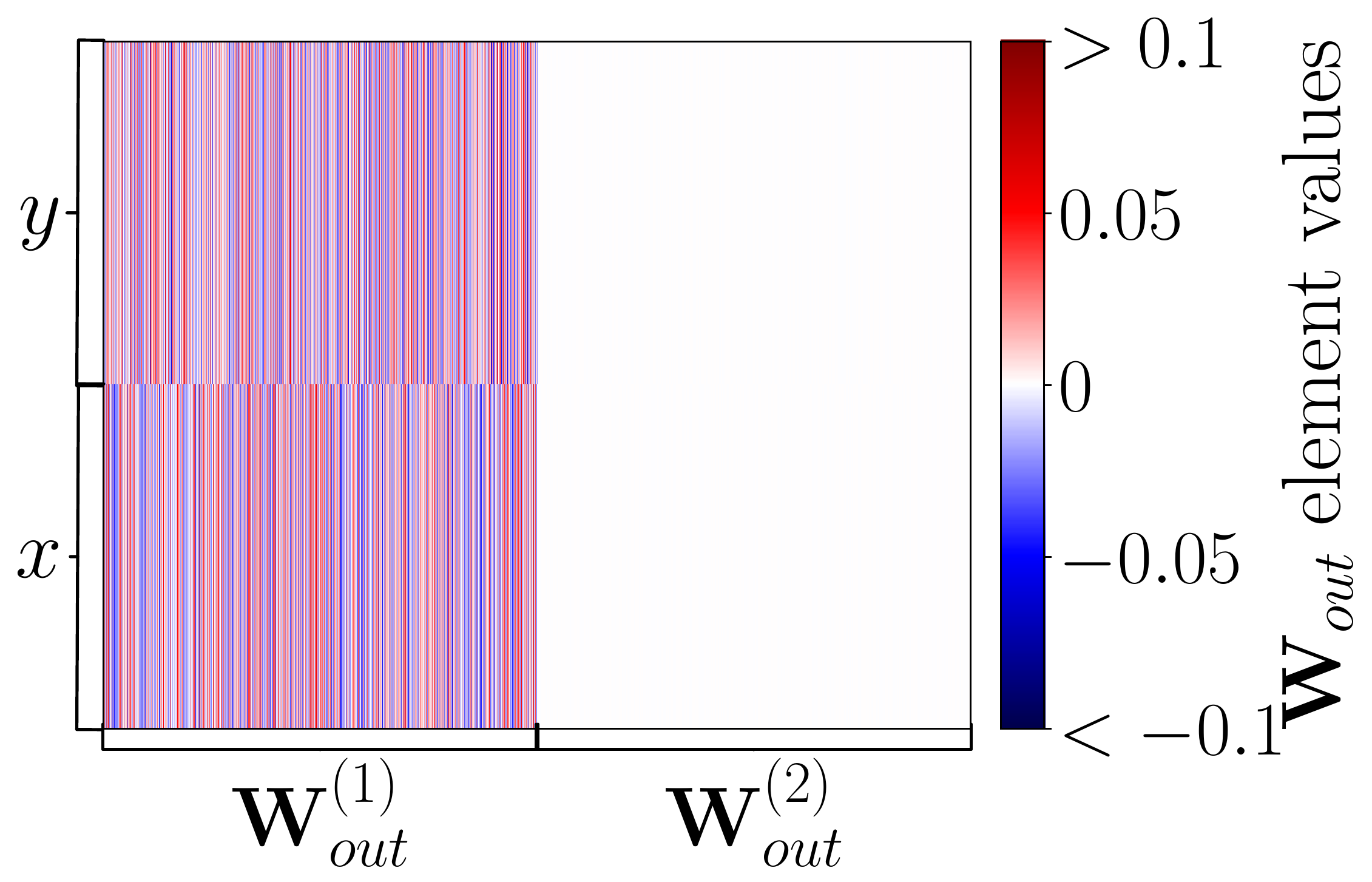}
        \caption{$\left( \textbf{W}_{out} \right)_{x_{cen}=3} + \left( \textbf{W}_{out} \right)_{x_{cen}=-3}$}
        \label{fig:SDWout_xcen3plus_m3}
    \end{subfigure}
    \begin{subfigure}{0.245\textwidth}
        \centering
        \includegraphics[width=\textwidth]{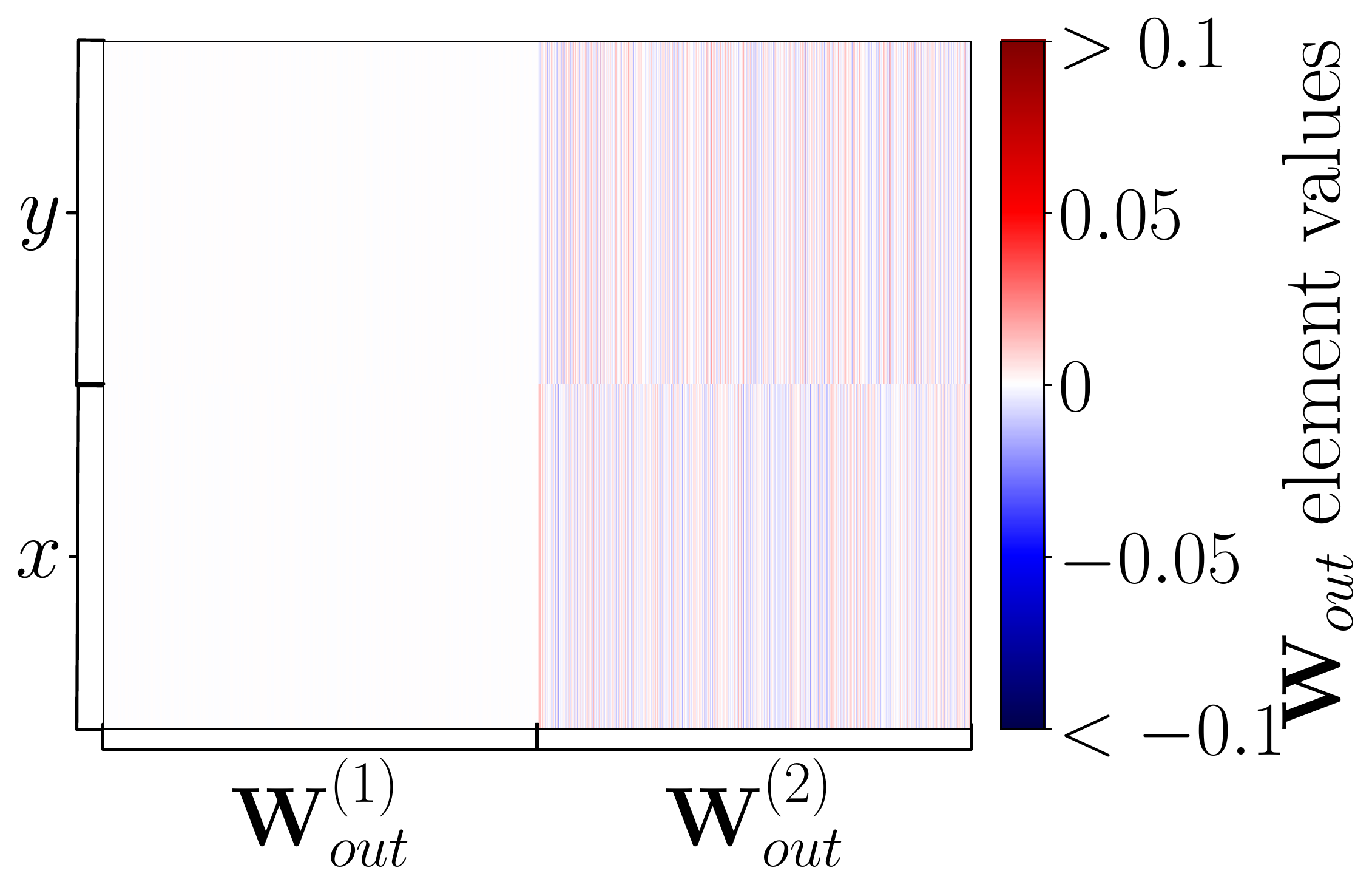}
        \caption{$\left( \textbf{W}_{out} \right)_{x_{cen}=3} - \left( \textbf{W}_{out} \right)_{x_{cen}=-3}$}
        \label{fig:SDWout_xcen3minus_m3}
    \end{subfigure}
    \caption{$\left( \textbf{W}_{out} \right)_{x_{cen}}$ elements with $\rho = 1.45$ for $x_{cen}=3$ in (a) and $x_{cen}=-3$ in (b). Result of element-wise addition of $\left( \textbf{W}_{out} \right)_{x_{cen}=3} + \left( \textbf{W}_{out} \right)_{x_{cen}=-3}$ in (c) and element-wise subtraction of $\left( \textbf{W}_{out} \right)_{x_{cen}=3} - \left( \textbf{W}_{out} \right)_{x_{cen}=-3}$ in (d). Each column represents the weight, $w_{i,j}$, given to the $j^{th}$ component of $\boldsymbol{q}(\boldsymbol{r}(t)) = ( \boldsymbol{r}(t),\boldsymbol{r}^{2}(t) )^{T}$. The distinction is made between which columns belong to $\textbf{W}_{out}^{(1)}$ and $\textbf{W}_{out}^{(2)}$. Note, the matrix plots may vary depending on the PDF viewer used, it is recommended to use `Adobe Reader' or `Chrome PDF Viewer' to see these images as intended.}
    \label{fig:SDWout_xcen_pm}
\end{figure*}

The above analytical results in Eq.\,\eqref{eq:W1W2pmxcen} are confirmed numerically in Fig.\,\ref{fig:SDWout_xcen_pm} by examining the behaviour of the resulting $\left( \textbf{W}_{out} \right)_{x_{cen}}$ matrices after training the RC with data from $\pzc_{A}$ and $\pzc_{B}$ for $x_{cen}= 3$ in Fig.\,\ref{fig:SDWout_xcen3} and $x_{cen}=-3$ in Fig.\,\ref{fig:SDWout_xcenm3} where $\rho = 1.45$ in both cases. The elements of the resulting matrix from the element-wise addition of $\left( \textbf{W}_{out} \right)_{x_{cen}=3}$ and $\left( \textbf{W}_{out} \right)_{x_{cen}=-3}$ is depicted in Fig.\,\ref{fig:SDWout_xcen3plus_m3} which shows that $\left( \textbf{W}_{out}^{(2)} \right)_{x_{cen}} + \left( \textbf{W}_{out}^{(2)} \right)_{-x_{cen}} = \textbf{0}$. Similarly, in Fig.\,\ref{fig:SDWout_xcen3minus_m3}, the element-wise subtraction of $\left( \textbf{W}_{out} \right)_{x_{cen}=-3}$ from $\left( \textbf{W}_{out} \right)_{x_{cen}=3}$ results in $\left( \textbf{W}_{out}^{(1)} \right)_{x_{cen}} - \left( \textbf{W}_{out}^{(1)} \right)_{-x_{cen}} = \textbf{0}$.

While these results are interesting from an analytical and numerical perspective, their influence has been a common factor throughout many of the results presented in this paper and arise due to the constraint in Eq.\,\eqref{eq:W1W2pmxcen} as imposed by the structure of the training data.

\section*{DATA AVAILABILITY}
The data that support the findings of this study are available from the corresponding author upon reasonable request.



%

\end{document}